\documentclass{amsart}

\usepackage{amsmath,amssymb,amsfonts}
\usepackage{amsthm}
\usepackage{pgfplots}
\pgfplotsset{compat=1.18}

\usepgfplotslibrary{groupplots,patchplots}

 \usepackage{enumitem}

\usepackage{tikz}
 \usepackage{graphicx}

\usetikzlibrary{patterns, arrows.meta,calc,decorations.pathreplacing}
\usepackage[hidelinks]{hyperref}
\usepackage{amsaddr}

\newtheorem{theorem}{Theorem}[section]
\newtheorem{lemma}[theorem]{Lemma}
\newtheorem{corollary}[theorem]{Corollary}
\newtheorem{proposition}[theorem]{Proposition}

\theoremstyle{definition}
\newtheorem{definition}[theorem]{Definition}
\newtheorem{question}[theorem]{Question}

\newtheorem{example}[theorem]{Example}

\theoremstyle{remark}
\newtheorem{remark}[theorem]{Remark}

\numberwithin{equation}{section}

\providecommand{\abs}[1]{\lvert#1\rvert}
\providecommand{\norm}[1]{\lVert#1\rVert}
\hypersetup{
  pdfauthor={Krzysztof J. Ciosmak},
  pdftitle={Localisation for the second-order Beckmann problem: bimartingale couplings and leaf decompositions}
}

\definecolor{orcidgreen}{cmyk}{.407,.004,.999,0}
\DeclareRobustCommand{\orcidlogo}{%
  \begin{tikzpicture}[x=1bp,y=1bp,baseline=0pt]
    \path[use as bounding box] (0,0) rectangle (9,9);
    \begin{scope}
      \clip (.06,0) rectangle (8.504,8.504);
      \fill[orcidgreen] (4.252,4.252) circle[radius=4.252bp];
    \end{scope}
    \fill[white] (2.356,2.319) rectangle (2.868,5.877);
    \fill[white] (2.612,6.617) circle[radius=.335bp];
    \fill[white,even odd rule] (3.617,5.877) -- (4.999,5.877)
      .. controls (6.314,5.877) and (6.893,4.937) .. (6.893,4.097)
      .. controls (6.893,3.185) and (6.178,2.318) .. (5.005,2.318)
      -- (3.617,2.318) -- cycle
      (4.129,2.777) -- (4.942,2.777)
      .. controls (6.1,2.777) and (6.365,3.659) .. (6.365,4.097)
      .. controls (6.365,4.813) and (5.913,5.417) .. (4.917,5.417)
      -- (4.129,5.417) -- cycle;
  \end{tikzpicture}}
\makeatletter
\def\@setemails{\mbox{\emailaddrname:\space}{\ttfamily\emails}.}
\makeatother

\begin{document}

\title[Localisation for the second-order Beckmann problem]{Localisation for the second-order Beckmann problem: bimartingale couplings and leaf decompositions}

\author[Localisation for the second-order Beckmann problem]{Krzysztof J. Ciosmak \texorpdfstring{\href{https://orcid.org/0000-0001-9571-1160}{\orcidlogo}}{}}

\address{Beijing Institute of Mathematical Sciences and Applications,\\
 No.~544, Hefangkou Village, Huaibei Town, Huairou District, Beijing 101408, China}
\address{Institute of Mathematics of the Polish Academy of Sciences,\\
 ul.~\'Sniadeckich 8, 00-656 Warszawa, Poland}
\email{kciosmak@bimsa.cn, kciosmak@impan.pl}
\thanks{
 \textbf{Acknowledgments:} the author would like to express his thanks to Karol Bo\l{}botowski for a discussion that introduced the author to the details of the Kantorovich--Rubinstein duality theory for the Hessian. A large part of this work was done while the author was a postdoctoral fellow at the University of Toronto, whose support is gratefully acknowledged.
 }

\keywords{second-order localisation; second-order Beckmann problem; bimartingale couplings; convex-concave order; leaf decompositions; optimal transport}

\subjclass[2020]{Primary: 28A50, 49K20, 49Q22, 60G42, 60G48; secondary: 52A20, 74P05, 90C25}

\date{17.11.2026}

 \begin{abstract}

 We develop a second-order localisation theory for optimal transport based on
leaf decompositions and bimartingale couplings, which impose
martingale constraints in opposite directions on complementary orthogonal
subspaces. It provides a second-order
analogue of the classical decomposition of optimal transport into transport
rays and monotone couplings: transport rays are replaced by the leaves of the
$1$-Lipschitz derivative map $Du$ of an optimal dual potential
$u\in C^{1,1}(\mathbb R^n)$, while monotone couplings are replaced by
bimartingale couplings.

We apply this framework to the three-marginal optimal transport problem
introduced by Bo\l{}botowski and Bouchitt\'e, whose relaxation is the second-order Beckmann problem.

We characterise the existence of bimartingale couplings through
a convex-concave order condition. This yields a generalisation of Strassen's
theorem from convex order to the convex-concave setting. Equivalently, the
dual problem associated with the second-order Beckmann problem admits an
optimiser whose derivative is an isometry.

For absolutely continuous measures with common barycentre, assuming the
existence of an optimal plan with absolutely continuous third marginal, we
prove that every optimal plan decomposes into a family of problems on the leaves of $Du$. On each leaf, all optimal plans are completely described by bimartingale couplings between the corresponding conditional measures.

Without this absolute continuity assumption on the third marginal, we
show that the leaf
decomposition persists in a more general form: optimal plans are mixtures of
plans concentrated on triples $(x,y,z)$ satisfying
$x\in \mathcal{S}_1$, $y\in \mathcal{S}_2$, and $z\in \mathcal{S}_1\cap \mathcal{S}_2$, where $\mathcal{S}_1$ and $\mathcal{S}_2$ are
neighbouring leaves of $Du$.

 \end{abstract}

\renewcommand{\leftmark}{\MakeUppercase{\shorttitle}}
\maketitle

\section{Introduction}

The objective of this paper is to establish a second-order localisation
principle for optimal transport and to apply it to the Hessian-constrained
transport problem introduced by Bo\l{}botowski and Bouchitt\'e
\cite{Bolbotowski2024}. The main result gives a geometric description of
optimal three-marginal transport plans through the leaves of an optimal dual
potential and identifies the corresponding conditional transport problems.

Our results may be viewed as a second-order analogue of the classical
Sudakov localisation principle. In the Kantorovich--Rubinstein theory,
optimal transports are localised on the transport rays of an optimal
$1$-Lipschitz potential. The residual one-dimensional transport problems are
then described by monotone couplings. In the present setting, transport rays
are replaced by the leaves of the $1$-Lipschitz derivative $Du$ of an optimal
$\mathcal C^{1,1}(\mathbb{R}^n)$ dual potential, while monotone couplings are replaced by
bimartingale couplings. Thus the geometry of the derivative of the dual
potential determines the localisation structure of the second-order problem.

The paper of Bo\l{}botowski and Bouchitt\'e \cite{Bolbotowski2024}
established an analogue of the Kantorovich--Rubinstein duality theory for the Hessian and the
existence of optimal solutions for the corresponding three-marginal problem.
The aim of the present work is different: we provide the first structural
description of all optimal plans. In particular, we introduce bimartingale
couplings, prove the characterisation of their existence through a convex-concave order, yielding a convex-concave analogue of Strassen's theorem, and show that they describe the conditional optimisers arising after localisation.

The proof has three main ingredients. First, we analyse the geometry of
leaves of $1$-Lipschitz derivatives. In particular, we prove that the
intersection of two leaves is a common convex face and establish the precise
compatibility conditions for neighbouring leaves. Second, we solve completely
the model case in which the derivative of the dual optimiser is an isometry.
In this case the optimal solutions to the problem are governed by bimartingale couplings. Third, combining these
ingredients with measurable disintegration gives the localisation theorem for general optimal plans, under the assumption that there exists an optimal plan with absolutely continuous marginals.

Without the absolute continuity assumption on the
third marginal, the same localisation principle remains valid in a more
general form: optimal plans may couple neighbouring leaves, but the third
point is forced to lie in their common face. Thus the additional
phenomenon is the interaction of adjacent leaves along their interfaces.

The decomposition of optimal transport into transport rays was first
proposed by Sudakov \cite{Sudakov1976}. Sudakov envisioned solving the Monge
problem by gluing optimal solutions on rays, an approach which was later
completed through the works of Caffarelli, Feldman and McCann
\cite{Caffarelli2002}, Trudinger and Wang \cite{Trudinger2001}, Ambrosio and
Pratelli \cite{Pratelli2001}, and Evans and Gangbo \cite{Evans1999}.
The localisation developed here should be understood as a higher-order
counterpart of this theory.

We refer the reader to \cite{Santambrogio2015}, \cite{Villani2003}, \cite{Villani2009}
for comprehensive presentations of optimal transport and, in particular, to
\cite[Chapter 18]{Maggi2023} for an exposition of the Sudakov's strategy for the Monge problem.

The second-order Beckmann problem is a relaxation of the three-marginal
problem considered here and, in dimension two, is related to the optimal
design of grillages supporting prescribed loads. The dual problem involves
the maximisation of integrals against the difference of two measures over the
class of functions with $1$-Lipschitz derivative. This gives rise to a metric
on probability measures with common barycentre, which is a particular case of
the Zolotarev ideal metrics; see \cite{Zolotarev1984}, \cite{Hanin1994},
and \cite{Rachev1998}.

\subsection{Classical Kantorovich--Rubinstein duality and related decomposition into transport rays}

In the optimal transport theory the Kantorovich--Rubinstein duality is a cornerstone of many developments. It tells that if $\mu,\nu\in\mathcal{P}_1(\mathbb{R}^n)$, that is, if $\mu,\nu$ are two Borel probability measures with finite first moments, then the Wasserstein distance between $\mu$ and $\nu$
\begin{equation}\label{eqn:lone}
    \inf\Big\{\int_{\mathbb{R}^n\times\mathbb{R}^n}\norm{x-y}\, d\pi(x,y)\mid \pi\in\Gamma(\mu,\nu)\Big\}
\end{equation}
is equal to
\begin{equation}\label{eqn:suplip}
    \sup\Big\{\int_{\mathbb{R}^n}v\, d(\nu-\mu)\mid v\colon\mathbb{R}^n\to\mathbb{R}\text{ is }1\text{-Lipschitz}\Big\}.
\end{equation}
Here $\norm{\cdot}$ denotes the Euclidean norm, and $\Gamma(\mu,\nu)$ stands for the set of Borel probability measures $\pi\in\mathcal{P}(\mathbb{R}^n\times\mathbb{R}^n)$ such that $\mathrm{P}_1\pi=\mu,\mathrm{P}_2\pi=\nu $, where $\mathrm{P}_i\pi$ for $i=1,2$ stand for the respective marginal measures of $\pi$, i.e., measures pushed forward by the $i$-th coordinate projections $p_i$. Any minimiser of (\ref{eqn:lone}) we shall call an \emph{optimal transport} for $\mu,\nu$. The existence of such minimisers follows by the Prokhorov theorem, see \cite{Villani2003}.
Any $1$-Lipschitz function $u\colon\mathbb{R}^n\to\mathbb{R}$ that attains the supremum (\ref{eqn:suplip}) we shall call an \emph{optimal potential} for $\mu,\nu$. Again, the existence is standard -- see, e.g.,  \cite[Theorem 2, p. 11]{Ciosmak20211}.

\subsubsection{Optimality conditions}

The duality of (\ref{eqn:lone}) and of (\ref{eqn:suplip}) shows that $\pi\in\Gamma(\mu,\nu)$ and a $1$-Lipschitz $u$ are optimal transport and optimal potential respectively if and only if
\begin{equation}\label{eqn:optimalityone}
    u(y)-u(x)=\norm{x-y}\text{ for }\pi\text{-almost every }(x,y)\in\mathbb{R}^n\times\mathbb{R}^n.
\end{equation}

\subsubsection{Transport rays}

A key concept that allows to characterise and analyse the optimal solutions of (\ref{eqn:lone}) is the notion of \emph{transport rays}. A transport ray of $u$ is a maximal set $\mathcal{T}\subset\mathbb{R}^n$ such that $u|_{\mathcal{T}}$ is an isometry. That is, for any $x,y\in\mathcal{T}$
\begin{equation*}
    \abs{u(x)-u(y)}=\norm{x-y},
\end{equation*}
and for any $z\notin\mathcal{T}$ there exists some $t\in\mathcal{T}$ such that
\begin{equation*}
    \abs{u(t)-u(z)}<\norm{t-z}.
\end{equation*}

\subsubsection{Decomposition}\label{s:transportrays}

For any $1$-Lipschitz function $u\colon\mathbb{R}^n\to\mathbb{R}$ the corresponding transport rays enjoy the important property that -- up to a set of Lebesgue measure zero -- any two distinct transport rays are disjoint, see, e.g., \cite{Ciosmak2021}.
The condition (\ref{eqn:optimalityone}) implies that $x,y$ belong to a common transport ray of an optimal $u$. Thus, any optimal transport has to occur within the transport rays of an optimal potential. There are, however, further restrictions on the optimal transports.

\subsubsection{Monotone order}\label{s:monotone}

Note that the transport rays are closed and convex, see \cite{Ciosmak2021}. Since optimal $u\colon \mathbb{R}^n\to\mathbb{R}$ is isometric  on each of the rays, they are necessarily at most one-dimensional. Let us analyse the optimal transport problem on a transport ray -- we assume that the two considered measures $\mu,\nu$ belong to $\mathcal{P}_1(\mathbb{R})$ and that the identity is an optimal potential. Thus, we may apply \cite[Proposition 4.4]{Ciosmak20212} to infer that such a pair of measures is in the \emph{non-decreasing order}. This is to say, for any bounded monotonically non-decreasing function $f\colon\mathbb{R}\to\mathbb{R}$ we have
\begin{equation}\label{eqn:increasing}
    \int_{\mathbb{R}}f\, d\mu\leq \int_{\mathbb{R}}f\, d\nu.
\end{equation}
Thus, this condition is equivalent to the existence of a coupling $\pi\in\mathcal{P}(\mathbb{R}\times\mathbb{R})$ such that
\begin{equation*}
   x\leq y\text{ for }\pi\text{-almost every }(x,y)\in\mathbb{R}^2.
\end{equation*}
Indeed, \cite{Ciosmak20212} shows that (\ref{eqn:increasing}) is equivalent to the identity
being optimal potential for $\mu,\nu$. For every optimal coupling,
(\ref{eqn:optimalityone}) then reads $y-x=\abs{y-x}$, $\pi$-almost surely.
Conversely, integration of $f(x)\leq f(y)$ against a coupling $\pi$ concentrated
on $x\leq y$ gives (\ref{eqn:increasing}). Any such transport we shall call a
non-decreasing coupling.

\subsubsection{Further partition}\label{s:fur}

We note that if for some non-decreasing $f\colon\mathbb{R}\to\mathbb{R}$ we had equality in (\ref{eqn:increasing}), then necessarily any non-decreasing coupling $\pi\in\Gamma(\mu,\nu)$ would be constrained within the sets on which $f$ is constant. This is to say, for such $f$ and such $\pi$, we would have
\begin{equation*}
    f(x)=f(y)\text{ for }\pi\text{-almost every }(x,y)\in\mathbb{R}^2.
\end{equation*}
Taking $f=\mathbf{1}_{[a,\infty)}$ for $a\in\mathbb{R}$ equality in
(\ref{eqn:increasing}) gives $\pi\{(x,y)\in\mathbb{R}^2\mid x<a\leq y\}=0$.
Such optimal transport induced partitions were studied in \cite{DiMarino2018} with an application to the entropic regularisation of the optimal transport.

\subsubsection{Characterisation of optimal transports}

Let $u$ be a $1$-Lipschitz optimal potential for absolutely continuous measures $\mu,\nu\in\mathcal{P}(\mathbb{R}^n)$. As follows e.g. by \cite{Ciosmak2021} the measures $\mu,\nu$ can be disintegrated with respect to the partition induced by the transport rays of $u$. Then a transport $\pi\in\Gamma(\mu,\nu)$ is optimal if and only if the transports on the rays are optimal. Since on each transport ray $u$ is still an optimal potential, this means that these transports have to be non-decreasing couplings. Conversely, if we find a measurable family of non-decreasing couplings of the conditional measures of $\mu$ and $\nu$ on the rays oriented by this $u$ and glue them up, we will obtain an optimal transport for $\mu$ and $\nu$.

We do not elaborate further on this observation and continue to the main results of the paper.

\subsection{Kantorovich--Rubinstein duality for the Hessian}

We shall now focus on providing an analogous description to the one above of optimal plans for the three-marginal optimal transport problem.

The analogy with the classical Kantorovich--Rubinstein theory may be
summarised as follows. In the first-order problem, optimal transports are
localised on transport rays of an optimal $1$-Lipschitz potential, and the
remaining one-dimensional problems are governed by monotone couplings. In the
second-order problem, the corresponding localisation occurs on the
leaves of the $1$-Lipschitz derivative $Du$ of an optimal dual potential
$u\in \mathcal{C}^{1,1}(\mathbb{R}^n)$. On each such leaf, optimal transports are governed by bimartingale couplings. Thus the results below may be viewed as a second-order analogue of the Sudakov's transport-ray decomposition.
\subsubsection{Second-order Beckmann problem}

The second-order Beckmann problem, introduced in \cite{Bolbotowski2024}, is concerned with the following variational task. Suppose that $\mu,\nu\in\mathcal{P}_2(\mathbb{R}^n)$ are two Borel probability measures with finite second moments and with common barycentre. We look for solutions to
\begin{equation}\label{eqn:beckmann}
    \inf\big\{\ \norm{\rho}_1\mid \rho\in\mathcal{M}(\mathbb{R}^n,\mathcal{S}^{n\times n}), \mathrm{div}^2\rho=\nu-\mu\big\},
\end{equation}
where $\mathcal{S}^{n\times n}$ denotes the set of all symmetric $n\times n$ matrices, $\mathcal{M}(\mathbb{R}^n,\mathcal{S}^{n\times n})$ is the space of all vector measures on $\mathbb{R}^n$ with values in $\mathcal{S}^{n\times n}$, $\norm{\rho}_1$ is the total variation norm with respect to the Schatten $1$-norm,\footnote{Let us remind that the Schatten $1$-norm of a matrix is the sum of its singular values.} and the condition $\mathrm{div}^2\rho=\nu-\mu$ reads that for all smooth, compactly supported $\phi\in\mathcal{C}_c^{\infty}(\mathbb{R}^n)$ we have
\begin{equation*}
    \int_{\mathbb{R}^n} \phi \, d(\nu-\mu)=\int_{\mathbb{R}^n}\langle D^2\phi, \,d\rho\rangle,
\end{equation*}
where $D^2\phi$ is the second derivative of $\phi$.
It is shown in \cite{Bolbotowski2024} that the infimum of (\ref{eqn:beckmann}) is attained and equal to the following infimum:
\begin{equation}\label{eqn:jot}
    \mathcal{J}(\mu,\nu)=\inf\Big\{\int_{(\mathbb{R}^n)^3}\frac12(\norm{x-z}^2+\norm{y-z}^2)\, d\sigma(x,y,z)\mid \sigma\in\Sigma(\mu,\nu)\Big\},
\end{equation}
where $\Sigma(\mu,\nu)$ is the set of all probabilities in $\mathcal{P}((\mathbb{R}^n)^3)$, that have $\mu$ and $\nu$ as their respective first and second marginal measures, and such that for all Borel, bounded $g,h\colon\mathbb{R}^n\to\mathbb{R}$
\begin{equation*}
    \int_{(\mathbb{R}^n)^3}g(x)(z-x)\, d\sigma(x,y,z)=\int_{(\mathbb{R}^n)^3}h(y)(z-y)\, d\sigma(x,y,z)=0.
\end{equation*}
In other words, the marginals $\mathrm{P}_{13}\sigma, \mathrm{P}_{23}\sigma$ are martingale couplings -- that is, they are distributions of one-step martingales, see Section~\ref{s:notation}.

As observed in \cite{Bolbotowski2024}, (\ref{eqn:beckmann}) is a relaxation of (\ref{eqn:jot}). Indeed, any plan $\sigma\in\Sigma(\mu,\nu)$ induces a measure $\rho\in\mathcal{M}(\mathbb{R}^n,\mathcal{S}^{n\times n})$ such that $\mathrm{div}^2\rho=\nu-\mu$ via the formula
\begin{equation}\label{eqn:grill}
    \rho=\int_{(\mathbb{R}^n)^3} \rho_{x,y,z}\, d\sigma(x,y,z),
\end{equation}
where $\rho_{x,y,z}\in\mathcal{M}(\mathbb{R}^n,\mathcal{S}^{n\times n})$ are given by specifying the densities
\begin{equation*}
    d\rho_{x,y,z}(\xi)=\norm{\xi-z}\Bigg(\frac{(x-z)(x-z)^*}{\norm{x-z}^2}d\mathcal{H}^1|_{[x,z]}(\xi)-\frac{(y-z)(y-z)^*}{\norm{y-z}^2}d\mathcal{H}^1|_{[y,z]}(\xi)\Bigg).
\end{equation*}
Here $\mathcal{H}^1$ is the one-dimensional Hausdorff measure on $\mathbb{R}^n$. Here $x^*$ stands for an adjoint map to $x$, that is for $y\in\mathbb{R}^n$, $x^*y=\langle x,y\rangle$. In general, if $T$ is linear, by $T^*$ we shall denote its adjoint.

\subsubsection{Dual problem and optimality conditions}\label{s:conditions}

The dual problem associated with (\ref{eqn:beckmann}) and with (\ref{eqn:jot}) is given by
\begin{equation}\label{eqn:dual}
    \mathcal{I}(\mu,\nu)=\sup\Big\{\int_{\mathbb{R}^n}u\, d(\nu-\mu)\mid u\in\mathcal{C}^{1,1}(\mathbb{R}^n)\text{ has }1\text{-Lipschitz derivative}\Big\}.
\end{equation}
In \cite[Theorem 1.1]{Bolbotowski2024} it is shown that  for any $\mu,\nu\in\mathcal{P}_2(\mathbb{R}^n)$ with common barycentre optimal $u\in\mathcal{C}^{1,1}(\mathbb{R}^n)$ and optimal $\sigma\in\Sigma(\mu,\nu)$ exist, and that $\mathcal{I}(\mu,\nu)=\mathcal{J}(\mu,\nu)$. Moreover, the optimisers of $\mathcal{I}(\mu,\nu)$ and of $\mathcal{J}(\mu,\nu)$ are characterised by the equality
\begin{equation}\label{eqn:optimality}
    u(y)+Du(y)(z-y)-u(x)-Du(x)(z-x)=\frac12(\norm{x-z}^2+\norm{y-z}^2)
\end{equation}
that has to hold for $\sigma$-almost every $(x,y,z)\in(\mathbb{R}^n)^3$.

In Section \ref{s:derivative}, Lemma \ref{lem:three} we  show that this optimality condition implies that the $1$-Lipschitz map $Du\colon\mathbb{R}^n\to\mathbb{R}^n$ is isometric on $\{x,z\}$ and on $\{y,z\}$ for any $(x,y,z)$ that satisfy (\ref{eqn:optimality}). Moreover, $Du(x)-Du(z)=x-z$ and $Du(y)-Du(z)=z-y$ for such points, as we will see in Lemma \ref{lem:three}.

\subsubsection{Leaf decompositions}

The last observation leads us to the topic of leaf decompositions, which will be our analogue of the partition into  transport rays considered in Section \ref{s:transportrays}.

For any $1$-Lipschitz map $v\colon\mathbb{R}^n\to\mathbb{R}^m$, the leaves of $v$ are defined as the maximal sets $\mathcal{S}$ such that $v|_{\mathcal{S}}$ is isometric. It is shown in \cite{Ciosmak2021} that such leaves are closed and convex sets and that they partition $\mathbb{R}^n$ -- up to a set of the Lebesgue measure zero.

As one of the contributions of the paper, we further the understanding of the leaves of any $1$-Lipschitz map, and show that the intersection of any two distinct leaves $\mathcal{S}_1,\mathcal{S}_2$ is not only contained in the intersection of their relative boundaries $\partial\mathcal{S}_1\cap\partial\mathcal{S}_2$, as proven in \cite[Lemma 2.11]{Ciosmak2021}, but is precisely the common convex face of the two convex sets $\mathcal{S}_1$ and $\mathcal{S}_2$.

We also characterise the local compatibility of two intersecting leaves $\mathcal{S}_1,\mathcal{S}_2$ of $1$-Lipschitz $Du$ with $u\in\mathcal{C}^{1,1}(\mathbb{R}^n)$. On each leaf $Du$ acts as a difference of two orthogonal projections.
The positive directions of either leaf of $\mathcal{S}_1,\mathcal{S}_2$ must form obtuse angles with the
negative directions of the other. We prove a converse -- Theorem \ref{thm:compat} -- every compatible pair of leafwise quadratic functions
extends to a global $\mathcal C^{1,1}(\mathbb{R}^n)$ function with $1$-Lipschitz
derivative. Figure~\ref{fig:leaf-graphs} illustrates two leaves whose transverse
curvatures have opposite signs. Figures~\ref{fig:two-equality-configurations}
and~\ref{fig:one-leaf} show contact triples on patches of two distinct leaves and on a single patch, respectively.

\begin{figure}[ht]
\centering
\begin{tikzpicture}
 \begin{groupplot}[
  group style={group size=1 by 1, horizontal sep=2.0cm},
   every axis/.append style={font=\small}
 ]

\nextgroupplot[
    width=11.1cm,
    height=8.25cm,
    view={-60}{40},
    xmin=-2.35, xmax=2.60,
    ymin=-2.25, ymax=2.25,
    zmin=0, zmax=3.0,
    axis lines=center,
     xlabel={$x_1$}, ylabel={$x_2$},
      zlabel={$x_3$},
    xlabel style={at={(ticklabel* cs:1.02)}, anchor=west},
    ylabel style={at={(ticklabel* cs:1.02)}, anchor=south},
     zlabel style={at={(ticklabel* cs:1.02)}, anchor=south},
    xtick=\empty, ytick=\empty, ztick=\empty,
    enlargelimits=false,
 clip=false,
    domain=-2:2,
    samples=19,
    samples y=11,
    colormap/blackwhite,
]

\addplot3[
    patch, patch type=polygon, vertex count=4,
    draw=black!34, line width=0.38pt,
    fill=gray!70,
      fill opacity=0.7,
] coordinates {(-2,0,0) (2,0,0) (2,2,0) (-2,2,0)};

\addplot3[
    patch, patch type=polygon, vertex count=4,
    draw=black!34, line width=0.38pt,
    fill=gray!50,
     fill opacity=0.7,
     postaction={pattern= north west lines, pattern color=black!90}
] coordinates {(-2,-2,0) (2,-2,0) (2,0,0) (-2,0,0)};

\foreach \x in {-1.90,-1.55,...,1.60} {
  \addplot3[black!17, line width=0.18pt]
  coordinates {(\x,-2,0) (\x+0.35,0,0)};
}
\addplot3[black!92, line width=0.92pt] coordinates {(-2,0,0) (2,0,0)};

\addplot3[
    surf,
    y domain=0:2,
    shader=flat,
    draw=black!40,
    fill opacity=0.1,
    opacity=1
] ({x},{y},{0.5*(x^2+y^2)+4});

\addplot3[
    surf,
    y domain=-2:0,
    shader=flat,
    draw=black!50,
    fill opacity=0.1,
    opacity=1
] ({x},{y},{0.5*(x^2-y^2)+4});

\addplot3[
    black!80,
    line width=1pt,
    smooth
]
coordinates {
(-2,0,6)
(-1.8,0,5.62)
(-1.6,0,5.28)
(-1.4,0,4.98)
(-1.2,0,4.72)
(-1,0,4.5)
(-0.8,0,4.32)
(-0.6,0,4.18)
(-0.4,0,4.08)
(-0.2,0,4.02)
(0,0,4)
(0.2,0,4.02)
(0.4,0,4.08)
(0.6,0,4.18)
(0.8,0,4.32)
(1,0,4.5)
(1.2,0,4.72)
(1.4,0,4.98)
(1.6,0,5.28)
(1.8,0,5.62)
(2,0,6)
};

 \addplot3[
    black!84,
    line width=1pt,
    smooth
]
coordinates {
(-2,2,8)
(-1.8,2,7.62)
(-1.6,2,7.28)
(-1.4,2,6.98)
(-1.2,2,6.72)
(-1.0,2,6.50)
(-0.8,2,6.32)
(-0.6,2,6.18)
(-0.4,2,6.08)
(-0.2,2,6.02)
(0,2,6)
(0.2,2,6.02)
(0.4,2,6.08)
(0.6,2,6.18)
(0.8,2,6.32)
(1.0,2,6.50)
(1.2,2,6.72)
(1.4,2,6.98)
(1.6,2,7.28)
(1.8,2,7.62)
(2,2,8)
};
\addplot3[
    black!84,
    line width=1pt,
    smooth
]
coordinates {
(-2,-2,4)
(-1.8,-2,3.62)
(-1.6,-2,3.28)
(-1.4,-2,2.98)
(-1.2,-2,2.72)
(-1.0,-2,2.50)
(-0.8,-2,2.32)
(-0.6,-2,2.18)
(-0.4,-2,2.08)
(-0.2,-2,2.02)
(0,-2,2)
(0.2,-2,2.02)
(0.4,-2,2.08)
(0.6,-2,2.18)
(0.8,-2,2.32)
(1.0,-2,2.50)
(1.2,-2,2.72)
(1.4,-2,2.98)
(1.6,-2,3.28)
(1.8,-2,3.62)
(2,-2,4)
};
\addplot3[black!84, line width=1pt, domain=0:2, samples=60]
({-2},{y},{0.5*(4+y^2)+4});
\addplot3[black!84, line width=1pt, domain=0:2, samples=60]
({2},{y},{0.5*(4+y^2)+4});
\addplot3[black!84, line width=1pt, domain=-2:0, samples=60]
({-2},{y},{0.5*(4-y^2)+4});
\addplot3[black!84, line width=1pt, domain=-2:0, samples=60]
({2},{y},{0.5*(4-y^2)+4});

\node[
    inner sep=0.8pt,
    font=\small
] at (axis cs:-1.10,0.4,0) {$\mathcal{S}_1$};

\node[
    inner sep=0.8pt,
    font=\small
] at (axis cs:-1.10,-0.4,0) {$\mathcal{S}_2$};

\node[
    anchor=north west,
    inner sep=0.8pt,
    font=\footnotesize
] at (axis cs:-3.4,0,0.0) {$\mathcal{F}_{12}=\mathcal{S}_1\cap \mathcal{S}_2$};

\node[
    inner sep=1.15pt,
    font=\footnotesize
] at (axis cs:3.1,1.0,7.0)
{$u|_{\mathcal{S}_1}(x)=\frac{1}{2}(x_1^2+x_2^2)$};

\node[
    inner sep=1.15pt,
    font=\footnotesize
] at (axis cs:3.35,-0.96,5.5)
{$u|_{\mathcal{S}_2}(x)=\frac{1}{2}(x_1^2-x_2^2)$};

 \end{groupplot}

\end{tikzpicture}

\caption{
Two adjacent leaves $\mathcal{S}_1$ and $\mathcal{S}_2$, which intersect along
the common face $\mathcal{F}_{12}=\mathcal{S}_1\cap\mathcal{S}_2$. The potential satisfies $Du|_{\mathcal{S}_1}=\mathrm{Id}$, $Du|_{\mathcal{S}_2}=P_{\mathbb{R}e_1}-P_{\mathbb{R}e_2}$.
}\label{fig:leaf-graphs}
\end{figure}

According to the observation made in Section \ref{s:conditions}, for any optimal  $\sigma\in\Sigma(\mu,\nu)$, $\{x,z\}$ and $\{y,z\}$ belong to leaves of $Du$, for $\sigma$-almost every $(x,y,z)$. Here $u$ is optimal for $\mathcal{I}(\mu,\nu)$.
Thanks to the results of \cite{Ciosmak2021} we can disintegrate $\mu,\nu$ with respect to the foliation induced by the leaves $\mathcal{S}$ of $Du$ and obtain Borel measures $\mu_{\mathcal{S}},\nu_{\mathcal{S}}$ for each such leaf.
Since the leaves are disjoint up to set of measure zero, we can show that any optimal $\sigma\in\Sigma(\mu,\nu)$ with all three marginals $\mu,\nu,\mathrm{P}_3\sigma$ absolutely continuous completely decomposes into optimal plans $\sigma_{\mathcal{S}}\in\Sigma(\mu_{\mathcal{S}},\nu_{\mathcal{S}})$. The same conclusion holds for any optimal  $\sigma\in\Sigma(\mu,\nu)$ if there exists an optimal $\sigma_0$ whose three marginals are absolutely continuous.

\subsubsection{Bimartingale couplings and convex-concave order}\label{s:bimartconvex}

As in Section \ref{s:monotone}, the fact that any optimal  $\sigma\in\Sigma(\mu,\nu)$ with absolutely continuous marginals is constrained within the leaves of $Du$ allows us to reduce the analysis of $\sigma$ to the analysis of   $\sigma_{\mathcal{S}}\in\Sigma(\mu_{\mathcal{S}},\nu_{\mathcal{S}})$ for each leaf $\mathcal{S}$ of $Du$.
The advantage of each of the reduced problems is that the optimiser of its dual problem has isometric derivative on the leaf.

We prove in Proposition \ref{pro:cc} and Theorem \ref{thm:bimart} that the fact that an optimiser of the dual problem (\ref{eqn:dual}) for a pair $\mu,\nu\in\mathcal{P}_2(\mathbb{R}^n)$ with common barycentre has isometric derivative is equivalent to
\begin{enumerate}[label=(\roman*)]
    \item $\mu\prec_{c-c}\nu$ being in convex-concave order,
    \item the existence of a \emph{bimartingale coupling} $\pi\in\Gamma(\mu,\nu)$.
\end{enumerate}
The first of the above conditions requires that there exist two mutually complementing, orthogonal subspaces $V_1,V_2\subset\mathbb{R}^n$ such that for any function $f\colon\mathbb{R}^n\to\mathbb{R}$ that is convex in the $V_1$-variable and concave in the $V_2$-variable, and of at most quadratic growth, we have
\begin{equation*}
    \int_{\mathbb{R}^n}f\, d\mu\leq\int_{\mathbb{R}^n}f\, d\nu.
\end{equation*}
A coupling  $\pi\in\Gamma(\mu,\nu)$ is called a bimartingale coupling with respect to $V_1,V_2$ if and only if for any pair of random vectors $(X,Y)\sim \pi$ we have
\begin{equation*}
    \mathbb{E}\big(P_{V_1}Y\mid \sigma(X)\big)=P_{V_1}X\text{ and } \mathbb{E}\big(P_{V_2}X\mid \sigma(Y)\big)=P_{V_2}Y.
\end{equation*}
Here $P_{V_1},P_{V_2}$ are the orthogonal projections onto $V_1,V_2$, respectively.

Moreover, if $\pi\in\Gamma(\mu,\nu)$ is a bimartingale coupling with respect to $V_1,V_2$, then its pushforward $R_{\#}\pi$ through the map $R\colon\mathbb{R}^n\times\mathbb{R}^n\to(\mathbb{R}^n)^3$ defined by the formula
\begin{equation*}
    R(x,y)=(x,y,P_{V_2}x+P_{V_1}y)\text{ for }x,y\in\mathbb{R}^n,
\end{equation*}
belongs to $\Sigma(\mu,\nu)$ and is an optimal plan for $\mathcal{J}(\mu,\nu)$, see Theorem \ref{thm:bimart}.

In particular, we complement a characterisation of a pair of measures in convex order --  these are measures $\mu,\nu\in\mathcal{P}_2(\mathbb{R}^n)$ for which
the map $\frac12\norm{\cdot}^2$ is an optimiser of
\begin{equation*}
    \sup\Big\{\int_{\mathbb{R}^n}u\, d(\nu-\mu)\mid u\in\mathcal{C}^{1,1}(\mathbb{R}^n) \text{ with }1\text{-Lipschitz derivative} \Big\}.
\end{equation*}
Note that if $V_2=\{0\}$ our observations retrieve  the classical theorem of Strassen \cite{Strassen1965} that shows that two probabilities admit a martingale coupling if and only if they are in convex order.

\subsubsection{Covariance matrices}

In general, if we are in the situation described in Section \ref{s:bimartconvex}, then the choice of subspaces $V_1,V_2\subset\mathbb{R}^n$ with respect to which $\mu,\nu\in\mathcal{P}_2(\mathbb{R}^n)$ are in convex-concave order need not be unique. However, this is not the case if the difference of the covariance matrices of $\nu$ and of $\mu$
\begin{equation*}
    C=\int_{\mathbb{R}^n}xx^*\, d(\nu-\mu)(x)
\end{equation*}
is non-degenerate, as shown by Lemma \ref{lem:isoderma} and Proposition \ref{pro:optimal}.
However, if $C$ is degenerate, we can consecutively partition the space, eventually arriving at a situation at which the corresponding matrices are non-degenerate. This procedure is described in Section \ref{s:refine}.

\subsubsection{Characterisation of optimal plans}

The above analysis gives a full description of all optimal plans in $\Sigma(\mu,\nu)$, under the assumption that there exists some optimal plan with absolutely continuous marginals. Note that this requires that $\mu,\nu\in\mathcal{P}_2(\mathbb{R}^n)$ are absolutely continuous measures with common barycentre.
Below by $CC(\mathbb{R}^n)$ we denote the set of all closed, convex and non-empty subsets of $\mathbb{R}^n$ and $\norm{\cdot}_1$ denotes the Schatten $1$-norm of a matrix.

\begin{theorem}\label{thm:descr}
    Let $\mu,\nu\in\mathcal{P}_2(\mathbb{R}^n)$ be absolutely continuous with respect to the Lebesgue measure and have common barycentre. Suppose that there exists an optimal plan $\sigma_0\in\Sigma(\mu,\nu)$ with absolutely continuous third marginal.
    Then there exist a family of closed and convex sets $\mathcal{S}\in CC(\mathbb{R}^n)$ such that each point outside a single Lebesgue-null set belongs to exactly one member of the family, a Borel probability measure $\theta\in\mathcal{P}(CC(\mathbb{R}^n))$ and Borel probability measures
    \begin{equation*}
        \mu_{\mathcal{S}},\nu_{\mathcal{S}}\in\mathcal{P}_2(\mathbb{R}^n) \text{ defined for all }\mathcal{S}\in CC(\mathbb{R}^n)
    \end{equation*}
    such that for $\theta$-almost every $\mathcal{S}$:
    \begin{enumerate}[label=(\roman*)]
        \item\label{i:concen} $\mu_{\mathcal{S}},\nu_{\mathcal{S}}$ are concentrated on $\mathcal{S}$
\item\label{i:bary} $\mu_{\mathcal{S}},\nu_{\mathcal{S}}$ have common barycentre,
   \item\label{i:covar} the difference of the covariance matrices
      \begin{equation*}
            C_{\mathcal{S}}=\int_{\mathbb{R}^n}xx^*\, d(\nu_{\mathcal{S}}-\mu_{\mathcal{S}})(x)
        \end{equation*}
        is non-degenerate on the tangent space $V(\mathcal{S})$ to $\mathcal{S}$,

\item\label{i:conconv} $\mu_{\mathcal{S}}\prec_{c-c}\nu_{\mathcal{S}}$ is in convex-concave order with respect to $V(\mathcal{S})_1,V(\mathcal{S})_2$, where $V(\mathcal{S})_1,V(\mathcal{S})_2$ are mutually complementing, orthogonal subspaces of $V(\mathcal{S})$,
         $V(\mathcal{S})_1$ is the sum of positive eigenspaces of the difference of the covariance matrices $C_{\mathcal{S}}$ and $V(\mathcal{S})_2$ is the  sum of negative eigenspaces of $C_{\mathcal{S}}$,
        \item\label{i:set} for any Borel $A\subset\mathbb{R}^n$ the map
    \begin{equation*}
        CC(\mathbb{R}^n)\ni\mathcal{S}\mapsto (\mu_{\mathcal{S}}(A),\nu_{\mathcal{S}}(A))\in\mathbb{R}\times\mathbb{R}
    \end{equation*}
    is Borel measurable and
    \begin{equation*}
        \mu(A)=\int_{CC(\mathbb{R}^n)}\mu_{\mathcal{S}}(A)\, d\theta(\mathcal{S})\text{ and }        \nu(A)=\int_{CC(\mathbb{R}^n)}\nu_{\mathcal{S}}(A)\, d\theta(\mathcal{S}).
    \end{equation*}
    \end{enumerate}
    Moreover, a plan $\sigma\in \mathcal{P}((\mathbb{R}^n)^3)$
    belongs to $ \Sigma(\mu,\nu)$  and is optimal for $\mathcal{J}(\mu,\nu)$ if and only if it is of the form\footnote{For a Borel map $T\colon X\to Y$ and a Borel measure $\mu$ on $X$, by $T_{\#}\mu$ we denote  the pushforward of $\mu$ via $T$.}
    \begin{equation}\label{eqn:sigmadis}
        \sigma=\int_{CC(\mathbb{R}^n)}(R_{\mathcal{S}})_{\#}\pi_{\mathcal{S}}\, d\theta(\mathcal{S}),
    \end{equation}
    for some map
    \begin{equation*}
        CC(\mathbb{R}^n)\ni \mathcal{S}\mapsto\pi_{\mathcal{S}}\in \mathcal{P}(\mathbb{R}^n\times\mathbb{R}^n)
    \end{equation*}
    such that for $\theta$-almost every $\mathcal{S}$:
    \begin{enumerate}[label=(\alph*)]
    \item\label{i:displan} for any Borel set $B\subset\mathbb{R}^n\times\mathbb{R}^n$ the map
    \begin{equation*}
        CC(\mathbb{R}^n)\ni \mathcal{S}\mapsto \pi_{\mathcal{S}}(B)\in\mathbb{R}
    \end{equation*}
    is Borel measurable,
        \item\label{i:bimartplan} $\pi_{\mathcal{S}}\in \Gamma_{bm}(\mu_{\mathcal{S}},\nu_{\mathcal{S}}, V(\mathcal{S})_1,V(\mathcal{S})_2)$ is a bimartingale coupling,
        \item\label{i:rpush} $R_{\mathcal{S}}(x,y)=(x,y,z(\mathcal{S})+P_{V(\mathcal{S})_2}(x-z(\mathcal{S}))+P_{V(\mathcal{S})_1}(y-z(\mathcal{S})))$ for all $x,y\in\mathbb{R}^n$ and some Borel map
        \begin{equation*}
            z\colon CC(\mathbb{R}^n)\to\mathbb{R}^n\text{ such that }z(\mathcal{S})\in\mathcal{S}.
        \end{equation*}
    \end{enumerate}
    Furthermore, the optimal cost equals
    \begin{equation}\label{eqn:cost}
       \mathcal{J}(\mu,\nu)=\frac12\int_{CC(\mathbb{R}^n)}\norm{C_{\mathcal{S}}}_1\, d\theta(\mathcal{S}).
    \end{equation}
\end{theorem}

The proof of the theorem relies first on decomposition of the space into leaves of the derivative of an optimiser of $\mathcal{I}(\mu,\nu)$. Second, each leaf we further decompose so as to ensure that the differences of the covariances are non-degenerate. Third, optimal plans for such measures are completely described through bimartingale couplings with respect to the subspaces on which the  difference of the covariances is positive definite and negative definite, respectively. The resulting partition may depend on the choice of the dual optimiser. The details of the proof are presented in Section \ref{s:descr}.

\subsubsection{Absolute continuity and measure contraction property}\label{s:mcp}

Let $\eta$ be a finite absolutely continuous measure on $\mathbb R^n$.
By \cite[Theorem 3.1]{ciosmakMCP2026}, for almost every isometric leaf
$\mathcal S$ of a $1$-Lipschitz map, the conditional measure $\eta_{\mathcal S}$
is absolutely continuous with respect to $\mathcal H^{\dim\mathcal S}|_{\mathcal S}$.
This applies in particular to the leaves of $Du$. If $\eta$ has positive
density almost everywhere on its convex support and satisfies
$\operatorname{MCP}(\kappa,N)$, with $1<N<\infty$ and $N\geq n$, then
almost every conditional measure satisfies $\operatorname{MCP}(\kappa,N)$
on its support, by \cite[Theorem A]{ciosmakMCP2026}.
See \cite{ohta2007mcp} for the definition of the measure contraction property $\operatorname{MCP}(\kappa,N)$. Both conclusions also hold for the conditionals $\mu_{\mathcal S},\nu_{\mathcal S}$ relative to the partition constructed in Theorem~\ref{thm:descr}, under the respective hypotheses on each marginal.

Indeed, each refinement slices a leaf $\mathcal S$ by translates of the range $\operatorname{ran} C_{\mathcal S}$. By~\eqref{eqn:leafcovariance}, $C_{\mathcal S}$ is the difference of the covariance matrices of $\nu_{\mathcal S}$ and $\mu_{\mathcal S}$. These conditional measures are obtained by disintegrating $\nu$ and $\mu$, respectively, with respect to the partition at that level, as in \hyperref[par:refinedconditionals]{Section~\ref*{s:proof}}. In the affine hull $\operatorname{aff}\mathcal S$, these translates are the isometric leaves of the orthogonal projection onto $\operatorname{ran} C_{\mathcal S}$. Thus \cite[Theorem 3.1]{ciosmakMCP2026} gives absolute continuity at the next step. Under the additional hypotheses above, \cite[Theorem A]{ciosmakMCP2026} also gives positive density almost everywhere on the intersection of the leaf with the convex support, and $\operatorname{MCP}(\kappa,N)$ there. This intersection is the convex support of the slice, and $N\geq n\geq\dim\mathcal S$, so the same theorem applies at each subsequent step. Since every proper refinement lowers the dimension, the assertion follows by induction.

For earlier related results on isometric leaves and on the faces of a convex
function, see also \cite{Ciosmak2021} and \cite{Caravenna2010}, respectively.

\subsubsection{Further decompositions and irreducible convex paving}\label{s:irreducible}

The decomposition of the space described in Theorem \ref{thm:descr} is in general not the finest one that yields constraints for all optimal plans, analogously to the situation in Section \ref{s:fur}.

Suppose that $\mu,\nu\in\mathcal{P}_2(\mathbb{R}^n)$ have common barycentre and that $\mu\prec_{c-c}\nu$ are in convex-concave order with respect to two mutually complementing, orthogonal subspaces $V_1,V_2$. Then for any bimartingale coupling $\pi\in\Gamma(\mu,\nu)$ with respect to $V_1,V_2$, the pushforwards
\begin{equation*}
    (P_{V_1},P_{V_1})_{\#}\pi\in \Gamma((P_{V_1})_{\#}\mu,(P_{V_1})_{\#}\nu)\text{ and }   S(P_{V_2},P_{V_2})_{\#}\pi\in \Gamma((P_{V_2})_{\#}\nu,(P_{V_2})_{\#}\mu)
\end{equation*}
are martingale transports. Here $S\colon V_2\times V_2\to V_2\times V_2$ is a map that swaps the variables.

The study of martingale transports -- conceived in \cite{Beiglbock2013} -- stems from financial mathematics. A crucial concept of the theory is the notion of irreducible components and the related convex paving, investigated in \cite{Juillet2016} and \cite{Nutz2017} in the one-dimensional case, in \cite{Touzi2019} in the multi-dimensional case, and vastly extended in \cite{Ciosmak2024} and in \cite{Ciosmak20241}.

It turns out that to any pair of probability measures in convex order one may associate the finest partition of the ambient space into convex sets, such that any  martingale coupling has to adhere to this paving. Since the aforementioned pushforwards  for any bimartingale coupling are martingale couplings, we see that any bimartingale coupling between $\mu$ and $\nu$ is constrained within the products of the corresponding irreducible components of martingale transports for the pairs $\big((P_{V_1})_{\#}\mu,(P_{V_1})_{\#}\nu\big)$ and $\big((P_{V_2})_{\#}\nu,(P_{V_2})_{\#}\mu\big)$.

\subsection{General absolutely continuous measures}

The proof of Theorem \ref{thm:descr} relies on the assumption that there exists an optimal $\sigma\in\Sigma(\mu,\nu)$ with absolutely continuous marginals. If no such assumption is made, the description need not be valid, see Example~\ref{exa:atomicthird}. However, in the general case, the partition of the space into leaves of $Du$ for an optimal $u\in\mathcal{C}^{1,1}(\mathbb{R}^n)$ with $1$-Lipschitz derivative is still relevant to the description of solutions to $\mathcal{J}(\mu,\nu)$. Theorem \ref{thm:noassu} shows that for general absolutely continuous measures with common barycentre, any solution of $\mathcal{J}(\mu,\nu)$ can be decomposed into parts that charge only triples of points in the neighbouring leaves of $Du$ and in their interfaces.

\subsection{Relation to optimal grillage}

Here we describe a relation of the considered problem to the problem of design of an optimal grillage that was already presented in \cite{Bolbotowski2024}.

Let us consider the two-dimensional case. Let $f\in \mathcal{M}_2(\mathbb{R}^2)$ be a signed, absolutely continuous measure with total mass zero,  with barycentre equal to zero and with finite second moments. Then the positive part and negative part of $f$, up to normalisation, satisfy the assumptions of Theorem \ref{thm:descr}, provided that there is an optimal plan with absolutely continuous third marginal. In a plate subjected to the load $f$,  any measure $\rho\in\mathcal{M}(\mathbb{R}^2,\mathcal{S}^{2\times 2})$ such that
\begin{equation*}
    \mathrm{div}^2\rho=f,
\end{equation*}
 encodes the distribution of bending moment tensor throughout the plate.
If the bending moment tensor decomposes as in (\ref{eqn:grill}), then it is called a grillage, as in this case the plate can be described in terms of straight bars.
A problem that has its applications in optimal design, see, e.g., \cite{Bolbotowski2018}, \cite{Prager1977}, \cite{Rozvazny1972}, is to find an optimal arrangement of the bars so as to minimise  a certain energy functional. The existence of optimisers was first proved by Bo\l{}botowski and Bouchitt\'e in \cite{Bolbotowski2024}. Our work builds upon \cite{Bolbotowski2024} and gives a fine structure theorem for optimal solutions. We also open a path towards establishing a structure theorem and characterisation of polar sets with respect to optimal solutions of $\mathcal{J}(\mu,\nu)$. Here a set is said to be polar with respect to a family of measures whenever it is negligible for any measure in the family. This result would provide a comprehensive understanding of the design flexibility available for optimal grillage configurations.

\subsection{Further research}

\subsubsection{Polar sets}\label{s:polar}

Section \ref{s:irreducible} provides a promising avenue for future research.

\begin{question}
Suppose that $\mu,\nu\in\mathcal{P}_2(\mathbb{R}^n)$ are in convex-concave order with respect to mutually complementing, orthogonal subspaces $V_1,V_2\subset\mathbb{R}^n$. What are the polar sets with respect to all bimartingale couplings in between $\mu$ and $\nu$ with respect to $V_1,V_2$?
\end{question}

Let us note that any bimartingale coupling between $\mu$ and $\nu$ is constrained within the products of the irreducible components pertaining to martingale transports between $(P_{V_1})_{\#}\mu$ and $(P_{V_1})_{\#}\nu$ and between $(P_{V_2})_{\#}\nu$ and $(P_{V_2})_{\#}\mu$ .

Providing an answer to the above question would also, thanks to Theorem \ref{thm:descr}, give a characterisation of polar sets with respect to all optimal couplings between two given probabilities $\mu,\nu\in\mathcal{P}_2(\mathbb{R}^n)$ with common barycentre which admit an optimal $\sigma_0\in\Sigma(\mu,\nu)$ with absolutely continuous marginals.
Related problems were studied in \cite{Ciosmak2024} and in \cite{Ciosmak20241}.

\subsubsection{Absolute continuity and measure contraction property for the fine paving}
Do the conclusions of Section~\ref{s:mcp} extend to the analogue of the
irreducible convex paving discussed in Section~\ref{s:irreducible} and in
Section~\ref{s:polar}? Suppose that such a paving gives a measurable
partition up to a Lebesgue-null set. If Lebesgue measure is disintegrated
with respect to this paving, is the conditional measure on almost every
cell $C$ absolutely continuous with respect to $\mathcal H^{\dim C}$?
More generally, if an ambient measure $\eta$ satisfies the density and
measure contraction hypotheses, does almost every conditional measure
satisfy $\operatorname{MCP}(\kappa,N)$ on its support?

\subsection{Literature}

Similar Beckmann-type functionals of the form
\begin{equation*}
    \sup\Big\{\int_{\mathbb{R}^n}u\, d(\nu-\mu)\mid u\in\mathcal{C}^{\infty}(\mathbb{R}^n,\mathbb{R}^m), \norm{Au}\leq 1\Big\},
\end{equation*}
where $\norm{\cdot}$ is a seminorm on the space of $m\times n$ matrices, $A$ is a linear differential operator, such that $Au\colon \mathbb{R}^n\to \mathbb{R}^{m\times n}$ for any $u\in\mathcal{C}^{\infty}(\mathbb{R}^n,\mathbb{R}^m)$,  were investigated e.g. in \cite{Bouchitte2001}, \cite{Bolbotowski2022}.

Other similar problems in the setting of vector measures were studied in \cite{Ciosmak20212} and in  \cite{Ciosmak20211}.
The mass balance condition studied in Section \ref{s:balance} fails in these settings. This failure is closely related to the continuity of extensions of $1$-Lipschitz vector-valued maps, a topic studied thoroughly in \cite{Ciosmak20213}, \cite{Ciosmak20242}.

\subsection{Organisation of the paper}

The paper is organised as follows. In Section \ref{s:par} we define the partitions associated to any $1$-Lipschitz map and study their properties, see Section \ref{s:leaves}. In Section \ref{s:derivative} the leaves of $1$-Lipschitz maps that are derivatives are considered.  Section \ref{s:notation} is devoted to the introduction of preliminary notation that we will use throughout the paper. In Section \ref{s:duality} we recall the results of \cite{Bolbotowski2024} on the Kantorovich--Rubinstein duality for the Hessian. In Section \ref{s:balance} we prove the mass and the moment balance conditions. In Section \ref{s:disintegrate} we show that in the case of existence of optimiser with absolutely continuous marginals,  the problem considered in this paper completely decomposes into a collection of problems on the leaves of an optimiser of the dual problem.

In Section \ref{s:bimart} we study in detail the case of a pair of measures for which an optimiser of the dual problem has isometric derivative. This includes the introduction and the study of the convex-concave order in Section \ref{s:convexconcave} and the notion of bimartingale coupling in Section \ref{s:couplings}. Their mutual relation is explored in Section \ref{s:relation}. Section \ref{s:covariance} comprises considerations on the covariance matrices.

Section \ref{s:descr} is concerned with the final proof of Theorem \ref{thm:descr}. This is preceded by the measurability considerations in Section \ref{s:refine}. Section \ref{s:proof} completes the argument.

Section \ref{s:general} in Theorem \ref{thm:noassu} discusses the relevance of the leaf decomposition to the solution of the $\mathcal{J}(\mu,\nu)$ problem in the case of general absolutely continuous measures.

\section{Partition, mass balance and moment balance conditions}\label{s:par}

\subsection{Leaves of 1-Lipschitz maps -- partitions}\label{s:leaves}
We start by recalling some definitions concerned with \emph{leaves} of $1$-Lipschitz maps. We refer the reader to \cite{Ciosmak2021} for the details.

\begin{definition}
    Let $v\colon\mathbb{R}^n\to\mathbb{R}^m$ be $1$-Lipschitz. A \emph{leaf} $\mathcal{S}\subset\mathbb{R}^n$ of $v$ is a maximal set with respect to inclusion such that for all $x,y\in\mathcal{S}$
    \begin{equation*}
        \norm{v(x)-v(y)}=\norm{x-y}.
    \end{equation*}
    A \emph{transport set} of $v$ is a Borel set in $\mathbb{R}^n$ that is a union of leaves.
\end{definition}

In \cite[Corollary 2.5]{Ciosmak2021} it is shown that each leaf of a $1$-Lipschitz map is closed and convex.

By $CC(\mathbb{R}^n)$, we denote the space of closed, convex, non-empty subsets of $\mathbb{R}^n$. It is a closed subspace of $CL(\mathbb{R}^n)$ -- the space of closed non-empty subsets of $\mathbb{R}^n$ equipped with the Wijsman topology -- see \cite{Wijsman1964}. The Wijsman topology is the weakest topology such that for any $x\in\mathbb{R}^n$ function
\begin{equation*}
A\mapsto \mathrm{dist}(x,A)
\end{equation*}
is continuous. By a result of Beer -- see \cite{Beer1991} -- the space $CL(\mathbb{R}^n)$, equipped with this topology, is Polish. Hence so is $CC(\mathbb{R}^n)$, as a closed subset of $CL(\mathbb{R}^n)$.

 We shall consider the map
 \begin{equation*}
    \mathcal{S}\colon\mathbb{R}^n\to CC(\mathbb{R}^n)
\end{equation*}
that assigns to $x\notin N(v)$ the unique leaf of $v$ that contains $x$, and such that $\mathcal{S}(x)=\{x\}$ for $x\in N(v)$, where $N(v)$ is the set of points in $\mathbb{R}^n$ at which $v$ is not differentiable. Note that  by \cite[Corollary 2.15]{Ciosmak2021} $B(v)\subset N(v)$, where $B(v)$ is a set of points in $\mathbb{R}^n$ that belong to at least two distinct leaves of $v$ and that $N(v)$ is of the Lebesgue measure zero by the Rademacher theorem.

In \cite[Lemma 2.11]{Ciosmak2021} it is shown that the leaves of $v$ are pairwise disjoint, up to the boundaries. That is, if $\mathcal{S}_1,\mathcal{S}_2$ are two distinct leaves of a $1$-Lipschitz map $v\colon\mathbb{R}^n\to\mathbb{R}^m$, then
\begin{equation*}
    \mathcal{S}_1\cap\mathcal{S}_2\subset\partial\mathcal{S}_1\cap\partial\mathcal{S}_2.
\end{equation*}

Below, for a finite-dimensional convex set $\mathcal{F}$, by $\mathrm{int}\mathcal{F}$ we will denote the relative interior of $\mathcal{F}$. That is, $\mathrm{int}\mathcal{F}$ is the interior of $\mathcal{F}$ with respect to the topology of the affine hull of $\mathcal{F}$. A face of a convex set $\mathcal{F}$ is a convex subset $\mathcal{G}$ of $\mathcal{F}$ such that if for some $t\in (0,1)$ and some $x,y\in\mathcal{F}$ we have $tx+(1-t)y\in\mathcal{G}$, then $x,y\in\mathcal{G}$. Note that for any $x\in\mathcal{F}$ there exists a unique face of $\mathcal{F}$ that contains $x$ in its relative interior; equivalently, this is the minimal face of $\mathcal{F}$ that contains $x$. We refer to \cite{Rockafellar1970} for an account of the subject and in particular to \cite[Theorem 18.2]{Rockafellar1970}.

We now provide a strengthening of \cite[Lemma 2.11]{Ciosmak2021}.

\begin{lemma}\label{lem:intersection}
    If $\mathcal{S}_1,\mathcal{S}_2$ are two leaves of a $1$-Lipschitz map $v\colon\mathbb{R}^n\to\mathbb{R}^m$ and $\mathcal{S}_1\cap\mathcal{S}_2\neq\emptyset$, then the intersection $\mathcal{S}_1\cap\mathcal{S}_2$ is a common face of $\mathcal{S}_1$ and of $\mathcal{S}_2$.
\end{lemma}
\begin{proof}
Let $z\in\mathcal{S}_{1}\cap\mathcal{S}_2$. Suppose that for some $t\in (0,1)$ and some $x,y\in\mathcal{S}_1$ we have
\begin{equation*}
    z=tx+(1-t)y.
\end{equation*}
Since $\mathcal{S}_1$ is convex, $v$ is affine on that set. Therefore
\begin{equation*}
    v(z)=tv(x)+(1-t)v(y).
\end{equation*}
Let $w\in\mathcal{S}_2$. We have
\begin{equation*}
    t\norm{v(x)-v(w)}^2+(1-t)\norm{v(y)-v(w)}^2=t(1-t)\norm{v(x)-v(y)}^2+\norm{v(z)-v(w)}^2.
\end{equation*}
and
\begin{equation*}
    t\norm{x-w}^2+(1-t)\norm{y-w}^2=t(1-t)\norm{x-y}^2+\norm{z-w}^2.
\end{equation*}
Since $v$ is isometric on $\{x,y\}$ and on $\{z,w\}$, we see that
\begin{equation*}
    t\norm{v(x)-v(w)}^2+(1-t)\norm{v(y)-v(w)}^2= t\norm{x-w}^2+(1-t)\norm{y-w}^2.
\end{equation*}
By the $1$-Lipschitzness of $v$, we infer that $v$ is isometric on $\{x,w\}$ and on $\{y,w\}$. This shows that $x,y\in\mathcal{S}_1\cap\mathcal{S}_2$ and completes the proof that the intersection is a face of $\mathcal{S}_1$. We show analogously that it is also a face of $\mathcal{S}_2$.
\end{proof}

Recall that the map $x\mapsto \mathcal{S}(x)$ that maps a point to a leaf that contains the point is Borel measurable with respect to the Borel $\sigma$-algebra generated by the Wijsman topology, see \cite[Section 5]{Ciosmak2021} and Section \ref{s:refine}.

\subsection{Leaves of the derivatives of \texorpdfstring{$\mathcal{C}^{1,1}(\mathbb{R}^n)$}{C1,1(Rn)} maps}\label{s:derivative}

\begin{lemma}\label{lem:isoderma}
Let $u\in\mathcal{C}^{1,1}(\mathbb{R}^n)$ have $1$-Lipschitz derivative. Let $\mathcal{S}$ be a leaf of $Du$. Then there exist orthogonal, mutually complementing linear subspaces $V_1,V_2$ of the tangent space $V$ to $\mathcal{S}$ such that
  for all $x,x_0\in\mathcal{S}$
  \begin{equation}\label{eqn:um}
      u(x)-u(x_0)-Du(x_0)(x-x_0)=\frac12 \big(\norm{P_{V_1}(x-x_0)}^2-\norm{P_{V_2}(x-x_0)}^2\big).
  \end{equation}
  Moreover, for all $x,y\in\mathcal{S}$
  \begin{equation}\label{eqn:isometrym}
      Du(x)-Du(y)=(P_{V_1}-P_{V_2})(x-y).
  \end{equation}
  Here $P_{V_1},P_{V_2}$ are the orthogonal projections onto $V_1$ and $V_2$ respectively.
  If $x,y\in\mathcal{S}$, then also
  \begin{equation}\label{eqn:square}
      z=x+P_{V_1}(y-x)=y-P_{V_2}(y-x)\in\mathcal{S}.
  \end{equation}
\end{lemma}
\begin{proof}
 Since $\mathcal{S}$ is a leaf, there exists a linear isometry $T\colon V\to\mathbb{R}^n$ such that for all $x,y\in\mathcal{S}$
  \begin{equation*}
      Du(x)-Du(y)=T(x-y).
  \end{equation*}
  We claim that $T(V)\subset V$. Let $x\in\mathrm{int}\mathcal{S}$ and let $y\in\mathcal{S}$. Let $v=y-x$, then $v\in V$. Let $w=Tv=Du(y)-Du(x)$. Let us define
  \begin{equation*}
      a=\frac{v+w}{2},\quad b=\frac{v-w}2, \quad z=x+a=y-b.
  \end{equation*}
  Then
  \begin{equation*}
      u(y)-u(x)-Du(x)(y-x)=\int_0^1\langle T(t(y-x)),y-x\rangle\, dt=\frac 12\langle T(y-x),y-x\rangle
  \end{equation*}
  so
  \begin{equation}\label{eqn:differ}
    u(y)-u(x)=\Big\langle v,\frac{Du(x)+Du(y)}2\Big\rangle.
  \end{equation}
Note that the $1$-Lipschitz condition for the derivative $Du$ implies that
  \begin{equation*}
      u(z)- u(x)-\langle Du(x),z-x\rangle=\int_0^1 \langle Du(x+t(z-x))-Du(x),z-x\rangle\, dt\leq \frac12\norm{z-x}^2,
  \end{equation*}
  and
   \begin{equation*}
      u(y)- u(z)-\langle Du(y),y-z\rangle=\int_0^1 \langle Du(z+t(y-z))-Du(y),y-z\rangle\, dt\leq \frac12\norm{z-y}^2,
  \end{equation*}
However, a direct calculation using (\ref{eqn:differ}) shows that adding up the two above inequalities we  get an equality:
  \begin{equation*}
      u(y)+\langle Du(y),z-y\rangle -u(x)-\langle Du(x),z-x\rangle=\frac12\norm{w}^2= \frac12\norm{z-x}^2+ \frac12\norm{z-y}^2.
  \end{equation*}
Thus, we must have equalities in both of the inequalities.
This is to say, for all $t\in [0,1]$,
\begin{equation*}
Du(x+t(z-x))-Du(x)=t(z-x)\text{ and } Du(z+t(y-z))-Du(y)=(1-t)(y-z).
\end{equation*}
Taking $t=0,1$ in the formulae above shows that
\begin{equation*}
    Du(z)-Du(x)=z-x, Du(z)-Du(y)=y-z,
\end{equation*}
and, in particular, $\norm{Du(z)-Du(x)}=\norm{z-x}$. Since $x\in\mathrm{int}\mathcal{S}$, we infer that $z\in\mathcal{S}$. Thus, $a\in V$, and, consequently, $w\in V$. We have proven the claim that $T(V)\subset V$.

    On $\mathrm{int}\mathcal{S}$ the restriction of $u$ to its affine hull has Hessian $T$, since its relative gradient is $P_V Du$ and $T(V)\subset V$. By the symmetry of the second derivative, we infer that $T\colon V\to V$ is symmetric, and thus diagonalisable. Since it is an isometry onto its range, it can only have as its eigenvalues $1$ and $-1$. Let $V_1$ denote the eigenspace of $T$ corresponding to the eigenvalue $1$ and $V_2$ the eigenspace of $T$ corresponding to the eigenvalue $-1$. We see that (\ref{eqn:isometrym}) holds true. Clearly, these subspaces are orthogonal and mutually complementing. Now $a=P_{V_1}(y-x)$ and $b=P_{V_2}(y-x)$, so $z\in\mathcal{S}$ proves (\ref{eqn:square}) for $x\in\mathrm{int}\mathcal{S}$. The general case follows by continuity and closedness of $\mathcal{S}$.
    Now, (\ref{eqn:um}) is a direct consequence of (\ref{eqn:isometrym}): if $x,x_0\in\mathcal{S}$,  then
    \begin{equation*}
           u(x)-u(x_0)-Du(x_0)(x-x_0)=\int_0^1\langle T(t(x-x_0)),x-x_0\rangle\, dt
    \end{equation*}
    which in turn equals
    \begin{equation*}
        \frac12  \big(\norm{P_{V_1}(x-x_0)}^2-\norm{P_{V_2}(x-x_0)}^2\big).
    \end{equation*}
\end{proof}

\begin{example}
    The maximality of the leaves is essential for Lemma \ref{lem:isoderma} to be valid. Indeed, consider $u\in\mathcal{C}^{1,1}(\mathbb{R}^n)$ defined by
    \begin{equation*}
        u(x_1,x_2,x_3)=x_1x_2\text{ for }x=(x_1,x_2,x_3)\in\mathbb{R}^3.
    \end{equation*}
    Then $Du(x)=(x_2,x_1,0)$ is isometric on the line $\mathbb{R}e_1$, which however is not preserved by the derivative. Note that this line is not a leaf, as $Du$ is isometric on the entire plane $\mathbb{R}^2\times \{0\}$. On the plane the eigenspaces are
    \begin{equation*}
        V_1=\mathbb{R}(e_1+e_2),\quad V_2=\mathbb{R}(e_1-e_2),
    \end{equation*}
    where $(e_1,e_2,e_3)$ is the standard orthonormal basis of $\mathbb{R}^3$.
\end{example}

The following corollary shows that leaves of $1$-Lipschitz derivatives enjoy a strong convexity property.

\begin{corollary}\label{col:square}
    Let $\mathcal{S}$ be a leaf of the $1$-Lipschitz derivative of $u\in\mathcal{C}^{1,1}(\mathbb{R}^n)$. Let $V=V_1\oplus V_2$ be the orthogonal decomposition of the tangent space $V$ of $\mathcal{S}$ as in Lemma \ref{lem:isoderma}. Then $\mathcal{S}$ is $(V_1,V_2)$-convex, that is, whenever $x,y\in\mathcal{S}$, then the box
    \begin{equation*}
      \{x+tP_{V_1}(y-x)+sP_{V_2}(y-x)\mid s,t\in [0,1]\}
    \end{equation*}
    is contained in $\mathcal{S}$.
\end{corollary}
\begin{proof}
    By Lemma \ref{lem:isoderma}, we see that $x+P_{V_1}(y-x)\in\mathcal{S}$. Replacing the role of $x,y$ we see that $y+P_{V_1}(x-y)\in\mathcal{S}$. The convex hull of the points $x,y,x+P_{V_1}(y-x),y+P_{V_1}(x-y)$ is the box.  Convexity of leaf $\mathcal{S}$ shows that that box is contained in $\mathcal{S}$.
\end{proof}

\begin{lemma}\label{lem:three}
    Suppose that $u\in\mathcal{C}^{1,1}(\mathbb{R}^n)$ has $1$-Lipschitz derivative. Suppose also that
    \begin{equation}\label{eqn:three}
\big(u(y)+Du(y)(z-y)\big)-\big(u(x)+Du(x)(z-x)\big) = \frac12 \Big(\norm{y-z}^2+\norm{x-z}^2\Big).
\end{equation}
for some $x,y,z\in\mathbb{R}^n$. Then
\begin{equation*}
    \norm{Du(y)-Du(z)}=\norm{y-z}\text{ and }\norm{Du(x)-Du(z)}=\norm{x-z}.
\end{equation*}
More specifically,
\begin{equation*}
   u(y)+Du(y)(z-y)-u(z)=\frac12\norm{y-z}^2,\quad Du(y)-Du(z)=z-y,
\end{equation*}
and
\begin{equation*}
     -u(x)-Du(x)(z-x)+u(z)=\frac12\norm{x-z}^2,\quad Du(x)-Du(z)=x-z.
\end{equation*}
\end{lemma}
\begin{proof}
    Following the lines of \cite[Proposition 2.4]{Gruyer2009} we see that
    \begin{equation*}
           u(x)-u(y)=  \int_0^1 Du(tx+(1-t)y)(x-y)\, dt
    \end{equation*}
    and
    \begin{equation*}
 \abs{u(x)+Du(x)(z-x)-u(z)}=\Big\lvert\int_0^1 \Big(Du(tx+(1-t)z)(x-z)-Du(x)(x-z)\Big)\, dt\Big\rvert,
    \end{equation*}
    which is at most
    \begin{equation*}
     \frac12\norm{z-x}^2.
    \end{equation*}
Now, this shows that
\begin{equation*}
   \big(u(y)+Du(y)(z-y)\big)-\big(u(x)+Du(x)(z-x)\big)\leq   \frac12\norm{z-y}^2+ \frac12\norm{z-x}^2.
\end{equation*}
Our assumption shows that we have equality here. Thus, we must have had equalities
\begin{equation*}
   \norm{ Du(tx+(1-t)z)-Du(x)}=(1-t)\norm{x-z}
\end{equation*}
and
\begin{equation*}
    \norm{ Du(ty+(1-t)z)-Du(y)}=(1-t)\norm{y-z}
\end{equation*}
for all $t\in [0,1]$. Taking $t=0$ yields the desired assertion.  The equality cases in the Cauchy--Schwarz inequality show the four equalities at the end of the formulation of the lemma.
\end{proof}

The above lemma shows that the condition (\ref{eqn:three}) implies that $\{x,z\}$ and $\{y,z\}$ belong to common leaves of $Du$. Choose leaves $\mathcal{S}_1,\mathcal{S}_2$ of $Du$ containing $\{x,z\}$ and $\{y,z\}$, respectively. They may coincide. Let $V(\mathcal{S}_i)$ be the tangent space to  $\mathcal{S}_i$ for $i=1,2$.
Let $\mathcal{F}_{12}$ be the common face of $\mathcal{S}_1$ and of $\mathcal{S}_2$, see Lemma \ref{lem:intersection}.

The following result gives the precise structure of the contact set and will be used repeatedly in the sequel.

\begin{corollary}\label{col:twozet}
     Suppose that $u\in\mathcal{C}^{1,1}(\mathbb{R}^n)$ has $1$-Lipschitz derivative and satisfies (\ref{eqn:three}) for some $x,y,z\in\mathbb{R}^n$. Let  for $i=1,2$, $V(\mathcal{S}_i)=V(\mathcal{S}_i)_1\oplus V(\mathcal{S}_i)_2$ be the orthogonal decompositions of the tangent space $V(\mathcal{S}_i)$ to $\mathcal{S}_i$, as in Lemma \ref{lem:isoderma}.

Then $z$ is the unique point of $\mathcal{F}_{12}$ such that $z-x\in V(\mathcal{S}_1)_1$ and $z-y\in V(\mathcal{S}_2)_2$. Moreover, if $z_0,x_0,y_0\in\mathbb{R}^n$ are some points such that
   \begin{equation*}
\big(u(y_0)+Du(y_0)(z_0-y_0)\big)-\big(u(x_0)+Du(x_0)(z_0-x_0)\big) = \frac12\Big(\norm{y_0-z_0}^2+\norm{x_0-z_0}^2\Big).
\end{equation*}
and $x_0,z_0\in\mathcal{S}_1$ and $y_0,z_0\in\mathcal{S}_2$, then
\begin{equation}\label{eqn:zetzero}
    z=z_0+P_{V(\mathcal{S}_1)_2}(x-x_0)+P_{V(\mathcal{S}_2)_1}(y-y_0)
\end{equation}
and
\begin{equation}\label{eqn:xpp}
   P_{V(\mathcal{S}_2)_1} P_{V(\mathcal{S}_1)_2}(x-x_0)=0,
   \end{equation}
   and
   \begin{equation}\label{eqn:ypp}
  P_{V(\mathcal{S}_1)_2} P_{V(\mathcal{S}_2)_1}(y-y_0)=0.
\end{equation}

Conversely, if $x\in\mathcal{S}_1,y\in\mathcal{S}_2,z\in\mathcal{F}_{12}$ satisfy (\ref{eqn:zetzero}) and (\ref{eqn:xpp}), (\ref{eqn:ypp}), then they satisfy (\ref{eqn:three}).
\end{corollary}
\begin{proof}
According to Lemma \ref{lem:three}
    \begin{equation*}
    \norm{Du(y)-Du(z)}=\norm{y-z}\text{ and }\norm{Du(x)-Du(z)}=\norm{x-z},
\end{equation*}
and
\begin{equation}\label{eqn:yp}
   u(y)+Du(y)(z-y)-u(z)=\frac12\norm{y-z}^2,
\end{equation}
and
\begin{equation}\label{eqn:xp}
     -u(x)-Du(x)(z-x)+u(z)=\frac12\norm{x-z}^2.
\end{equation}
This implies that $y,z\in\mathcal{S}_2$ and $x,z\in\mathcal{S}_1$ belong to common leaves of $Du$.  Let $V_{12}$ denote the tangent space of $\mathcal{F}_{12}$.
Now (\ref{eqn:yp}) and (\ref{eqn:xp}), together with Lemma \ref{lem:isoderma}, show that
\begin{equation}\label{eqn:zett}
    P_{V(\mathcal{S}_2)_1}(z-y)=0\text{, and }P_{V(\mathcal{S}_1)_2}(z-x)=0.
\end{equation}
Note that since $Du$ is well-defined and continuous, we have
\begin{equation}\label{eqn:agree}
    P_{V(\mathcal{S}_1)_1}-P_{V(\mathcal{S}_1)_2}=P_{V(\mathcal{S}_2)_1}-P_{V(\mathcal{S}_2)_2}\text{ on }V_{12}.
\end{equation}
Note also that
 \begin{equation}\label{eqn:orto}
       P_{V(\mathcal{S}_1)_1}+P_{V(\mathcal{S}_1)_2}=P_{V(\mathcal{S}_1)}\text{ and } P_{V(\mathcal{S}_2)_1}+P_{V(\mathcal{S}_2)_2}=P_{V(\mathcal{S}_2)}.
 \end{equation}
By (\ref{eqn:agree}), (\ref{eqn:orto}) we see that
\begin{equation}\label{eqn:projections}
    P_{V(\mathcal{S}_1)_1}=P_{V(\mathcal{S}_2)_1}\text{ and }P_{V(\mathcal{S}_1)_2}=P_{V(\mathcal{S}_2)_2}\text{ on }V_{12}.
\end{equation}
Suppose  that $z'$ is another point on $\mathcal{F}_{12}$ such that (\ref{eqn:zett}) holds true. Then
\begin{equation*}
    z-z'\in V_{12}
\end{equation*}
and thus, by virtue of (\ref{eqn:zett}) and of (\ref{eqn:projections}),
\begin{equation*}
    P_{V(\mathcal{S}_2)_1}(z-z')=P_{V(\mathcal{S}_1)_2}(z-z')=0,\text{ and }P_{V(\mathcal{S}_1)_1}(z-z')=0.
\end{equation*}
Since $z-z'\in V(\mathcal{S}_1)$, we infer that $z-z'=0$. This shows the asserted uniqueness.
The formula (\ref{eqn:zetzero}) follows along similar lines. Since $z-z_0\in V_{12}$ and
\begin{equation*}
    P_{V(\mathcal{S}_1)_2}(z-z_0)=P_{V(\mathcal{S}_1)_2}(x-x_0)\text{ and } P_{V(\mathcal{S}_2)_1}(z-z_0)=P_{V(\mathcal{S}_2)_1}(y-y_0),
\end{equation*}
then, as $P_{V(\mathcal{S}_1)_1}(z-z_0)=P_{V(\mathcal{S}_2)_1}(z-z_0)=P_{V(\mathcal{S}_2)_1}(y-y_0)$, we get
\begin{equation*}
    z-z_0=P_{V(\mathcal{S}_1)_2}(x-x_0)+P_{V(\mathcal{S}_1)_1}(z-z_0)=P_{V(\mathcal{S}_1)_2}(x-x_0)+P_{V(\mathcal{S}_2)_1}(y-y_0).
\end{equation*}
This implies that
\begin{equation*}
   z-x=z_0-x_0-P_{V(\mathcal{S}_1)_1}(x-x_0)+P_{V(\mathcal{S}_2)_1}(y-y_0),
\end{equation*}
and
\begin{equation*}
     z-y=z_0-y_0-P_{V(\mathcal{S}_2)_2}(y-y_0)+P_{V(\mathcal{S}_1)_2}(x-x_0).
\end{equation*}
The fact that $z-x,z_0-x_0\in V(\mathcal{S}_1)_1,z-y,z_0-y_0\in V(\mathcal{S}_2)_2$ shows that
\begin{equation*}
    P_{V(\mathcal{S}_1)_2}P_{V(\mathcal{S}_2)_1}(y-y_0)=0\text{ and }P_{V(\mathcal{S}_2)_1}P_{V(\mathcal{S}_1)_2}(x-x_0)=0.
\end{equation*}
To prove the converse, observe that by (\ref{eqn:xpp}) and (\ref{eqn:ypp}), and the fact that
\begin{equation*}
    z_0-y_0\in V(\mathcal{S}_2)_2, z_0-x_0\in V(\mathcal{S}_1)_1,
\end{equation*}
we get that $z-x\in V(\mathcal{S}_1)_1,z-y\in V(\mathcal{S}_2)_2$.
\end{proof}

\begin{remark}
    If $\mathcal{S}_1=\mathcal{S}_2$, then (\ref{eqn:xpp}) and (\ref{eqn:ypp}) are satisfied trivially. In this case (\ref{eqn:zetzero}) reads $z-z_0=P_{V_2}(x-x_0)+P_{V_1}(y-y_0)$.
\end{remark}

\begin{corollary}\label{col:three}
    Suppose that $u\in\mathcal{C}^{1,1}(\mathbb{R}^n)$ and that (\ref{eqn:three}) holds true for some $x,y,z\in\mathbb{R}^n$ that belong to a common leaf $\mathcal{S}$ of $Du$. Then for every $z_0\in\mathcal{S}$
    \begin{equation}\label{eqn:formu}
        z-z_0=P_{V_2}(x-z_0)+P_{V_1}(y-z_0).
    \end{equation}
\end{corollary}
\begin{proof}
    For $z$ defined by the formula (\ref{eqn:formu}) we see that $z-x\in V_1$ and $z-y\in V_2$.
\end{proof}

\begin{corollary}\label{col:subspaces}
    Suppose that $\mathcal{S}_1,\mathcal{S}_2$ are two leaves of $1$-Lipschitz $Du$, such that $\mathcal{S}_1\cap\mathcal{S}_2\neq\emptyset$, and  $u\in\mathcal{C}^{1,1}(\mathbb{R}^n)$. Let $V(\mathcal{S}_i)=V(\mathcal{S}_i)_1\oplus V(\mathcal{S}_i)_2$, $i=1,2$, be the corresponding orthogonal decompositions of the tangent spaces $V(\mathcal{S}_i)$ to $\mathcal{S}_i$, $i=1,2$, as in Lemma \ref{lem:isoderma}. Suppose that for some $x_0\in\mathcal{S}_1,y_0\in\mathcal{S}_2$ that belong to unique leaves of $Du$ and for some $z_0\in\mathcal{S}_1\cap\mathcal{S}_2$ we have
       \begin{equation*}
\big(u(y_0)+Du(y_0)(z_0-y_0)\big)-\big(u(x_0)+Du(x_0)(z_0-x_0)\big) = \frac12\Big(\norm{y_0-z_0}^2+\norm{x_0-z_0}^2\Big).
\end{equation*}
Then
   \begin{equation*}
 V(\mathcal{S}_1)_2\subset V(\mathcal{S}_2)_2\text{ and }  V(\mathcal{S}_2)_1\subset V(\mathcal{S}_1)_1.
 \end{equation*}
\end{corollary}
\begin{proof}
   By Corollary \ref{col:square}, we see that for $x\in\mathcal{S}_1,y\in\mathcal{S}_2$
   \begin{equation*}
       x'=z_0+P_{V(\mathcal{S}_1)_2}(x-z_0)\in\mathcal{S}_1\text{ and }y'=z_0+P_{V(\mathcal{S}_2)_1}(y-z_0)\in\mathcal{S}_2.
   \end{equation*}
Moreover, $x'-z_0\in V(\mathcal{S}_1)_2$ and $y'-z_0\in V(\mathcal{S}_2)_1$, so that
\begin{equation*}
    Du(x')-Du(z_0)=-(x'-z_0)\text{ and }Du(y')-Du(z_0)=y'-z_0.
\end{equation*}
By Lemma \ref{lem:three}, $Du(x_0)-Du(z_0)=x_0-z_0$ and $Du(y_0)-Du(z_0)=-(y_0-z_0)$. Subtracting we get
\begin{equation*}
    Du(x')-Du(y_0)=y_0-x'\text{ and }Du(y')-Du(x_0)=y'-x_0,
\end{equation*}
so that $Du$ is isometric on $\{y_0,x'\}$ and on $\{x_0,y'\}$. Since $x_0,y_0$ belong to unique leaves of $Du$, we see that $y'\in\mathcal{S}_1,x'\in\mathcal{S}_2$. This implies that
\begin{equation*}
   x'-z_0=P_{V(\mathcal{S}_1)_2}(x-z_0)\in V(\mathcal{S}_2)_2\text{ and }y'-z_0=P_{V(\mathcal{S}_2)_1}(y-z_0)\in V(\mathcal{S}_1)_1.
\end{equation*}
  Note that
  \begin{equation*}
      \mathrm{span}\{x-z_0\mid x\in\mathcal{S}_1\}=V(\mathcal{S}_1)\text{ and }\mathrm{span}\{y-z_0\mid y\in\mathcal{S}_2\}=V(\mathcal{S}_2).
  \end{equation*}
The claim follows readily.
\end{proof}

\begin{remark}
 Note that Corollary \ref{col:twozet}, together with the optimality condition (\ref{eqn:optimality}) of \cite[Theorem 1.1]{Bolbotowski2024}, show that optimal $\sigma\in\Sigma(\mu,\nu)$ is determined by its projection onto the first two copies of $\mathbb{R}^n$, i.e., by a coupling of $\mu$ and $\nu$.
    In \cite[Corollary 1.2]{Bolbotowski2024} it is already proven that for the optimal $\sigma$, its third marginal is uniquely determined by a coupling of $\mu$ and $\nu$. Namely, it is shown that for $\sigma$-almost every $(x,y,z)$ we have
    \begin{equation*}
        z=\frac{x+y}{2}+\frac{Du(y)-Du(x)}2,
    \end{equation*}
    where $u$ is an optimal $\mathcal{C}^{1,1}(\mathbb{R}^n)$ function for $\mathcal{I}(\mu,\nu)$. The corollary of this is the fact that  the third marginal of an optimal $\sigma$ is supported on
    \begin{equation*}
        \bigcup\Big\{B\Big(\frac{x+y}2,\frac{\norm{x-y}}2\Big)\mid x\in\mathrm{supp}\mu,y\in\mathrm{supp}\nu\Big\}
    \end{equation*}
\end{remark}

The theorem below characterises the exact conditions for the leaves $\mathcal{S}_1,\mathcal{S}_2$  of the $1$-Lipschitz derivative of a $\mathcal{C}^{1,1}(\mathbb{R}^n)$ function and their tangent spaces decompositions.

\begin{theorem}\label{thm:compat}
Let $u\in\mathcal C^{1,1}(\mathbb R^n)$ have $1$-Lipschitz
derivative, and let $\mathcal S_1,\mathcal S_2$ be two distinct leaves of
$Du$ intersecting along their common face $\mathcal{F}_{12}$.
For $i=1,2$, let $ V(\mathcal S_i)=V(\mathcal S_i)_1\oplus V(\mathcal S_i)_2$ be the orthogonal decomposition given by Lemma~\ref{lem:isoderma}. Let $z_0\in\mathcal F_{12}$. Then, for every
$x\in\mathcal S_1$, $y\in\mathcal S_2$, and $z\in\mathcal{F}_{12}$
\begin{equation*}
    \left\langle
        P_{V(\mathcal S_1)_1}(x-z_0),
        P_{V(\mathcal S_2)_2}(y-z_0)
    \right\rangle
    \le 0,
\end{equation*}
\begin{equation*}
    \left\langle
        P_{V(\mathcal S_1)_2}(x-z_0),
        P_{V(\mathcal S_2)_1}(y-z_0)
    \right\rangle
    \le 0
\end{equation*}
and
\begin{equation*}
    P_{V(\mathcal S_1)_1}(z-z_0)=  P_{V(\mathcal S_2)_1}(z-z_0),\quad  P_{V(\mathcal S_1)_2}(z-z_0)=  P_{V(\mathcal S_2)_2}(z-z_0).
\end{equation*}
Moreover, for each $x\in\mathrm{int}\mathcal{S}_1$ and $y\in\mathrm{int}\mathcal{S}_2$ at least one of the above inequalities is strict.

Conversely, suppose that the above conditions for the convex sets $\mathcal{S}_1,\mathcal{S}_2$ intersecting along their common face $\mathcal{F}_{12}=\partial\mathcal{S}_1\cap\partial \mathcal{S}_2$, and for their tangent subspaces $ V(\mathcal S_i)=V(\mathcal S_i)_1\oplus V(\mathcal S_i)_2$, $i=1,2$, are satisfied.  For $i=1,2$, set $T_i=P_{V(\mathcal S_i)_1}-P_{V(\mathcal S_i)_2}$. If on each $\mathcal S_i$ we prescribe the quadratic data
\begin{equation*}
 u(s)=\tfrac12\langle T_i(s-z_0),s-z_0\rangle,\qquad Du(s)=T_i(s-z_0)\quad\text{for }s\in\mathcal S_i,
\end{equation*}
then there exists an extension $\tilde{u}\in\mathcal{C}^{1,1}(\mathbb{R}^n)$ of these data with $1$-Lipschitz derivative such that $\mathcal{S}_1,\mathcal{S}_2$ are contained in two distinct leaves of $D\tilde{u}$.
\end{theorem}
\begin{proof}
    Corollary \ref{col:square} shows that
\begin{equation*}
 x'=z_0+P_{V(\mathcal{S}_1)_1}(x-z_0)\in\mathcal{S}_1\quad\text{and}\quad
 y'=z_0+P_{V(\mathcal{S}_2)_2}(y-z_0)\in\mathcal{S}_2
\end{equation*}
and that
\begin{equation*}
 x''=z_0+P_{V(\mathcal{S}_1)_2}(x-z_0)\in\mathcal{S}_1\quad\text{and}\quad
 y''=z_0+P_{V(\mathcal{S}_2)_1}(y-z_0)\in\mathcal{S}_2.
\end{equation*}
The proof of the inequalities is an application of the $1$-Lipschitz condition of $Du$ at points $x',y'$ and at $x'',y''$, respectively, taking into account that the linear part of $Du$ is given by the difference of projectors on each leaf. The proof of the equalities is an application of the fact that $Du$ is well-defined on $\mathcal{F}_{12}$. The fact that $Du$ is not isometric on the union of $\mathcal{S}_1\cup\mathcal{S}_2$ implies that at least one of the inequalities is strict for $x\in\mathrm{int}\mathcal{S}_1,y\in\mathrm{int}\mathcal{S}_2$.

Let us prove the converse.
Let $z_0\in\mathcal{F}_{12}$. Let us define
\begin{equation*}
    u(x)=\frac12 \norm{P_{V(\mathcal{S}_1)_1}(x-z_0)}^2-\frac12\norm{P_{V(\mathcal{S}_1)_2}(x-z_0)}^2\text{ for }x\in\mathcal{S}_1
\end{equation*}
and
\begin{equation*}
    u(y)=\frac12 \norm{P_{V(\mathcal{S}_2)_1}(y-z_0)}^2-\frac12\norm{P_{V(\mathcal{S}_2)_2}(y-z_0)}^2\text{ for }y\in\mathcal{S}_2.
\end{equation*}
These definitions agree on $\mathcal{F}_{12}$. Moreover, the compatibility conditions ensure that the derivatives coincide on $\mathcal{F}_{12}$.

Let $x\in\mathcal{S}_1$ and $y\in\mathcal{S}_2$, with $x\neq y$.
A direct calculation shows that
\begin{equation*}
\frac{2\bigl(u(x)-u(y)\bigr)+\langle Du(x)+Du(y),y-x\rangle}{\norm{x-y}^{2}}
\end{equation*}
equals
\begin{equation*}
\frac{2}{\norm{x-y}^{2}}\Bigl(\left\langle P_{V(\mathcal{S}_1)_1}(x-z_{0}),P_{V(\mathcal{S}_2)_2}(y-z_{0})\right\rangle -\left\langle P_{V(\mathcal{S}_1)_2}(x-z_{0}),P_{V(\mathcal{S}_2)_1}(y-z_{0})\right\rangle\Bigr),
\end{equation*}
and
\begin{equation*}
\frac{\norm{Du(x)-Du(y)}^{2}}{\norm{x-y}^{2}}
\end{equation*}
equals
\begin{equation*}
    \begin{aligned}
 1+\frac{4}{\norm{x-y}^{2}}\Bigl(&\left\langle P_{V(\mathcal{S}_1)_1}(x-z_{0}),P_{V(\mathcal{S}_2)_2}(y-z_{0})\right\rangle\\
 &+\left\langle P_{V(\mathcal{S}_1)_2}(x-z_{0}),P_{V(\mathcal{S}_2)_1}(y-z_{0})\right\rangle\Bigr).
\end{aligned}
\end{equation*}

Substituting the two preceding identities into the formula of
\cite[Proposition~2.2]{Gruyer2009}, we see that the corresponding Le Gruyer constant is at most one if and only if
\begin{equation*}
\left\langle P_{V(\mathcal{S}_1)_1}(x-z_{0}),P_{V(\mathcal{S}_2)_2}(y-z_{0})\right\rangle\leq 0
\end{equation*}
and
\begin{equation*}
\left\langle P_{V(\mathcal{S}_1)_2}(x-z_{0}),P_{V(\mathcal{S}_2)_1}(y-z_{0})\right\rangle\leq 0.
\end{equation*}
Indeed, after substitution, the required inequality is equivalent to
the sum of the two scalar products plus the absolute value of their
difference being non-positive, which is twice their maximum. Note that the bound of the constant for two points in the same leaf is immediate, by \cite[Proposition 2.4]{Gruyer2009}.

Thus, the assumption and \cite[Theorem 1.1]{Gruyer2009} show that there exists $\tilde{u}\in\mathcal{C}^{1,1}(\mathbb{R}^n)$ with $1$-Lipschitz derivative that extends $u$. This $\tilde{u}$ has isometric derivative on $\mathcal{S}_1$ and on $\mathcal{S}_2$, and therefore $\mathcal{S}_1,\mathcal{S}_2$ are contained in two leaves of $D\tilde{u}$. These leaves are necessarily distinct, since at least one of the inequalities for $x\in\mathcal{S}_1,y\in\mathcal{S}_2$ is strict and therefore $D\tilde{u}$ cannot be isometric on $\{x,y\}$.
\end{proof}

\begin{remark}
    The conditions of Corollary \ref{col:subspaces} show that for $x,y,z$ that satisfy (\ref{eqn:three}) the second of the conditions of Theorem \ref{thm:compat} is saturated. Moreover, the conditions of Theorem \ref{thm:compat} are equivalent to obtuseness of the angles between the pairs of tangent cones at $z_0$ to $\mathcal{S}_1$ and to $\mathcal{S}_2$.
\end{remark}

\begin{example}\label{exa:patches}
Let $(e_1,e_2,e_3)$ be the canonical basis of $\mathbb R^3$ and set $w:=\frac{e_2+e_3}{\sqrt2}$.
Consider the two rectangles
\begin{equation*}
    \mathcal{S}_1=\{(x_1,x_2,0)\in\mathbb R^3\mid
|x_1|\le1,\ -1\le x_2\le0\},
\end{equation*}
and
\begin{equation*}
    \mathcal S_2=\{x_1e_1+t w:\ |x_1|\le1,\ 0\le t\le1\}.
\end{equation*}
The two sets meet along the common face
\begin{equation*}
    \mathcal F_{12}
=
\{(x_1,0,0):|x_1|\le1\}.
\end{equation*}
On $\mathcal S_1$ we prescribe
\begin{equation*}
 Du=P_{\mathbb Re_1}-P_{\mathbb Re_2},
\end{equation*}
and
\begin{equation*}
 u(x_1,x_2,0)=\frac12(x_1^2-x_2^2),\qquad |x_1|\le1,\ -1\le x_2\le0 .
\end{equation*}
On $\mathcal S_2$ we prescribe
\begin{equation*}
 Du=P_{\mathbb Re_1}+P_{\mathbb Rw},
\end{equation*}
and
\begin{equation*}
 u(x_1,t/\sqrt2,t/\sqrt2)=\frac12(x_1^2+t^2),\qquad |x_1|\le1,\ 0\le t\le1 .
\end{equation*}
Both prescriptions coincide on $\mathcal F_{12}$ with
\begin{equation*}
 u(x_1,0,0)=\frac12 x_1^2,\qquad Du(x_1,0,0)=x_1e_1 .
\end{equation*}
Moreover, for $x=(x_1,x_2,0)\in \mathcal{S}_1$, $y=(y_1,t/\sqrt2,t/\sqrt2)\in\mathcal{S}_2$, we have
\begin{equation*}
    \norm{Du(x)-Du(y)}^2-\norm{x-y}^2=2\sqrt2 x_2t\leq0,
\end{equation*}
since $x_2\leq 0$, $t\geq 0$. With $z_0=0$, the two scalar products in Theorem~\ref{thm:compat} are $0$ and $x_2t/\sqrt2$, the latter being strictly negative in the relative interiors. The converse in that theorem gives an extension of the prescribed data to a $\mathcal C^{1,1}(\mathbb R^3)$ function with $1$-Lipschitz derivative. The rectangles $\mathcal S_1,\mathcal S_2$ are local patches contained in distinct leaves of its derivative.

\begin{figure}[htbp]
 \centering

\begin{tikzpicture}[
    scale=0.92,
    transform shape,
    x={(1.65cm,0cm)},
    y={(-0.90cm,0.80cm)},
    z={(0cm,3.10cm)},
    >=Latex,
    line cap=round,
    line join=round,
    font=\small,
    leaf edge/.style={
        draw=black!68,
        line width=0.55pt
    },
    interface/.style={
        draw=black,
        line width=1.05pt
    },
    construction/.style={
        draw=black!82,
        line width=0.62pt,
        -{Latex[length=1.8mm,width=1.2mm]}
    },
    point/.style={
        circle,
        fill=black,
        inner sep=1.45pt
    },
    label box/.style={
        text opacity=1,
        inner sep=1.0pt
    }
]

\newcommand{\XYZ}[3]{xyz cs:x=#1,y=#2,z=#3}

\begin{scope}[shift={(0,0)}]

    \coordinate (A) at (\XYZ{-1}{0}{0});
    \coordinate (B) at (\XYZ{ 1}{0}{0});

    \coordinate (C) at (\XYZ{ 1}{-1}{0});
    \coordinate (D) at (\XYZ{-1}{-1}{0});

    \coordinate (E) at (\XYZ{ 1}{0.7071}{0.7071});
    \coordinate (F) at (\XYZ{-1}{0.7071}{0.7071});

    \fill[
        black!36,
        fill opacity=0.12
    ]
        (A)--(B)--(E)--(F)--cycle;

    \draw[leaf edge]
        (A)--(B)--(E)--(F)--cycle;

    \foreach \t in {0.18,0.38,0.58,0.78}{
        \draw[
            black!43,
            line width=0.30pt
        ]
            ($ (A)!\t!(F) $)
            --
            ($ (B)!\t!(E) $);
    }

    \fill[
        black!16,
        fill opacity=0.23
    ]
        (A)--(B)--(C)--(D)--cycle;

    \draw[leaf edge]
        (A)--(B)--(C)--(D)--cycle;

    \draw[interface] (A)--(B);

    \coordinate (X) at (\XYZ{-0.62}{0}{0});
    \coordinate (Z) at (\XYZ{ 0.52}{0}{0});

    \node[point] at (X) {};
    \node[point] at (Z) {};

    \node[
        label box,
        anchor=north east,
        yshift=-1pt,
    ] at (X) {$x$};

    \node[
        label box,
        anchor=south west,
        yshift=5pt,
        xshift=-8pt
    ] at (Z) {$y=z$};

    \draw[construction]
        ($(X)+(0,3pt)$)
        --
        ($(Z)+(0,3pt)$);

    \node[label box]
        at (\XYZ{0.15}{0.56}{0.63})
        {$\mathcal{S}_2$};

    \node[
        label box,
        anchor=north
    ]
        at (\XYZ{0}{-0.085}{0})
        {$\mathcal{F}_{12}$};

    \draw[construction]
        (\XYZ{-0.82}{0.18}{0.18})
        --
        (\XYZ{-0.82}{0.52}{0.52})
        node[
            midway,
            left,
            label box
        ]
        {$w$};

    \node[
        anchor=south
    ] at (-0.6cm,2.8cm)
    {\textup{(a)} \(x\in\mathcal S_1,\ y\in\mathcal S_2\)};

\end{scope}

\begin{scope}[shift={(7.85cm,0)}]

    \coordinate (A) at (\XYZ{-1}{0}{0});
    \coordinate (B) at (\XYZ{ 1}{0}{0});

    \coordinate (C) at (\XYZ{ 1}{-1}{0});
    \coordinate (D) at (\XYZ{-1}{-1}{0});

    \coordinate (E) at (\XYZ{ 1}{0.7071}{0.7071});
    \coordinate (F) at (\XYZ{-1}{0.7071}{0.7071});

    \fill[
        black!36,
        fill opacity=0.12
    ]
        (A)--(B)--(E)--(F)--cycle;

    \draw[leaf edge]
        (A)--(B)--(E)--(F)--cycle;

    \foreach \t in {0.18,0.38,0.58,0.78}{
        \draw[
            black!43,
            line width=0.30pt
        ]
            ($ (A)!\t!(F) $)
            --
            ($ (B)!\t!(E) $);
    }

    \fill[
        black!16,
        fill opacity=0.23
    ]
        (A)--(B)--(C)--(D)--cycle;

    \draw[leaf edge]
        (A)--(B)--(C)--(D)--cycle;

    \draw[interface] (A)--(B);

    \coordinate (X) at (\XYZ{-0.52}{ 0.48}{0.48});
    \coordinate (Z) at (\XYZ{ 0.25}{ 0   }{0});
    \coordinate (Y) at (\XYZ{ 0.25}{-0.68}{0});

    \node[point] at (X) {};
    \node[point] at (Y) {};
    \node[point] at (Z) {};

    \node[
        label box,
        anchor=south west,
        xshift=2pt,
        yshift=-1.5pt,
    ] at (X) {$x'$};

    \node[
        label box,
        anchor=north east,
        xshift=-2.5pt,
        yshift=4pt,
    ] at (Y) {$y'$};

    \node[
        label box,
        anchor=south west,
        xshift=2pt,
        yshift=2pt
    ] at (Z) {$z'$};

 \draw[construction] (X)--(Z);

    \draw[construction]
        (Y)--(Z)
        node[
            pos=0.55,
            right,
            xshift=5pt,
            label box
        ]
        { };

    \node[label box]
        at (\XYZ{0.68}{-0.74}{0})
        {$\mathcal S_1$};

    \node[label box]
        at (\XYZ{0.15}{0.56}{0.63})
        {$\mathcal S_2$};

    \node[
        label box,
        anchor=north
    ]
        at (\XYZ{-0.48}{-0.085}{0})
        {$\mathcal F_{12}$};

    \node[
        anchor=south
    ] at (-0.6cm,2.8cm)
    {\textup{(b)} \(x'\in\mathcal S_2,\ y'\in\mathcal S_1\)};

\end{scope}

\end{tikzpicture}

\caption{}\label{fig:two-equality-configurations}
\end{figure}

For $p,z$ in the same rectangle the prescribed data give
\begin{equation*}
 2\big(u(z)-u(p)-Du(p)(z-p)\big)=
 \begin{cases}
 (z_1-p_1)^2-(z_2-p_2)^2,& p,z\in\mathcal S_1,\\
 \norm{z-p}^2,& p,z\in\mathcal S_2.
 \end{cases}
\end{equation*}
For a triple $x\in\mathcal{S}_1$, $y\in\mathcal{S}_2$, $z\in\mathcal{F}_{12}$, substitution in (\ref{eqn:three}) shows that equality holds precisely when
\begin{equation*}
    y=z\text{ and }z-x=se_1\text{ for some }s\in\mathbb{R}.
\end{equation*}
If $x'\in\mathcal{S}_2$, $y'\in\mathcal{S}_1$, $z'\in\mathcal{F}_{12}$, the same calculation gives
\begin{equation*}
    z'-x'\in\mathrm{span}\{e_1,w\},\qquad z'-y'\in\mathbb{R}e_2.
\end{equation*}
Thus, there are no restrictions for $x'\in\mathcal{S}_2$, and for $y'$ we have $y'=z'+re_2$ for some $r\in\mathbb{R}$. These two situations are illustrated by Figure \ref{fig:two-equality-configurations}.

If $x'',y'',z''\in\mathcal{S}_1$, the same formula shows that (\ref{eqn:three}) holds precisely when
\begin{equation*}
    z''-z_0=P_{\mathbb{R}e_2}(x''-z_0)+P_{\mathbb{R}e_1}(y''-z_0).
\end{equation*}
This is presented in Figure \ref{fig:one-leaf}.

\begin{figure}[h]
\centering

\begin{tikzpicture}[
    x={(4.00cm,0cm)},
    y={(1.42cm,-2.75cm)},
    >=Latex,
    line cap=round,
    line join=round,
    font=\small,
    leaf edge/.style={
        draw=black!72,
        line width=0.55pt
    },
    interface/.style={
        draw=black,
        line width=1.00pt
    },
    vector/.style={
        draw=black!88,
        line width=0.68pt,
        -{Latex[length=1.9mm,width=1.2mm]}
    },
    guide/.style={
        draw=black!48,
        line width=0.42pt,
        dash pattern=on 2.4pt off 1.8pt
    },
    projection guide/.style={
        draw=black!42,
        line width=0.36pt,
        densely dotted
    },
    point/.style={
        circle,
        fill=black,
        inner sep=1.45pt
    },
    projection point/.style={
        circle,
        draw=black,
        fill=black,
        line width=0.55pt,
        inner sep=1.30pt
    },
    label box/.style={
        text opacity=1,
        inner sep=0.9pt
    }
]

\coordinate (A) at (-1,0);
\coordinate (B) at ( 1,0);
\coordinate (C) at ( 1,1);
\coordinate (D) at (-1,1);
\coordinate (O) at ( 0,0);

\fill[
    black!18,
    fill opacity=0.20
]
    (A)--(B)--(C)--(D)--cycle;

\draw[leaf edge]
    (A)--(B)--(C)--(D)--cycle;

\draw[interface]
    (A)--(B);

\coordinate (X)  at (-0.70,0.76);
\coordinate (Y)  at ( 0.58,0.30);
\coordinate (PX) at ( 0,0.76);
\coordinate (PY) at ( 0.58,0);
\coordinate (Z)  at ( 0.58,0.76);

\draw[projection guide]
    (X)--(PX);

\draw[projection guide]
    (Y)--(PY);

\draw[vector]
    (O)--(PX)
    node[
        pos=0.48,
        left,
        xshift=-3pt,
        label box
    ]
    {$P_{\mathbb{R}e_2}(x''-z_0)$};

\draw[vector]
    ($(O)+(0,-2.8pt)$)
    --
    ($(PY)+(0,-2.8pt)$)
    node[
        pos=0.47,
        above,
        yshift=4pt,
        label box
    ]
    {$P_{\mathbb{R}e_1}(y''-z_0)$};

\draw[guide]
    (PX)--(Z);

\draw[guide]
    (PY)--(Z);

\node[point] at (O) {};
\node[point] at (X) {};
\node[point] at (Y) {};
\node[point] at (Z) {};

\node[projection point] at (PX) {};
\node[projection point] at (PY) {};

\node[
    label box,
    anchor=north east,
    xshift=-1pt,
    yshift=-1pt
] at (O) {$z_0$};

\node[
    label box,
    anchor=north east,
    xshift=-1pt,
    yshift=-1pt
] at (X) {$x''$};

\node[
    label box,
    anchor=south west,
    xshift=1pt,
    yshift=1pt
] at (Y) {$y''$};

\node[
    label box,
    anchor=north west,
    xshift=2pt,
    yshift=-1pt
] at (Z) {$z''$};

\node[label box]
    at (0.80,0.72)
    {$\mathcal S_1$};

\node[
    label box,
    anchor=north
] at (-0.62,0.05)
    {$\mathcal F_{12}$};

 \end{tikzpicture}

\caption{}\label{fig:one-leaf}
\end{figure}
\end{example}

\subsection{Notation}\label{s:notation}
We specify the notation used in the introduction and recall the martingale condition. By $\mathcal{P}(\mathbb{R}^n)$ we shall denote the set of all Borel probability measures on $\mathbb{R}^n$. For $p\in [1,\infty)$ by  $\mathcal{P}_p(\mathbb{R}^n)$ we denote the set of all  probabilities $\mu\in\mathcal{P}(\mathbb{R}^n)$ such that
\begin{equation*}
    \int_{\mathbb{R}^n}\norm{x}^p\, d\mu(x)<\infty.
\end{equation*}

If $\mu\in\mathcal{P}(\mathbb{R}^n)$ has finite first moments, i.e., if $\mu\in\mathcal{P}_1(\mathbb{R}^n)$, then the point
\begin{equation*}
    \int_{\mathbb{R}^n}x\, d\mu(x)
\end{equation*}
we call the \emph{barycentre} of $\mu$.

If $\sigma\in\mathcal{P}((\mathbb{R}^n)^3)$, then by $\mathrm{P}_i\sigma$ we denote the pushforward of $\sigma$ via the projection on the $i$\textsuperscript{th} coordinate, for $i=1,2,3$. The measure $\mathrm{P}_i\sigma$ is also called the $i$\textsuperscript{th} marginal of $\sigma$. By $\mathrm{P}_{ij}\sigma$ we denote the pushforward of $\sigma$ on the $ij$\textsuperscript{th} coordinates for $i\neq j$, $i,j=1,2,3$. We call $\mathrm{P}_{ij}\sigma$ the $ij$\textsuperscript{th} marginal of $\sigma$. Of course, similar definitions pertain to the case of measures in $\mathcal{P}((\mathbb{R}^n)^2)$.

For a Borel measure $\mu\in\mathcal{P}(X)$ and a Borel set $A\subset X$ such that $\mu(A)>0$ we define $\mu|_A\in\mathcal{P}(X)$ via the formula
\begin{equation*}
    \mu|_A(B)=\frac{\mu(A\cap B)}{\mu(A)}\text{ for every Borel }B\subset X.
\end{equation*}

A Borel probability measure $\pi\in\mathcal{P}(\mathbb{R}^n\times\mathbb{R}^n)$ is called a \emph{martingale transport} whenever it  is a distribution of a one-step martingale $(X_1,X_2)$ on $\mathbb{R}^n$, i.e., $(X_1,X_2)\sim\pi$. This is equivalent to assuming that the marginals $\mathrm{P}_i\pi\in\mathcal{P}_1(\mathbb{R}^n)$ for $i=1,2$ and that for all bounded and Borel measurable $g\colon\mathbb{R}^n\to\mathbb{R}$
\begin{equation*}
    \int_{\mathbb{R}^n\times\mathbb{R}^n} (y-x) g(x) \, d\pi(x,y)= 0 .
\end{equation*}

If for $\pi\in\mathcal{P}(\mathbb{R}^n\times\mathbb{R}^n)$ its first and second marginal are equal $\mu\in\mathcal{P}(\mathbb{R}^n)$ and $\nu\in\mathcal{P}(\mathbb{R}^n)$ respectively, then we write
\begin{equation*}
    \pi\in\Gamma(\mu,\nu).
\end{equation*}

\subsection{Kantorovich--Rubinstein duality for the Hessian}\label{s:duality}

Let us now recall the definitions and results of \cite{Bolbotowski2024}.
\begin{definition}
    For $\mu,\nu\in\mathcal{P}_1(\mathbb{R}^n)$ with common barycentre, by $\Sigma(\mu,\nu)$  we shall denote the set of Borel probability measures $\sigma\in \mathcal{P}(\mathbb{R}^n\times\mathbb{R}^n\times \mathbb{R}^n)$ such that
    \begin{equation*}
        \mathrm{P}_1\sigma=\mu,\mathrm{P}_2\sigma=\nu\text{, and }\mathrm{P}_3\sigma\in\mathcal{P}_1(\mathbb{R}^n),
    \end{equation*}
    and such that for all bounded and Borel measurable $g,h\colon\mathbb{R}^n\to\mathbb{R}$
\begin{equation}\label{eqn:mart}
    \int_{(\mathbb{R}^n)^3} (z-x) g(x) \, d\sigma(x,y,z)=   \int_{(\mathbb{R}^n)^3} (z-y)h(y) \, d\sigma(x,y,z)=0,
\end{equation}
i.e., the marginals $\mathrm{P}_{13}\sigma,\mathrm{P}_{23}\sigma$ are martingale transports.
\end{definition}

In \cite{Bolbotowski2024} the following three-marginal optimal transport problem is studied: for two Borel probability measures $\mu,\nu\in\mathcal{P}_2(\mathbb{R}^n)$ with common barycentre, we look for $\sigma\in\Sigma(\mu,\nu)$ for  which the integral
\begin{equation*}
    \int_{(\mathbb{R}^n)^3} \frac12 (\norm{z-x}^2+\norm{z-y}^2)\, d\sigma(x,y,z)
\end{equation*}
is minimal. The minimal value is denoted by
\begin{equation}\label{eqn:j}
    \mathcal{J}(\mu,\nu)=\inf\Big\{\int_{(\mathbb{R}^n)^3} \frac12 (\norm{z-x}^2+\norm{z-y}^2)\, d\sigma(x,y,z)\mid \sigma\in\Sigma(\mu,\nu)\Big\}.
\end{equation}
In \cite[Theorem 1.1]{Bolbotowski2024} it is shown that an optimal $\sigma\in\Sigma(\mu,\nu)$ exists. Moreover, $\mathcal{J}(\mu,\nu)$ coincides with
\begin{equation}\label{eqn:i}
    \mathcal{I}(\mu,\nu)=\sup\Big\{\int_{\mathbb{R}^n}v\, d(\nu-\mu)\mid v\in\mathcal{C}^{1,1}(\mathbb{R}^n)\text{ has }1\text{-Lipschitz derivative}\Big\}.
\end{equation}
If $u\in\mathcal{C}^{1,1}(\mathbb{R}^n)$ has $1$-Lipschitz derivative and attains (\ref{eqn:i}), then by \cite[Theorem 1.1., (ii)]{Bolbotowski2024} for $\sigma$-almost every $(x,y,z)\in(\mathbb{R}^n)^3$
\begin{equation*}
    u(y)+Du(y)(z-y)-(u(x)+Du(x)(z-x))=\frac12(\norm{z-y}^2+\norm{z-x}^2).
\end{equation*}

Following \cite{Bolbotowski2024} we define $c\colon\mathbb{R}^n\times\mathbb{R}^n\times\mathbb{R}^n\to\mathbb{R}$ by the formula
\begin{equation*}
    c(x,y,z)=\frac12(\norm{z-y}^2+\norm{z-x}^2)\text{ for }x,y,z\in\mathbb{R}^n.
\end{equation*}

We note that the assumptions that $\mu,\nu$ have finite second moments and that they have common barycentre guarantee that $\mathcal{J}(\mu,\nu)$ and $\mathcal{I}(\mu,\nu)$ are finite.

\subsection{Mass balance and moment balance conditions}\label{s:balance}

The proof of the theorem in this section is based on the idea of \cite[Theorem 4]{Ciosmak20211}. The theorem is an important ingredient in proving that the problem (\ref{eqn:j}) for a pair of two probability measures can be completely reduced to the collection of such problems on the leaves of $Du$ for an optimal $u\in\mathcal{C}^{1,1}(\mathbb{R}^n)$ with $1$-Lipschitz derivative.

\begin{theorem}\label{thm:deco}
    Suppose that $\mu,\nu\ll\lambda$ are two Borel probabilities in $\mathcal{P}_2(\mathbb{R}^n)$ with common barycentre, absolutely continuous with respect to the Lebesgue measure $\lambda$.
    Let $u\in\mathcal{C}^{1,1}(\mathbb{R}^n)$ be a function with $1$-Lipschitz derivative that attains the supremum
    \begin{equation*}
    \mathcal{I}(\mu,\nu)=\sup\Big\{\int_{\mathbb{R}^n}v\, d(\nu-\mu)\mid v\in\mathcal{C}^{1,1}(\mathbb{R}^n)\text{ has }1\text{-Lipschitz derivative}\Big\}
\end{equation*}
and let $\sigma\in \Sigma(\mu,\nu)$ attain the infimum
\begin{equation*}
    \mathcal{J}(\mu,\nu)=\inf\Big\{\int_{(\mathbb{R}^n)^3} \frac12 (\norm{z-x}^2+\norm{z-y}^2)\, d\sigma(x,y,z)\mid \sigma\in\Sigma(\mu,\nu)\Big\}.
\end{equation*}
    Let $A\subset\mathbb{R}^n$ be a transport set of $Du$.
    Then
    \begin{equation*}
       \abs{ \mu(A)-\nu(A)}\leq 2\mathrm{P}_3\sigma(N(Du))
       \end{equation*}
       and
       \begin{equation*}
       \Big\lVert \int_A x\, d\mu(x)-\int_A y\, d\nu(y)\Big\rVert\leq2\int_{N(Du)} \norm{x}\, d \mathrm{P}_3\sigma(x).
    \end{equation*}
    In particular, if $\sigma\in\Sigma(\mu,\nu)$ has absolutely continuous third marginal, then the mass balance and the moment balance conditions hold:
      \begin{equation*}
        \mu(A)=\nu(A)\text{ and }\int_A x\, d\mu(x)=\int_A y\, d\nu(y).
    \end{equation*}
    Moreover, in this case, if $\mu(A),\nu(A)>0$, then $u$ and $\sigma|_{A^3}$  are also optimal for the problems $\mathcal{I}(\mu|_A,\nu|_A)$ and $\mathcal{J}(\mu|_A,\nu|_A)$ for the pair of measures $\mu|_A,\nu|_A\in\mathcal{P}_2(\mathbb{R}^n)$.

\end{theorem}

\begin{remark}
 Let us remark that by \cite[Theorem 1.3]{Bolbotowski2024} a plan $\sigma\in\Sigma(\mu,\nu)$ is optimal for $\mathcal{J}(\mu,\nu)$ if and  only if its third marginal $ \mathrm{P}_3\sigma$  is  optimal for the problem
\begin{equation*}
    \mathcal{V}(\mu,\nu)=\inf\Big\{\int_{\mathbb{R}^n}\norm{x-x_0}^2\,d\rho(x)\mid \rho\in\mathcal{P}_2(\mathbb{R}^n), \mu\prec_c\rho,\nu\prec_c\rho\Big\},
\end{equation*}
where  $x_0\in\mathbb{R}^n$ is the common barycentre of $\mu,\nu,\rho$.
Thus the assumption that there exists an optimal $\sigma\in\Sigma(\mu,\nu)$ with absolutely continuous third marginal is equivalent to the existence of absolutely continuous minimiser of $\mathcal{V}(\mu,\nu)$.
\end{remark}

\begin{proof}[Proof of Theorem \ref{thm:deco}]
Let $\sigma\in\Sigma(\mu,\nu)$, $u\in\mathcal{C}^{1,1}(\mathbb{R}^n)$ be as  in the assumptions of the theorem.
By \cite[Theorem 1.1, (ii)]{Bolbotowski2024} for $\sigma$-almost every $(x,y,z)$ we have
    \begin{equation}\label{eqn:bolbo}
   \big(u(y)+Du(y)(z-y)\big)-\big(u(x)+Du(x)(z-x)\big)=\frac12(\norm{y-z}^2+\norm{x-z}^2)
\end{equation}
Note that $B(Du)$ -- the set of points that belong to at least two distinct leaves of $Du$ -- is contained in the non-differentiability set $N(Du)$ of $Du$, and thus is of the Lebesgue measure zero by \cite[Corollary 2.15]{Ciosmak2021}.
Fix a triple in the $\sigma$-full set where \eqref{eqn:bolbo} holds. By Lemma~\ref{lem:three}, $Du$ is isometric on $\{x,z\}$ and on $\{y,z\}$. Each pair is contained in a leaf by the Kuratowski--Zorn lemma, since the union of a chain of isometric sets is isometric. If $x\in A\setminus B(Du)$, its unique leaf is one of the leaves making up $A$ and contains $z$. The same argument with $x,z$ exchanged applies when $z\in A\setminus B(Du)$. Hence
\begin{equation*}
       \mathbf{1}_{A\setminus B(Du)}(x)\leq    \mathbf{1}_{A}(z)\leq    \mathbf{1}_{A}(x)+   \mathbf{1}_{A\cap B(Du)}(z).
\end{equation*}
This shows that
\begin{equation}\label{eqn:x}
      \abs{ \mathbf{1}_{A\setminus B(Du)}(x)-   \mathbf{1}_{A\setminus B(Du)}(z)}\leq  \mathbf{1}_{B(Du)}(z)+ \mathbf{1}_{B(Du)}(x),
      \end{equation}
      and similarly
      \begin{equation}\label{eqn:y}
      \abs{ \mathbf{1}_{A\setminus B(Du)}(y)-   \mathbf{1}_{A\setminus B(Du)}(z)}\leq  \mathbf{1}_{B(Du)}(z)+ \mathbf{1}_{B(Du)}(y).
\end{equation}
for $\sigma$-almost every $(x,y,z)\in(\mathbb{R}^n)^3$.

 Now, the first two marginals of $\sigma$ are absolutely continuous. We have
\begin{align*}
   & \mu(A)-\nu(A)=\int_{(\mathbb{R}^n)^3}\Big(\mathbf{1}_A(x)-\mathbf{1}_{A}(y)\Big)\, d\sigma(x,y,z)=\\
   &\int_{(\mathbb{R}^n)^3}\Big(\big(\mathbf{1}_{A\setminus B(Du)}(x)-\mathbf{1}_{A\setminus B(Du)}(z)\big)-\big(\mathbf{1}_{A\setminus B(Du)}(y)-\mathbf{1}_{A\setminus B(Du)}(z)\big)\Big)\, d\sigma(x,y,z).
\end{align*}
Therefore the first bound asserted in the theorem follows by absolute continuity of $\mu,\nu$, by (\ref{eqn:x}), and by (\ref{eqn:y}).
To prove the bound for the moments we proceed analogously. That is, we observe that
\begin{align*}
 &   \int_A x\, d\mu(x)-\int_A y\, d\nu(y)=\int_{(\mathbb{R}^n)^3} \big(\mathbf{1}_A(x)x-\mathbf{1}_A(y)y\big)\, d\sigma(x,y,z)=\\
 & \int_{(\mathbb{R}^n)^3} \Big(\big(\mathbf{1}_{A\setminus B(Du)}(x)x-\mathbf{1}_{A\setminus B(Du)}(z)z\big)\phantom{\Big)}\\
 &\hspace{6em}{}-\big(\mathbf{1}_{A\setminus B(Du)}(y)y-\mathbf{1}_{A\setminus B(Du)}(z)z\big)\Big)\, d\sigma(x,y,z).
\end{align*}
Note that by (\ref{eqn:mart})  we have
\begin{equation*}
    \int_{(\mathbb{R}^n)^3} \Big(\mathbf{1}_{A\setminus B(Du)}(x)(x-z)-\mathbf{1}_{A\setminus B(Du)}(y)(y-z)\Big)\, d\sigma(x,y,z)=0.
\end{equation*}
Thus the bound for the moments follows.
Suppose now that $\sigma\in\Sigma(\mu,\nu)$ has absolutely continuous third marginal and that $\mu(A),\nu(A)>0$. Then the just proven bounds show that  the barycentres and the masses of $\mu|_A$ and of $\nu|_A$ agree. Therefore the problem (\ref{eqn:j}) for normalised  measures $\mu|_A$ and $\nu|_A$ admits a solution.
By (\ref{eqn:x}) and (\ref{eqn:y}), we see that
\begin{equation*}
    \mathbf{1}_A(x)=\mathbf{1}_A(y)=\mathbf{1}_A(z)\text{ for }\sigma\text{-almost every }(x,y,z)\in(\mathbb{R}^n)^3.
\end{equation*}
Thus
\begin{equation*}
    \sigma|_{A^3}\in \Sigma (\mu|_A,\nu|_A).
\end{equation*}
Indeed, for Borel $B\subset\mathbb{R}^n$
\begin{equation*}
   \sigma((A\cap B) \times A^2)=   \sigma((A\cap B) \times (\mathbb{R}^n)^2)=  \mu(A\cap B).
      \end{equation*}
  Similarly,
  \begin{equation*}
   \sigma(A\times (A\cap B) \times A)= \sigma(\mathbb{R}^n\times (A\cap B)\times\mathbb{R}^n)=\nu(A\cap B).
   \end{equation*}
The martingale condition for $\sigma|_{A^3}$ is a direct consequence of the condition for $\sigma$.

Thus the optimality claim follows by \cite[Theorem 1.1, (ii)]{Bolbotowski2024}.

\end{proof}

\subsection{Disintegration}\label{s:disintegrate}

Let us recall a theorem that follows readily from \cite[Example 10.4.11, Definition 10.4.2]{Bogachev20072}.

\begin{theorem}\label{thm:disintegration}
Let $X,Y$ be two Polish spaces. Let $p\colon X\to Y$ be a Borel map and let $\rho$ be a non-negative finite Borel measure on $X$. Let $\theta$ be the push-forward of measure $\rho$ via $p$. Then there exist Borel measures $(\rho_y)_{y\in Y}$ on $X$ such that
\begin{enumerate}[label=(\roman*)]
\item  for every Borel set $B\subset X$ the function $y\mapsto\rho_y(B)$ is Borel measurable,
\item\label{i:disintegrationfibre} for $\theta$-almost every $y\in p(X)$ the measure $\rho_y$ is concentrated on $p^{-1}(y)$,
\item\label{i:pro} for every Borel sets $B\subset X$ and $E\subset Y$ there is
\begin{equation*}
\rho(B\cap p^{-1}(E))=\int_E \rho_y(B) d\theta(y).
\end{equation*}
\end{enumerate}
\end{theorem}

\begin{theorem}\label{thm:partim}
Let $\mu,\nu\ll\lambda$ be two Borel probability measures in $\mathcal{P}_2(\mathbb{R}^n)$. Let $\sigma\in\Sigma(\mu,\nu)$ be optimal for $\mathcal{J}(\mu,\nu)$ and $u\in\mathcal{C}^{1,1}(\mathbb{R}^n)$ be optimal for $\mathcal{I}(\mu,\nu)$. Suppose that $\sigma\in\Sigma(\mu,\nu)$ has absolutely continuous third marginal. There exists a Borel measure $\theta$ on $CC(\mathbb{R}^n)$ supported on the leaves of $Du$, and Borel measures $\mu_{\mathcal{S}},\nu_{\mathcal{S}}$, such that
\begin{enumerate}[label=(\roman*)]
    \item\label{i:mees} for any Borel set $A\subset\mathbb{R}^n$ the maps
    \begin{equation*}
        CC(\mathbb{R}^n)\ni \mathcal{S}\mapsto \mu_{\mathcal{S}}(A)\in\mathbb{R}\text{, and }CC(\mathbb{R}^n)\ni \mathcal{S}\mapsto \nu_{\mathcal{S}}(A)\in\mathbb{R}
    \end{equation*}
    are  Borel measurable,
    \item for $\theta$-almost every $\mathcal{S}$, $\mu_{\mathcal{S}}$, $\nu_{\mathcal{S}}$ are concentrated on $\mathcal{S}$, i.e.,
    \begin{equation*}
    \mu_{\mathcal{S}}(\mathcal{S}^c)=0,    \nu_{\mathcal{S}}(\mathcal{S}^c)=0,
    \end{equation*}
    \item for every Borel set $A\subset\mathbb{R}^n$
    \begin{equation*}
        \mu(A)=\int_{CC(\mathbb{R}^n)}\mu_{\mathcal{S}}(A)\, d\theta(\mathcal{S}), \nu(A)=\int_{CC(\mathbb{R}^n)}\nu_{\mathcal{S}}(A)\, d\theta(\mathcal{S}).
    \end{equation*}
\end{enumerate}
\end{theorem}
\begin{proof}
  The theorem follows from an application of Theorem \ref{thm:disintegration}.  The only non-trivial fact to check is that the pushforward $\theta\in\mathcal{P}(CC(\mathbb{R}^n))$ of $\mu$ with respect to the map $\mathcal{S}$ coincides with the pushforward of $\nu$ with respect to the map $\mathcal{S}$. As in the proof of Theorem \ref{thm:deco}, $x,z$ and $y,z$ belong pairwise to common leaves for $\sigma$-almost every $(x,y,z)$. All three marginals are absolutely continuous, so the leaves at $x,y,z$ are unique there and $\mathcal S(x)=\mathcal S(z)=\mathcal S(y)$ outside a $\sigma$-null set. Thus for any Borel set $E\subset CC(\mathbb R^n)$ and $A=\mathcal S^{-1}(E)$ we have
\begin{equation*}
    \mu(A)=\nu(A).
\end{equation*}
This is however exactly what we were claiming.
\end{proof}

For any optimal plan $\sigma\in\Sigma(\mu,\nu)$, we define probability measures $\sigma_{\mathcal{S}}\in\mathcal{P}((\mathbb{R}^n)^3)$ as the conditional measures of $\sigma$ obtained through Theorem \ref{thm:disintegration} and the map
\begin{equation*}
    (\mathbb{R}^n)^3 \ni (x_1,x_2,x_3)\mapsto \mathcal{S}(x_1)\in CC(\mathbb{R}^n).
\end{equation*}
Since $x\mapsto\mathcal{S}(x)$ is Borel measurable, so is the above defined map.
Note that the pushforward of $\sigma$ via this map coincides with the pushforward $\theta$ of $\mu$ via the map  $x\mapsto \mathcal{S}(x)$, as the first marginal of $\sigma$ is $\mu$. The conditional measures $\sigma_{\mathcal{S}}$ enjoy the following properties:
\begin{enumerate}[label=(\roman*)]
    \item for any Borel set $B\subset(\mathbb{R}^n)^3$ the map
    \begin{equation*}
        CC(\mathbb{R}^n)\ni \mathcal{S}\mapsto\sigma_{\mathcal{S}}(B)\text{ is Borel measurable,}
    \end{equation*}
    and
    \begin{equation*}
        \sigma(B)=\int_{CC(\mathbb{R}^n)}\sigma_{\mathcal{S}}(B)\, d\theta(\mathcal{S}),
    \end{equation*}
    \item for $\theta$-almost every $\mathcal{S}$ the measure $\sigma_{\mathcal{S}}$ is concentrated on $\mathcal{S}\times(\mathbb{R}^n)^2$.
\end{enumerate}

If we disintegrate $\mu,\nu$ with respect to the partition induced by the leaves of an optimal $u$, then $u$ will still be an optimal map for almost every such pair of conditional measures. I.e., for almost every leaf $\mathcal{S}$ of $Du$, $Du$ is isometric on the leaf and $u$ is an optimal map for the problem $\mathcal{I}(\mu_{\mathcal{S}},\nu_{\mathcal{S}})$ for the conditional measures. This is the essence of the following theorem. Moreover, if $\sigma\in\Sigma(\mu,\nu)$ is optimal for $\mathcal{J}(\mu,\nu)$ and some optimal plan has absolutely continuous third marginal, then the conditional measures $\sigma_{\mathcal{S}}\in\Sigma(\mu_{\mathcal{S}},\nu_{\mathcal{S}})$ will be optimal for $\mathcal{J}(\mu_{\mathcal{S}},\nu_{\mathcal{S}})$.

\begin{corollary}\label{col:oneleaf}
    Let $\mu,\nu\ll\lambda$ be two Borel probability measures in $\mathcal{P}_2(\mathbb{R}^n)$. Let $\sigma\in\Sigma(\mu,\nu)$ be optimal for $\mathcal{J}(\mu,\nu)$ and $u\in\mathcal{C}^{1,1}(\mathbb{R}^n)$ be optimal for $\mathcal{I}(\mu,\nu)$. Suppose that $\mathcal{S}_{\#}\mu=\mathcal{S}_{\#}\nu$.
    Then for $\sigma$-almost every $(x,y,z)\in(\mathbb{R}^n)^3$, $x,y,z$ belong to a common leaf of $Du$.
\end{corollary}
\begin{proof}
    Let $\pi=\mathrm{P}_{12}\sigma$. By Corollary \ref{col:subspaces}, we see that for $\pi$-almost every $(x,y)$ we have
    \begin{equation*}
        V(\mathcal{S}(x))_2\subset V(\mathcal{S}(y))_2\text{ and }V(\mathcal{S}(y))_1\subset V(\mathcal{S}(x))_1.
    \end{equation*}
    Therefore
    \begin{equation}\label{eqn:dimensions}
          \mathrm{dim}V(\mathcal{S}(x))_2\leq\mathrm{dim} V(\mathcal{S}(y))_2\text{ and }\mathrm{dim}V(\mathcal{S}(y))_1\leq \mathrm{dim} V(\mathcal{S}(x))_1.
    \end{equation}
    By Lemma \ref{lem:v} these functions are all measurable.
    Integrating the inequalities (\ref{eqn:dimensions}) against $\pi$, we get
    \begin{equation*}
       \int_{\mathbb{R}^n} \mathrm{dim}(V(\mathcal{S}(x))_2)\, d\mu(x)\leq\int_{\mathbb{R}^n}\mathrm{dim}(V(\mathcal{S}(y))_2)\, d\nu(y),
    \end{equation*}
    and
        \begin{equation*}
       \int_{\mathbb{R}^n} \mathrm{dim}(V(\mathcal{S}(y))_1)\, d\nu(y)\leq\int_{\mathbb{R}^n}\mathrm{dim}(V(\mathcal{S}(x))_1)\, d\mu(x),
    \end{equation*}
    Since $\mathcal{S}_{\#}\mu=\mathcal{S}_{\#}\nu$, we see that we get an equality. Therefore for $\pi$-almost every $(x,y)$ we have
        \begin{equation*}
        \mathrm{dim}(V(\mathcal{S}(x))_1)=\mathrm{dim}(V(\mathcal{S}(y))_1)\text{ and }\mathrm{dim}(V(\mathcal{S}(y))_2)=\mathrm{dim}( V(\mathcal{S}(x))_2).
    \end{equation*}
    We infer that
      \begin{equation}\label{eqn:subspaces}
      V(\mathcal{S}(x))_1= V(\mathcal{S}(y))_1\text{ and }V(\mathcal{S}(y))_2= V(\mathcal{S}(x))_2\text{ for }\pi\text{-almost every }(x,y).
    \end{equation}
Thus, for $\sigma$-almost every $(x,y,z)$, $Du$ is isometric on $\{x,y,z\}$. Indeed, $z\in\mathcal{S}(x)\cap\mathcal{S}(y)$ by Lemma \ref{lem:three} and uniqueness of the leaves at $x,y$. By Lemma \ref{lem:isoderma} and (\ref{eqn:subspaces}), we get
    \begin{align*}
       & Du(x)-Du(y)=Du(x)-Du(z)-(Du(y)-Du(z))=\\
       &(P_{V(\mathcal{S}(x))_1}-P_{V(\mathcal{S}(x))_2})(x-z)-(P_{V(\mathcal{S}(y))_1}-P_{V(\mathcal{S}(y))_2})(y-z)=\\
       &(P_{V(\mathcal{S}(x))_1}-P_{V(\mathcal{S}(x))_2})(x-y).
    \end{align*}
    Since the leaves that contain $x,y$ are determined uniquely, almost everywhere, we see that $(x,y,z)$ belong to a common leaf of $Du$.
\end{proof}

\begin{lemma}\label{lem:v}
    The map
    \begin{equation*}
        x\mapsto \mathrm{dim}V(\mathcal{S}(x))_i\text{ for }i=1,2,
    \end{equation*}
    is Borel measurable.
\end{lemma}
\begin{proof}
    Note that the leaf $\mathcal{S}(x)$ containing $x\notin N(Du)$ is unique, and that $\mathcal{S}(x)=\{x\}$ on $N(Du)$. The dimension of $V(\mathcal{S}(x))_1$ is at least $k\in\mathbb{N}$ if and only if there exist $x_1,\dotsc,x_k\in\mathcal{S}(x)$ such that
    \begin{equation*}
        (x_j-x)+(Du(x_j)-Du(x))\text{ for }j=1,\dotsc,k\text{ are linearly independent.}
    \end{equation*}
This is equivalent to
  \begin{equation*}
       \mathrm{det}\big( \langle x_s-x+Du(x_s)-Du(x),x_t-x+Du(x_t)-Du(x)\rangle \big)_{s,t=1}^k\neq 0.
    \end{equation*}
    The set $D$ of $(k+1)$-tuples $(x,x_1,\dotsc,x_k)\in (\mathbb{R}^n)^{k+1}$ for which the above determinant is non-zero is an open set, as $Du$ is continuous. Note that the function
    \begin{equation*}
        x\mapsto \mathcal{S}(x)^k\in CC\big((\mathbb{R}^n)^k\big)
    \end{equation*}
    is Borel measurable, as the map $\mathbb{R}^n\ni x\mapsto \mathcal{S}(x)\in CC(\mathbb{R}^n)$ is Borel.
    Therefore, also set
    \begin{equation*}
        \{x\in\mathbb{R}^n\mid (\{x\}\times\mathcal{S}(x)^k)\cap D\neq\emptyset\}
    \end{equation*}
    is Borel, as measurability in the Wijsman topology is equivalent to measurability as a multifunction. The latter set is equal to the set of points $x$ for which $\mathrm{dim}V(\mathcal{S}(x))_1\geq k$.
    The measurability of $x\mapsto \mathrm{dim}V(\mathcal{S}(x))_2$ we prove analogously; it suffices to notice that the $\mathrm{dim}V(\mathcal{S}(x))_2\geq k$ for some $k\in\mathbb{N}$ if and only if
     there exist $x_1,\dotsc,x_k\in\mathcal{S}(x)$ such that
    \begin{equation*}
        (x_j-x)-(Du(x_j)-Du(x))\text{ for }j=1,\dotsc,k\text{ are linearly independent.}
    \end{equation*}
\end{proof}

\begin{theorem}\label{thm:parti}
   Let $\mu,\nu\in\mathcal{P}_2(\mathbb{R}^n)$ be absolutely continuous with respect to the Lebesgue measure, and let $u\in\mathcal{C}^{1,1}(\mathbb{R}^n)$ be an optimiser of $\mathcal{I}(\mu,\nu)$ and $\sigma\in\Sigma(\mu,\nu)$ be an optimiser of $\mathcal{J}(\mu,\nu)$. Suppose that there exists an optimal plan $\sigma_0\in\Sigma(\mu,\nu)$ for $\mathcal{J}(\mu,\nu)$ with absolutely continuous third marginal. Let $\theta\in\mathcal{P}(CC(\mathbb{R}^n))$, $\mu_{\mathcal{S}},\nu_{\mathcal{S}}\in\mathcal{P}(\mathbb{R}^n)$ defined for $\theta$-almost every leaf $\mathcal{S}$ of $Du$, be as in Theorem \ref{thm:partim}.
Then,
\begin{equation}\label{eqn:bary}
  \mu_{\mathcal{S}},\nu_{\mathcal{S}}\in\mathcal{P}_2(\mathbb{R}^n),  \int_{\mathbb{R}^n}x\,d\mu_{\mathcal{S}}(x)=\int_{\mathbb{R}^n}x\,d\nu_{\mathcal{S}}(x), \text{ for }\theta\text{-almost every }\mathcal{S},
\end{equation}
and
\begin{equation}\label{eqn:sigmapierw}
    \sigma_{\mathcal{S}}\in\Sigma(\mu_{\mathcal{S}},\nu_{\mathcal{S}})\text{ for }\theta\text{-almost every }\mathcal{S}.
\end{equation}
Moreover $u$ is optimal for  $\mathcal{I}(\mu_{\mathcal{S}},\nu_{\mathcal{S}})$ and
 $\sigma_{\mathcal{S}}$ is optimal for $\mathcal{J}(\mu_{\mathcal{S}},\nu_{\mathcal{S}})$.
\end{theorem}
\begin{proof}
By Theorem \ref{thm:partim} applied to $\sigma_0$, we have $\mathcal{S}_{\#}\mu=\mathcal{S}_{\#}\nu$. By Corollary \ref{col:oneleaf}, $x,y,z$ belong to a common leaf of $Du$ for $\sigma$-almost every $(x,y,z)$. Since $\mu,\nu\ll\lambda$, the leaves at $x$ and $y$ are unique there. Hence $\mathcal{S}(x)=\mathcal{S}(y)$ and $z\in\mathcal{S}(x)$ for $\sigma$-almost every $(x,y,z)$.
The integrability of second moments of $\mu_{\mathcal{S}},\nu_{\mathcal{S}}$ follows directly by the assumption that $\mu,\nu\in\mathcal{P}_2(\mathbb{R}^n)$ and Theorem \ref{thm:partim}. Finiteness of the cost gives a finite second moment for $\mathrm{P}_3\sigma$, hence also for $\mathrm{P}_3\sigma_{\mathcal S}$ for $\theta$-almost every $\mathcal S$.

Now, thanks to \ref{i:pro} of Theorem \ref{thm:disintegration}, for any Borel set $E\subset CC(\mathbb{R}^n)$ and any continuous bounded $f\colon\mathbb{R}^n\to\mathbb{R}$ we have
\begin{equation*}
\int_{\mathcal{S}^{-1}(E)\times(\mathbb{R}^n)^2} f(x)\,d\sigma(x,y,z)
=\int_E\int_{(\mathbb{R}^n)^3} f(x)\,d\sigma_{\mathcal{S}}(x,y,z)\,d\theta(\mathcal{S}).
\end{equation*}
Since the first marginal of $\sigma$ is $\mu$, it follows that
\begin{equation*}
    \int_{E}  \int_{(\mathbb{R}^n)^3} f(x)d\sigma_{\mathcal{S}}(x,y,z)\, d\theta(\mathcal{S})=\int_{E}  \int_{\mathbb{R}^n} fd\mu_{\mathcal{S}}\, d\theta(\mathcal{S}).
\end{equation*}
Since these equalities are valid for any Borel set $E\subset CC(\mathbb{R}^n)$, we infer that for any continuous and bounded $f$
\begin{equation*}
     \int_{(\mathbb{R}^n)^3} f(x)d\sigma_{\mathcal{S}}(x,y,z)= \int_{\mathbb{R}^n} fd\mu_{\mathcal{S}}\text{ for }\theta\text{-almost every }\mathcal{S}.
\end{equation*}
Thanks to separability of the space of continuous and compactly supported functions on $\mathbb{R}^n$ we infer that
\begin{equation*}
    \mathrm{P}_1\sigma_{\mathcal{S}}=\mu_{\mathcal{S}}\text{ for }\theta\text{-almost every }\mathcal{S}.
\end{equation*}
Using $\mathcal{S}(x)=\mathcal{S}(y)$, we get
\begin{equation*}
    \mathbf{1}_{\mathcal{S}^{-1}(E)}(x)=\mathbf{1}_{\mathcal{S}^{-1}(E)}(y)\text{ for }\sigma\text{-almost every }(x,y,z)\in (\mathbb{R}^n)^3.
\end{equation*}
The same argument with $f$ depending on the second coordinate gives
\begin{equation*}
    \mathrm{P}_2\sigma_{\mathcal{S}}=\nu_{\mathcal{S}}\text{ for }\theta\text{-almost every }\mathcal{S}.
\end{equation*}
Testing (\ref{eqn:mart}) with $g\mathbf 1_{\mathcal S^{-1}(E)}$ and $h\mathbf 1_{\mathcal S^{-1}(E)}$, respectively, and disintegrating gives
\begin{align*}
 \int_E\!\int_{(\mathbb R^n)^3}(z-x)g(x)\,d\sigma_{\mathcal S}(x,y,z)\,d\theta(\mathcal S)&=0,\\
 \int_E\!\int_{(\mathbb R^n)^3}(z-y)h(y)\,d\sigma_{\mathcal S}(x,y,z)\,d\theta(\mathcal S)&=0.
\end{align*}
Since $E$ is arbitrary, both inner integrals vanish for $\theta$-almost every $\mathcal S$. Taking $g,h$ first to be indicators from a countable generating algebra and then applying the monotone class theorem gives (\ref{eqn:mart}) for all bounded Borel $g,h$ on a common set of full $\theta$-measure. Taking $g=h=1$ gives the common barycentre. This proves (\ref{eqn:bary}) and (\ref{eqn:sigmapierw}).

We know that
\begin{equation*}
    \int_{\mathbb{R}^n}u\, d(\nu-\mu)=\int_{(\mathbb{R}^n)^3} c(x,y,z)\, d\sigma(x,y,z)
\end{equation*}
and the disintegration implies that
\begin{equation*}
 \int_{CC(\mathbb R^n)}\!\int_{\mathbb R^n} u\,d(\nu_{\mathcal S}-\mu_{\mathcal S})\,d\theta(\mathcal S)
 =\int_{CC(\mathbb R^n)}\!\int_{(\mathbb R^n)^3} c\,d\sigma_{\mathcal S}\,d\theta(\mathcal S).
\end{equation*}
Since $u\in\mathcal{C}^{1,1}(\mathbb{R}^n)$ has $1$-Lipschitz derivative, and $\sigma_{\mathcal{S}}\in\Sigma(\mu_{\mathcal{S}},\nu_{\mathcal{S}})$, we infer, thanks to weak duality, that  $u$ is optimal for  $\mathcal{I}(\mu_{\mathcal{S}},\nu_{\mathcal{S}})$ and
 $\sigma_{\mathcal{S}}$ is optimal for $\mathcal{J}(\mu_{\mathcal{S}},\nu_{\mathcal{S}})$.
\end{proof}

\section{Convex-concave order, bimartingale couplings and case of isometric derivative}\label{s:bimart}

We shall consider the case when for optimal $u$ its derivative $Du$ is an isometry. For any absolutely continuous measures $\mu,\nu\in\mathcal{P}_2(\mathbb{R}^n)$ with common barycentre for which some optimal plan has absolutely continuous third marginal, the problem $\mathcal{J}(\mu,\nu)$ can be reduced to a collection of problems with optimal $Du$ being an isometry, as follows by the results of Section \ref{s:disintegrate}, see Theorem \ref{thm:parti}.

\subsection{Convex-concave order}\label{s:convexconcave}

In order to analyse the case when the derivative of an optimal $u$ is isometric, we shall introduce a notion of convex-concave order for pairs of probability measures.

We recall that two subspaces $V_1,V_2\subset\mathbb{R}^n$ are said to be mutually complementing if $V_1\oplus V_2=\mathbb{R}^n$.

\begin{definition}\label{def:cc}
    Let $V_1, V_2$ be two mutually complementing orthogonal subspaces. We shall say that $f\colon\mathbb{R}^n\to\mathbb{R}$ is \emph{convex-concave} with respect to $V_1,V_2$ whenever for each $v_1\in V_1$ the function
\begin{equation*}
 V_2  \ni v_2\mapsto f(v_1+v_2)\in\mathbb{R}
\end{equation*}
 is concave and for any $v_2\in V_2$ the function
 \begin{equation*}
 V_1 \ni  v_1\mapsto f(v_1+v_2)\in\mathbb{R}
\end{equation*}
is convex.  We shall also say that $f$ has at most \emph{quadratic growth} if there is $C\geq 0$ such that for all $v\in\mathbb{R}^n$
\begin{equation*}
    \abs{f(v)}\leq C(1+\norm{v}^2).
\end{equation*}
\end{definition}

Below, for a matrix $A\in\mathbb{R}^{n\times n }$ we write $A\geq 0$ if $A$ is positive semi-definite.
\begin{lemma}\label{lem:cc}
   Let $V_1, V_2$ be two mutually complementing orthogonal subspaces. Suppose that $f\in\mathcal{C}^{1,1}(\mathbb{R}^n)$. Then $f$ is convex-concave with respect to $V_1,V_2$ if and only if
   \begin{equation*}
     0\leq  P_{V_1} D^2fP_{V_1}\text{ and } 0\geq  P_{V_2} D^2fP_{V_2}\text{ almost everywhere}
   \end{equation*}
   where $P_{V_i}$ are orthogonal projections onto $V_i$ for $i=1,2.$
\end{lemma}
\begin{proof}
    This follows directly by Definition \ref{def:cc} and the characterisation of convexity by means of positive semi-definiteness of second derivative almost everywhere.
\end{proof}

\begin{definition}\label{def:convex}
    Let $\mu,\nu\in\mathcal{P}_2(\mathbb{R}^n)$ be two Borel probability measures with finite second moments.
We shall say that $\mu,\nu$ are in \emph{convex-concave order} and write $\mu\prec_{c-c}\nu$ if there exist mutually complementing orthogonal subspaces $V_1,V_2$ such that for any $f\colon\mathbb{R}^n\to\mathbb{R}$ that is convex-concave with respect to $V_1,V_2$  and of at most quadratic growth
\begin{equation*}
    \int_{\mathbb{R}^n}f\, d\mu\leq\int_{\mathbb{R}^n}f\, d\nu.
\end{equation*}
We also say that $\mu,\nu$ are in \emph{convex-concave order} with respect to $V_1,V_2$, when the above holds for specific orthogonal, mutually complementing linear subspaces $V_1,V_2$.
\end{definition}

Note that the convex-concave order with respect to $V_1,V_2$ is transitive. However, this is no longer true if we allow to vary the subspaces.

\begin{remark}
Any $f\in\mathcal{C}^{1,1}(\mathbb{R}^n)$ has at most quadratic growth and therefore the integrals in Definition \ref{def:convex} are well-defined for any convex-concave $f\in\mathcal{C}^{1,1}(\mathbb{R}^n)$.
\end{remark}

We cite the following classical result after \cite[Lemma 3.4]{Bolbotowski2024}.
\begin{lemma}\label{lem:afterb}
 Let $v\colon\mathbb{R}^n\to\mathbb{R}$ be continuous. Then the following conditions are equivalent:
 \begin{enumerate}[label=(\roman*)]
     \item $v\in\mathcal{C}^{1,1}(\mathbb{R}^n)$ and $v$ has $1$-Lipschitz derivative,
     \item $v\in \mathcal{C}^{1,1}(\mathbb{R}^n)$ and $-\mathrm{Id}\leq D^2v\leq\mathrm{Id}$ almost everywhere,
     \item\label{i:plusminus} the functions $\frac12\norm{\cdot}^2\pm v$ are convex.
 \end{enumerate}
\end{lemma}

\begin{remark}
    Let us note that the condition \ref{i:plusminus} can be conveniently rephrased as a requirement that for all $t\in [0,1]$, $x,y\in\mathbb{R}^n$
    \begin{equation*}
        \Big\lvert v\big(tx+(1-t)y\big)-\big(tv(x)+(1-t)v(y)\big)\Big\rvert\leq \frac12 t(1-t)\norm{x-y}^2.
    \end{equation*}
    This follows directly from the definition of convexity.
 Another equivalent condition is that for each $\mu,\nu\in\mathcal{P}_2(\mathbb{R}^n)$  with common barycentre $b\in\mathbb{R}^n$
 \begin{equation*}
     \Big\lvert \int_{\mathbb{R}^n}v\, d(\nu-\mu)\Big\rvert\leq \frac12\int_{\mathbb{R}^n} \norm{\cdot-b}^2\, d(\mu+\nu).
 \end{equation*}
 The proof follows by the Jensen inequality.
\end{remark}

\begin{proposition}\label{pro:cc}
      Let $\mu,\nu\in\mathcal{P}_2(\mathbb{R}^n)$ be two Borel probability measures with common barycentre. Then the following conditions are equivalent:
      \begin{enumerate}[label=(\roman*)]
      \item\label{i:cc}$\mu\prec_{c-c}\nu$ with respect to some mutually complementing, orthogonal subspaces $V_1,V_2$,
          \item\label{i:cclip}there exist some mutually complementing, orthogonal subspaces $V_1,V_2$ such that for any $f\in\mathcal{C}^{1,1}(\mathbb{R}^n)$ convex-concave with respect to $V_1,V_2$
          \begin{equation*}
              \int_{\mathbb{R}^n}f\, d\mu\leq\int_{\mathbb{R}^n}f\, d\nu,
          \end{equation*}
                \item\label{i:isom} the supremum in the definition of  $\mathcal{I}(\mu,\nu)$ is attained by  $u\in\mathcal{C}^{1,1}(\mathbb{R}^n)$ whose derivative is an isometry.
          \item\label{i:projec} there exist  mutually complementing, orthogonal subspaces $V_1,V_2$ such that the supremum in $\mathcal{I}(\mu,\nu)$
        is attained by the function  $u_{V_1,V_2}\colon\mathbb{R}^n\to\mathbb{R}$ given by the  formula
          \begin{equation*}
              u_{V_1,V_2}(v)=\frac12\norm{P_{V_1}v}^2-\frac12\norm{P_{V_2}v}^2\text{ for all }v\in\mathbb{R}^n.
          \end{equation*}
          Here $P_{V_i}$ denotes the orthogonal projection onto $V_i$ for $i=1,2$.

      \end{enumerate}
\end{proposition}
\begin{proof}
    Suppose that there exist orthogonal, mutually complementing subspaces $V_1,V_2$ for which $u_{V_1,V_2}$ attains the supremum in $\mathcal{I}(\mu,\nu)$, i.e., that \ref{i:projec} holds true.

Let $f\in \mathcal{C}^{1,1}(\mathbb{R}^n)$ have $1$-Lipschitz derivative and let it be convex-concave with respect to $V_1,V_2$.

Then by Lemma \ref{lem:cc} and Lemma \ref{lem:afterb} almost everywhere we have
    \begin{equation}\label{eqn:bounds}
      0\leq  P_{V_1}D^2fP_{V_1} \leq P_{V_1}\text{ and }      -P_{V_2}\leq  P_{V_2}D^2fP_{V_2} \leq 0.
    \end{equation}
    Let $A=P_{V_1}-P_{V_2}$. Calculation yields that
    \begin{equation*}
        AD^2f+D^2fA=2P_{V_1}D^2fP_{V_1}-2P_{V_2}D^2fP_{V_2}\geq 0.
    \end{equation*}
    Then for $t>0$
    \begin{equation*}
        (A-tD^2f)^2=\mathrm{Id}+t^2(D^2f)^2-t(AD^2f+D^2fA)\leq \mathrm{Id}+t^2(D^2f)^2
    \end{equation*}
    Since $A-tD^2f$ is symmetric, it follows that its operator norm is bounded
    \begin{equation}\label{eqn:oneone}
        \norm{A-tD^2f}\leq \sqrt{1+t^2\norm{D^2f}^2}\leq \sqrt{1+t^2}.
    \end{equation}
    For $t>0$, consider the function $u_{tf}\colon\mathbb{R}^n\to\mathbb{R}$ defined by the formula
\begin{equation*}
    u_{tf}(v)=\frac1{\sqrt{1+t^2}}\bigg(\frac12\norm{P_{V_1}v}^2-\frac12\norm{P_{V_2}v}^2-tf(v)\bigg)\text{ for }v\in\mathbb{R}^n.
\end{equation*}
Then almost everywhere
\begin{equation*}
    D^2u_{tf}=\frac{P_{V_1}-P_{V_2}-tD^2f}{\sqrt{1+t^2}}.
\end{equation*}
Thus, (\ref{eqn:oneone}) shows that $u_{tf}$ has $1$-Lipschitz derivative.

This shows that
\begin{equation*}
    \int_{\mathbb{R}^n}u_{tf}\, d(\nu-\mu)\leq \int_{\mathbb{R}^n}u_{V_1,V_2}\, d(\nu-\mu),
\end{equation*}
and consequently, after rearrangement and taking limit $t\to 0$, we get
\begin{equation*}
    \int_{\mathbb{R}^n}f\, d(\nu-\mu)\geq 0.
\end{equation*}
By homogeneity we see that the same inequality is valid for any $f\in\mathcal{C}^{1,1}(\mathbb{R}^n)$. We conclude that \ref{i:cclip} holds.

Conversely, \ref{i:cclip}. Consider for any $u\in\mathcal{C}^{1,1}(\mathbb{R}^n)$ with $1$-Lipschitz derivative a function $f\colon\mathbb{R}^n\to\mathbb{R}$ given by the formula
\begin{equation*}
    f(v)=\frac12\norm{P_{V_1}v}^2-\frac12\norm{P_{V_2}v}^2-u(v)\text{ for }v\in\mathbb{R}^n.
\end{equation*}
Then almost everywhere
\begin{equation*}
    D^2f=P_{V_1}-P_{V_2}-D^2u,
\end{equation*}
which is positive semi-definite on $V_1$ and negative semi-definite on $V_2$. By Lemma \ref{lem:cc}, $f$ is convex-concave with respect to $V_1,V_2$. Therefore
\begin{equation*}
    \int_{\mathbb{R}^n}f\, d(\nu-\mu)\geq 0,
\end{equation*}
which implies that indeed $u_{V_1,V_2}$ attains the supremum in the definition of $\mathcal{I}(\mu,\nu)$. This is to say, \ref{i:cclip} and \ref{i:projec} are equivalent.

Now, Lemma \ref{lem:isoderma} shows that \ref{i:isom} implies \ref{i:projec}. Indeed, since $\mu,\nu$ have common mass and common barycentre, by Lemma \ref{lem:isoderma}, if \ref{i:isom} holds true, then there exist mutually complementing, orthogonal subspaces $V_1,V_2$ such that $u_{V_1,V_2}$ attains the supremum in (\ref{eqn:i}). The converse implication is trivial.
The fact that \ref{i:projec} implies \ref{i:cc} is a consequence of the existence of bimartingale couplings, see Theorem \ref{thm:bimart} and Corollary \ref{col:twozet}.
\end{proof}

\subsection{Bimartingale couplings}\label{s:couplings}

\begin{definition}\label{def:bimart}
     Let $\mu,\nu\in\mathcal{P}_1(\mathbb{R}^n)$ be two Borel probability measures with finite moments. Let $V_1,V_2$ be two orthogonal, mutually complementing subspaces. We say that a Borel probability measure $\pi\in\mathcal{P}(\mathbb{R}^n\times\mathbb{R}^n)$ is a \emph{bimartingale coupling} between $\mu$ and $\nu$ with respect to $V_1,V_2$ if and only if its respective marginals are equal to $\mu,\nu$, i.e., $\pi\in\Gamma(\mu,\nu)$, and for any bounded Borel functions $g,h\colon \mathbb{R}^n\to\mathbb{R}$
      \begin{equation}\label{eqn:bary1}
          \int_{\mathbb{R}^n\times\mathbb{R}^n}P_{V_1}(y)g(x)\, d\pi(x,y)=     \int_{\mathbb{R}^n\times\mathbb{R}^n}P_{V_1}(x)g(x)\, d\pi(x,y)
      \end{equation}
      and
            \begin{equation}\label{eqn:bary2}
          \int_{\mathbb{R}^n\times\mathbb{R}^n}P_{V_2}(x)h(y)\, d\pi(x,y)=     \int_{\mathbb{R}^n\times\mathbb{R}^n}P_{V_2}(y)h(y)\, d\pi(x,y).
      \end{equation}
      Here $P_{V_i}$ is the orthogonal projection onto $V_i$ for $i=1,2$.
      By $\Gamma_{bm}(\mu,\nu,V_1,V_2)$ we denote the set of bimartingale couplings between $\mu$ and $\nu$ with respect to $V_1,V_2$.
      We say that $\pi$ is a bimartingale coupling between $\mu$ and $\nu$ if there exist orthogonal, mutually complementing subspaces $V_1,V_2$ such that $\pi\in\Gamma_{bm}(\mu,\nu,V_1,V_2)$. By $\Gamma_{bm}(\mu,\nu)$ we denote the set of bimartingale couplings between $\mu$ and $\nu$.
\end{definition}

\begin{remark}
    If $V_1=\mathbb{R}^n$ and $V_2=\{0\}$ then any bimartingale coupling  between $\mu$ and $\nu$ with respect to $V_1,V_2$ is a martingale coupling. Indeed, the conditions of Definition \ref{def:bimart} for a coupling $\pi$ between $\mu$ and $\nu$ are that for any bounded Borel function  $g\colon\mathbb{R}^n\to\mathbb{R}$
     \begin{equation*}
          \int_{\mathbb{R}^n\times\mathbb{R}^n}yg(x)\, d\pi(x,y)=     \int_{\mathbb{R}^n\times\mathbb{R}^n}xg(x)\, d\pi(x,y),
      \end{equation*}
      which is precisely the martingale coupling condition for $\pi$.
\end{remark}

\begin{proposition}\label{pro:swaps}
If $\pi\in\Gamma_{bm}(\mu,\nu,V_1,V_2)$ and if $S\colon V_2\times V_2\to\ V_2\times V_2$ swaps the variables
\begin{equation*}
    S(w_1,w_2)=(w_2,w_1)\text{ for }w_1,w_2\in V_2,
\end{equation*}
then the pushforwards
\begin{equation*}
    (P_{V_1},P_{V_1})_{\#}\pi,(S(P_{V_2},P_{V_2}))_{\#}\pi
\end{equation*}
are martingale couplings in $\Gamma((P_{V_1})_{\#}\mu,(P_{V_1})_{\#}\nu)$ and in $\Gamma((P_{V_2})_{\#}\nu,(P_{V_2})_{\#}\mu)$ respectively.
\end{proposition}
\begin{proof}
    The condition that
   \begin{equation}\label{eqn:warun}
       (P_{V_1},P_{V_1})_{\#}\pi\text{ and } S(P_{V_2},P_{V_2})_{\#}\pi\text{ are martingale transports }
   \end{equation}
   is equivalent to assuming that for all bounded Borel $g\colon V_1\to\mathbb{R},h\colon V_2\to\mathbb{R}$ we have
    \begin{equation*}
        \int_{\mathbb{R}^n\times \mathbb{R}^n}P_{V_1}(y)g(P_{V_1}x)\, d\pi(x,y)=\int_{\mathbb{R}^n\times\mathbb{R}^n}P_{V_1}(x)g(P_{V_1}x) \, d\pi(x,y)
    \end{equation*}
    and
        \begin{equation*}
        \int_{\mathbb{R}^n\times \mathbb{R}^n}P_{V_2}(x)h(P_{V_2}y)\, d\pi(x,y)=\int_{\mathbb{R}^n\times\mathbb{R}^n}P_{V_2}(y)\,h(P_{V_2}y) d\pi(x,y).
    \end{equation*}
    Thus, the condition (\ref{eqn:warun})  follows immediately from (\ref{eqn:bary1}) and from (\ref{eqn:bary2}).
\end{proof}

 \begin{example}
    It is not true that the condition of Proposition \ref{pro:swaps} that the pushforwards
\begin{equation*}
    (P_{V_1},P_{V_1})_{\#}\pi\in\Gamma((P_{V_1})_{\#}\mu,(P_{V_1})_{\#}\nu)\text{ and }(S(P_{V_2},P_{V_2}))_{\#}\pi\in\Gamma((P_{V_2})_{\#}\nu,(P_{V_2})_{\#}\mu),
\end{equation*}
are martingale couplings is equivalent to $\pi\in\Gamma_{bm}(\mu,\nu,V_1,V_2)$. As an example, let $V_1$ be spanned by the first standard vector and $V_2$ be spanned by the second standard vector and consider the measure $\pi\in\mathcal{P}(\mathbb{R}^2\times\mathbb{R}^2)$ given by the formula
\begin{equation*}
    \pi= \frac12\big(\delta_{(0,1,-1,1)}+\delta_{(0,-1,1,-1)}\big).
\end{equation*}
Then
\begin{equation*}
    \mathrm{P}_1\pi=\mu=\frac12\big(\delta_{(0,1)}+\delta_{(0,-1)}\big)\text{ and }\mathrm{P}_2\pi=\nu=\frac12\big(\delta_{(-1,1)}+\delta_{(1,-1)}\big).
\end{equation*}
Moreover
\begin{equation*}
    (P_{V_1})_{\#}\mu=\delta_0,(P_{V_1})_{\#}\nu=\frac12(\delta_{-1}+\delta_1),  (P_{V_2})_{\#}\mu=\frac12(\delta_{-1}+\delta_1)=(P_{V_2})_{\#}\nu
\end{equation*}
and
\begin{equation*}
      (P_{V_1},P_{V_1})_{\#}\pi=\frac12\big(\delta_{(0,-1)}+\delta_{(0,1)}\big)\text{ and }   S(P_{V_2},P_{V_2})_{\#}\pi=\frac12\big(\delta_{(1,1)}+\delta_{(-1,-1)}\big).
\end{equation*}
Thus indeed the pushforwards are martingale couplings between their respective marginals. However, $\pi$ is not a bimartingale coupling. Indeed, for $g\colon\mathbb{R}^2\to\mathbb{R}$ given by $g=\mathbf{1}_{(0,1)}$ we have
\begin{equation*}
    \int_{\mathbb{R}^2\times\mathbb{R}^2}P_{V_1}(y)g(x)\, d\pi(x,y)=-\frac12 \neq 0= \int_{\mathbb{R}^2\times\mathbb{R}^2}P_{V_1}(x)g(x)\, d\pi(x,y),
\end{equation*}
so that (\ref{eqn:bary1}) is not satisfied.
\end{example}

\subsection{Relation of bimartingale couplings to convex-concave order -- the Strassen theorem for convex-concave order}\label{s:relation}

\begin{remark}
Lemma \ref{lem:three}, Corollary \ref{col:three} and \cite[Theorem 1.1, (ii)]{Bolbotowski2024}  show that a symmetric isometry $T\colon\mathbb{R}^n\to\mathbb{R}^n$ is the derivative of an optimal $u$ in the problem  (\ref{eqn:i}) if and only if conditions (\ref{eqn:yp}) and (\ref{eqn:xp}) are satisfied for $\sigma$-almost every $(x,y,z)\in(\mathbb{R}^n)^3$, where $\sigma\in\Sigma(\mu,\nu)$ is optimal for $\mathcal{J}(\mu,\nu)$, where $u(x)=\frac12\langle Tx,x\rangle$ for $x\in\mathbb{R}^n$.
In \cite[Example 4.2]{Bolbotowski2024}, a similar observation has been made in the case of Gaussian measures.
\end{remark}

The following theorem is the analogue of Strassen's theorem for the convex-concave order.
The two martingale constraints correspond to the two opposite curvatures of the dual potential. In the convex directions the martingale constraint is imposed from the first marginal to the second, whereas in the concave directions the roles of the marginals are reversed.

\begin{theorem}\label{thm:bimart}
      Let $\mu,\nu\in\mathcal{P}_2(\mathbb{R}^n)$ be two Borel probability measures with finite second moments and the same barycentre. Then $\mu\prec_{c-c}\nu$ if and only if there exists a bimartingale coupling $\pi$ between $\mu$ and $\nu$.

      More precisely, if $V_1,V_2$ are orthogonal, mutually complementing subspaces, then there exists $\pi\in\Gamma_{bm}(\mu,\nu, V_1,V_2)$ if and only if $\mu\prec_{c-c}\nu$ with respect to $V_1,V_2$.
      If $\sigma\in\Sigma(\mu,\nu)$ is optimal for $\mathcal{J}(\mu,\nu)$ and  $\mu\prec_{c-c}\nu$ with respect to $V_1,V_2$, then $\mathrm{P}_{12}\sigma\in\Gamma_{bm}(\mu,\nu,V_1,V_2)$ and $\sigma=R_{\#}\mathrm{P}_{12}\sigma$, where
      \begin{equation*}
          R(x,y)=(x,y,P_{V_2}x+P_{V_1}y)\text{ for }x,y\in\mathbb{R}^n.
      \end{equation*}
\end{theorem}
\begin{proof}
    Suppose first that there exists a bimartingale coupling $\pi\in\Gamma_{bm}(\mu,\nu)$ with respect to $V_1,V_2$. Let $f\colon\mathbb{R}^n\to\mathbb{R}$ be a convex-concave function with respect to $V_1,V_2$ with at most quadratic growth. Then by convexity of $f$ with respect to $V_1$, by the Jensen inequality and by the barycentre condition (\ref{eqn:bary1})  we get
    \begin{equation*}
\int_{\mathbb{R}^n}f\, d\mu=\int_{\mathbb{R}^n\times\mathbb{R}^n} f(x)\, d\pi(x,y)\leq\int_{\mathbb{R}^n\times\mathbb{R}^n} f(P_{V_1}(y)+P_{V_2}(x))\, d\pi(x,y) .
    \end{equation*}
    Now,  by concavity of $f$ with respect to $V_2$, the Jensen inequality and the barycentre condition (\ref{eqn:bary2}) we get that
        \begin{equation*}
\int_{\mathbb{R}^n}f\, d\nu=\int_{\mathbb{R}^n\times\mathbb{R}^n} f(y)\, d\pi(x,y)\geq\int_{\mathbb{R}^n\times\mathbb{R}^n} f(P_{V_1}(y)+P_{V_2}(x))\, d\pi(x,y) .
    \end{equation*}
    Thus
    \begin{equation*}
        \int_{\mathbb{R}^n}f\, d\mu\leq\int_{\mathbb{R}^n}f\, d\nu.
    \end{equation*}
    We see that indeed $\mu\prec_{c-c}\nu$.

Let us now suppose that $\mu\prec_{c-c}\nu$. Proposition \ref{pro:cc}, \ref{i:cc} to \ref{i:projec}, shows that for some mutually complementing, orthogonal subspaces $V_1,V_2$ the function $u_{V_1,V_2}$ attains the supremum  in the definition of $\mathcal{I}(\mu,\nu)$.

Now, by \cite[Theorem 1.1]{Bolbotowski2024}, choose an optimal plan $\sigma\in\Sigma(\mu,\nu)$ for $\mathcal{J}(\mu,\nu)$.

Moreover by \cite[Theorem 1.1, (ii)]{Bolbotowski2024}, for $\sigma$-almost every $(x,y,z)\in(\mathbb{R}^n)^3$ we have
 \begin{equation*}
 \medmuskip=1mu \thickmuskip=2mu \big(u_{V_1,V_2}(y)+Du_{V_1,V_2}(y)(z-y)\big)-\big(u_{V_1,V_2}(x)+Du_{V_1,V_2}(x)(z-x)\big)=\frac12\Big(\norm{y-z}^2+\norm{x-z}^2\Big).
\end{equation*}
Corollary \ref{col:three}, with $z_0=0$, shows that
\begin{equation}\label{eqn:zetka}
    z=P_{V_2}x+P_{V_1}y\text{ for }\sigma\text{-almost every }(x,y,z)\in(\mathbb{R}^n)^3.
\end{equation}
This shows that $\sigma$ is determined by its projection onto the first two coordinates: $\sigma=R_{\#}\mathrm{P}_{12}\sigma$.
The condition that $\mathrm{P}_{13}\sigma, \mathrm{P}_{23}\sigma$ are martingale transports is equivalent to that for any bounded, Borel functions $g,h\colon\mathbb{R}^n\to\mathbb{R}$ we have
      \begin{equation}\label{eqn:bary1three}
          \int_{(\mathbb{R}^n)^3}(z-x)g(x)\, d\sigma(x,y,z)= 0
      \end{equation}
      and
            \begin{equation}\label{eqn:bary2three}
          \int_{(\mathbb{R}^n)^3}(z-y)h(y)\, d\sigma(x,y,z)=0.
      \end{equation}
    By (\ref{eqn:zetka}) we see that
    \begin{equation*}
        z-x=P_{V_1}(y-x), z-y=P_{V_2}(x-y)\text{ for }\sigma\text{-almost every }(x,y,z)\in(\mathbb{R}^n)^3.
    \end{equation*}
    Let $\pi=\mathrm{P}_{12}\sigma\in\Gamma(\mu,\nu)$. We see that (\ref{eqn:bary1three}) and (\ref{eqn:bary2three}) are equivalent respectively to
       \begin{equation*}
          \int_{\mathbb{R}^n\times\mathbb{R}^n}P_{V_1}(y-x)g(x)\, d\pi(x,y)= 0
      \end{equation*}
      and to
     \begin{equation*}
          \int_{\mathbb{R}^n\times\mathbb{R}^n}P_{V_2}(x-y)h(y)\, d\pi(x,y)=0.
      \end{equation*}
      This shows that $\pi$ is a bimartingale coupling between $\mu$ and $\nu$ with respect to $V_1,V_2$.
\end{proof}

\begin{remark}
    Suppose that $\mathcal{I}(\mu,\nu)$ is attained by a function with isometric derivative. The proof of Theorem \ref{thm:bimart} shows that to any $\sigma\in\mathcal{P}((\mathbb{R}^n)^3)$ that attains the infimum $\mathcal{J}(\mu,\nu)$ we may associate a bimartingale coupling $\pi\in\Gamma_{bm}(\mu,\nu)$. However, the converse is also true, which is an assertion of the following proposition.
\end{remark}

\begin{proposition}\label{pro:bimart}
  Let $\mu,\nu\in\mathcal{P}_2(\mathbb{R}^n)$ have common barycentre. Suppose that $\mathcal{I}(\mu,\nu)$ is attained by a function $u\in\mathcal{C}^{1,1}(\mathbb{R}^n)$ with isometric derivative, which, after modifying by an affine function can be taken to be of the form $u=u_{V_1,V_2}$. Suppose that $\pi\in\Gamma_{bm}(\mu,\nu,V_1,V_2)$.
Let $\sigma\in\mathcal{P}((\mathbb{R}^n)^3)$ be the pushforward of $\pi$ via the map
  \begin{equation*}
 R\colon   \mathbb{R}^n\times\mathbb{R}^n\to \mathbb{R}^n\times\mathbb{R}^n\times\mathbb{R}^n
  \end{equation*}
  defined by the formula
  \begin{equation*}
      R(x,y)= (x,y,P_{V_2}x+P_{V_1}y)\text{ for }x,y\in\mathbb{R}^n.
  \end{equation*}
  Then $\sigma\in\Sigma(\mu,\nu)$ is  optimal for $\mathcal{J}(\mu,\nu)$.
\end{proposition}
\begin{proof}
By \cite[Theorem 1.1, (ii)]{Bolbotowski2024} and Lemma \ref{lem:three} if $\sigma\in \Sigma(\mu,\nu)$, then it is an optimiser if and only if the conditions (\ref{eqn:yp}) and (\ref{eqn:xp})  hold for $\sigma$-almost every $(x,y,z)\in(\mathbb{R}^n)^3$.
By Corollary \ref{col:three} we see that these are equivalent to
\begin{equation*}
    z=P_{V_2}x+P_{V_1}y\text{ for }\sigma\text{-almost every }(x,y,z)\in(\mathbb{R}^n)^3.
\end{equation*}
Therefore it suffices to check that $\sigma\in\Sigma(\mu,\nu)$. This condition is equivalent to $\pi$ being a bimartingale coupling between $\mu$ and $\nu$ with respect to $V_1,V_2$, as observed in the proof of Theorem \ref{thm:bimart}.
\end{proof}

\begin{remark}\label{rem:decom}
  Suppose that $V_1,V_2$ are mutually orthogonal, complementing subspaces and that optimal $u$ for $\mathcal{I}(\mu,\nu)$ has isometric derivative $T=P_{V_1}-P_{V_2}$.  If $z=P_{V_2}x+P_{V_1}y$, then
    \begin{equation*}
        \norm{z-x}^2+\norm{z-y}^2=\norm{x-y}^2.
    \end{equation*}
    Therefore, in view of Theorem \ref{thm:bimart} and Proposition \ref{pro:bimart}, the variational problem
    \begin{equation*}
        \inf\Big\{\int_{(\mathbb{R}^n)^3} \frac12(\norm{z-x}^2+\norm{z-y}^2)\, d\sigma(x,y,z)\mid \sigma\in\Sigma(\mu,\nu)\Big\}
    \end{equation*}
is equivalent to
        \begin{equation*}
        \inf\Big\{\int_{(\mathbb{R}^n)^2} \frac12\norm{x-y}^2\, d\pi(x,y)\mid \pi\in\Gamma_{bm}(\mu,\nu,V_1,V_2)\Big\}.
    \end{equation*}
    The latter problem gives the same value for any $\pi\in\Gamma_{bm}(\mu,\nu,V_1,V_2)$.
    Specifically, for  $\pi\in\Gamma_{bm}(\mu,\nu,V_1,V_2)$
    \begin{align*}
    & \int_{(\mathbb{R}^n)^2} \frac12\norm{x-y}^2\, d\pi(x,y)=  \\
    &\int_{V_1\times V_1} \frac12\norm{x-y}^2\, d(P_{V_1},P_{V_1})_{\#}\pi(x,y)+\int_{V_2\times V_2} \frac12\norm{x-y}^2\, d(P_{V_2},P_{V_2})_{\#}\pi(x,y)=\\
    &\frac12\Big( \mathrm{var}(P_{V_1})_{\#}\nu-\mathrm{var}(P_{V_1})_{\#}\mu+\mathrm{var}(P_{V_2})_{\#}\mu-\mathrm{var}(P_{V_2})_{\#}\nu\Big).
    \end{align*}
    Here for a measure $\rho\in\mathcal{P}_2(\mathbb{R}^n)$ we write
    \begin{equation*}
        \mathrm{var}\rho=\int_{\mathbb{R}^n}\Big\lVert x-\int_{\mathbb{R}^n}y\, d\rho(y)\Big\rVert^2\, d\rho(x).
    \end{equation*}
    We will also observe in Proposition \ref{pro:optimal} that the latter equals the half of the Schatten $1$-norm of the difference of the covariance matrices of $\mu$ and $\nu$.
\end{remark}

\begin{remark}
    In \cite[Section 4.3]{Bolbotowski2024}, the problem $\mathcal{J}(\mu,\nu)$ is examined in the particular case of $\mu,\nu$ each supported on two points.  The reason for the difference of \cite[Section 4.3, (B)]{Bolbotowski2024}, or \cite[Figure 3, (f)]{Bolbotowski2024}  with \cite[Section 4.3, (A)]{Bolbotowski2024} is that in the first case there is no decomposition of $\mathbb{R}^2$ into two orthogonal one-dimensional subspaces $V_1,V_2$ such that there exists a bimartingale coupling between $\mu$ and $\nu$ with respect to these subspaces. Indeed, this would imply that the support of $(P_{V_1})_{\#}\mu$ has to lie in the convex hull of the support of $(P_{V_1})_{\#}\nu$ and the support of $(P_{V_2})_{\#}\nu$ has to lie in the convex hull of the support of $(P_{V_2})_{\#}\mu$. For the first condition to be satisfied we have to have that the $V_1$ is at least perpendicular to the segment $[z_0,x_2]$ and at most perpendicular to $[y_1,z_0]$. But for such choices of $V_1$  the analogous condition for $V_2$ is not satisfied. This gives a simple argument that allows to infer the observations of \cite[Section 4.3]{Bolbotowski2024}.
\end{remark}

\subsection{Covariance matrices}\label{s:covariance}

The examples in \cite{Bolbotowski2024} suggest that if an optimiser for $\mathcal{I}(\mu,\nu)$ has isometric derivative, then the optimal projections and the subspaces should depend on the difference of the covariance matrices for $\nu$ and for $\mu$. Indeed, we have the following proposition.

In the statement of the proposition below for a subspace $W\subset\mathbb{R}^n$ by $P_W$ we denote the orthogonal projection onto $W$. We write $W_1\perp W_2$ if two subspaces $W_1,W_2$ are orthogonal. Below, $\norm{\cdot}_1$ denotes the Schatten $1$-norm of a matrix, i.e., the sum of all of its singular values. We also note that the Schatten $1$-norm of a matrix $A\in\mathbb{R}^{m\times n}$ -- see e.g. \cite[Remark 2.12]{Ciosmak20212} -- is equal to
\begin{equation}\label{eqn:schatten}
    \norm{A}_1=\sup\{\langle T,A\rangle\mid T\colon \mathbb{R}^n\to\mathbb{R}^m\text{ is a partial isometry}\}.
\end{equation}

\begin{proposition}\label{pro:optimal}
Let $\mu,\nu\in\mathcal{P}_2(\mathbb{R}^n)$ have common barycentre and let
\begin{equation*}
C=\int_{\mathbb{R}^n}xx^* \,d (\nu-\mu)(x)
\end{equation*}
be the covariance matrix of $\nu-\mu$. Let $C=C_1-C_2$ with  $C_1,C_2$ being symmetric positive semi-definite matrices, with $C_1=C P_{V_1}, C_2=-C P_{V_2}$, where $V_1,V_2$ are the sums of positive and negative eigenspaces of $C$. Then
\begin{equation}\label{eqn:superpro}
    \sup\Big\{\int_{\mathbb{R}^n} \frac12\big( \norm{P_{W_1}x}^2-\norm{P_{W_2}x}^2\big)\, d(\nu-\mu)(x)\mid  W_1\oplus W_2=\mathbb{R}^n, W_1\perp W_2\Big\}
\end{equation}
is attained by mutually complementing, orthogonal subspaces $W_1,W_2$ if and only if $W_1\supset V_1$ and $W_2\supset V_2$.
Moreover, the supremum is equal to $\frac12\norm{C}_1$.
\end{proposition}
\begin{proof}
Let $W_1,W_2$ be two mutually complementing, orthogonal subspaces. For an orthogonal projection $P$ and $x\in\mathbb{R}^n$ we have
\begin{equation*}
    \norm{Px}^2=\langle Px,Px\rangle=\langle Px,x\rangle=\mathrm{tr} (Pxx^*).
\end{equation*}
Therefore by linearity of the projections $P_{W_1},P_{W_2}$ and by linearity of the trace we get
\begin{align*}
    &  \int_{\mathbb{R}^n} \frac12\big( \norm{P_{W_1}x}^2-\norm{P_{W_2}x}^2\big)\, d(\nu-\mu)(x)=\\
    &\frac12 \mathrm{tr}\Big( P_{W_1}\int_{\mathbb{R}^n}  xx^*\, d(\nu-\mu)(x) -  P_{W_2}\int_{\mathbb{R}^n} xx^*\, d(\nu-\mu)(x)\Big)=\\
    &\frac12\mathrm{tr} (P_{W_1}-P_{W_2})(C_1-C_2).
\end{align*}
Since $P_{W_1}+P_{W_2}=\mathrm{Id}$, the optimal $P_{W_1}$ maximises
\begin{equation*}
    \mathrm{tr}P_{W_1}(C_1-C_2)
\end{equation*}
and thus it is necessarily a projection onto a subspace containing $V_1$ and perpendicular to $V_2$. Consequently, $W_2$ is  a subspace containing $V_2$ and perpendicular to $V_1$. Therefore, for optimal $P_{W_1}$, and $P_{W_2}$
\begin{equation*}
     \mathrm{tr}P_{W_1}(C_1-C_2)= \mathrm{tr}C_1,  \mathrm{tr}P_{W_2}(C_1-C_2)= -\mathrm{tr}C_2
\end{equation*}
Since $C_1,C_2$ are both symmetric and positive semi-definite, their traces equal the respective sums of the singular values.  Thus
\begin{equation*}
    \mathrm{tr} (P_{W_1}-P_{W_2})(C_1-C_2)=\norm{C}_1.
\end{equation*}
\end{proof}

\begin{proposition}\label{pro:reduction}
  Let $\mu,\nu\in\mathcal{P}_2(\mathbb{R}^n)$ have common barycentre and be such that  there exists  an optimal $u\in\mathcal{C}^{1,1}(\mathbb{R}^n)$ for $\mathcal{I}(\mu,\nu)$ with isometric derivative. Let $\sigma\in\Sigma(\mu,\nu)$ be optimal for $\mathcal{J}(\mu,\nu)$. Let
  \begin{equation*}
      C=\int_{\mathbb{R}^n}xx^*\, d(\nu-\mu)(x)
  \end{equation*}
  be the difference of the covariance matrices of $\nu$ and of $\mu$. Let $V=(\mathrm{ker}C)^{\perp}$ and let $V_1,V_2$ be two mutually complementing, orthogonal subspaces of $V$ -- $V_1$ the sum of positive eigenspaces of $C$ and $V_2$ the sum of negative eigenspaces of $C$.
  There exists a Borel measure $\theta_V\in\mathcal{P}(CC(\mathbb{R}^n))$ concentrated on the set of affine subspaces parallel to $V$, and Borel measures $\mu_{\mathcal{S},V},\nu_{\mathcal{S},V}\in\mathcal{P}_2(\mathbb{R}^n)$, and Borel measures $\sigma_{\mathcal{S},V}\in\mathcal{P}((\mathbb{R}^n)^3)$ defined for $\theta_{V}$-almost every $\mathcal{S}$, such that
\begin{enumerate}[label=(\roman*)]
    \item\label{i:measur} for any Borel set $A\subset\mathbb{R}^n$ the maps
    \begin{equation*}
        CC(\mathbb{R}^n)\ni \mathcal{S}\mapsto \mu_{\mathcal{S},V}(A)\in\mathbb{R}\text{, and }CC(\mathbb{R}^n)\ni \mathcal{S}\mapsto \nu_{\mathcal{S},V}(A)\in\mathbb{R}
    \end{equation*}
    are  Borel measurable, and
    \begin{equation*}
        \mu(A)=\int_{CC(\mathbb{R}^n)}\mu_{\mathcal{S},V}(A)\, d\theta_V(\mathcal{S}), \nu(A)=\int_{CC(\mathbb{R}^n)}\nu_{\mathcal{S},V}(A)\, d\theta_V(\mathcal{S}),
        \end{equation*}
    \item for $\theta_V$-almost every $\mathcal{S}$, $\mu_{\mathcal{S},V}$, $\nu_{\mathcal{S},V}$ are concentrated on $\mathcal{S}$, i.e.,
    \begin{equation*}
    \mu_{\mathcal{S},V}(\mathcal{S}^c)=0,    \nu_{\mathcal{S},V}(\mathcal{S}^c)=0,
    \end{equation*}
        \item\label{i:moments}  for $\theta_V$-almost every $\mathcal{S}$
 \begin{equation*}  \int_{\mathbb{R}^n}x\,d\,\mu_{\mathcal{S},V}(x)=\int_{\mathbb{R}^n}x\,d\nu_{\mathcal{S},V}(x),
  \end{equation*}
\item $
    \sigma_{\mathcal{S},V}\in\Sigma(\mu_{\mathcal{S},V},\nu_{\mathcal{S},V})$ for $\theta_V$-almost every $\mathcal{S}$,
    \item for any Borel set $A\subset(\mathbb{R}^n)^3$ the map
    \begin{equation*}
        CC(\mathbb{R}^n)\ni\mathcal{S}\mapsto \sigma_{\mathcal{S},V}(A)\in\mathbb{R}
    \end{equation*}
    is Borel measurable and
    \begin{equation*}
        \sigma(A)=\int_{CC(\mathbb{R}^n)}\sigma_{\mathcal{S},V}(A)\, d\theta_V(\mathcal{S}),
    \end{equation*}
\item  for $\theta_V$-almost every $\mathcal{S}$, $\sigma_{\mathcal{S},V}$ is concentrated on $\mathcal{S}^3$, that is
\begin{equation*}
    \sigma_{\mathcal{S},V}((\mathcal{S}^3)^c)=0,
\end{equation*}
\item\label{i:optimal}  $u_{V_1,V_2}$ is optimal for  $\mathcal{I}(\mu_{\mathcal{S},V},\nu_{\mathcal{S},V})$ and
 $\sigma_{\mathcal{S},V}$ is optimal for $\mathcal{J}(\mu_{\mathcal{S},V},\nu_{\mathcal{S},V})$.
\end{enumerate}
\end{proposition}
\begin{proof}
    By the assumption, $\mathcal{I}(\mu,\nu)$ admits an optimal $u$ with isometric derivative. By Proposition \ref{pro:cc}, $\mathcal{I}(\mu,\nu)$ is equal to the supremum considered in (\ref{eqn:superpro}). By Proposition \ref{pro:optimal} the map
\begin{equation*}
    u_{V_1,V_2}(x)=\frac12\big(\norm{P_{V_1}x}^2-\norm{P_{V_2}x}^2\big)
\end{equation*}
is optimal for $\mathcal{I}(\mu,\nu)$, since its integral against $\nu-\mu$ equals $\tfrac12\norm{C}_1$. Its derivative at a point $x\in\mathbb{R}^n$ equals $(P_{V_1}-P_{V_2})x$. This is a partial isometry on $\mathbb R^n$, isometric on $V$ and zero on $\ker C$. The leaves of the derivative are precisely the affine spaces that are translations of $V_1\oplus V_2$.

Let $q=P_{\ker C}$. By (\ref{eqn:optimality}) and Lemma \ref{lem:three}, $z-x\in V_1$ and $z-y\in V_2$ for $\sigma$-almost every $(x,y,z)$. Thus $q(x)=q(y)=q(z)$. Apply Theorem \ref{thm:disintegration} to $\sigma$ and the map $(x,y,z)\mapsto q(x)+V$, and denote its pushforward by $\theta_V$. For $\theta_V$-almost every $\mathcal S$, the resulting measure $\sigma_{\mathcal S,V}$ is concentrated on $\mathcal S^3$. By uniqueness of disintegration, its first and second marginals are the disintegrations $\mu_{\mathcal S,V}$ and $\nu_{\mathcal S,V}$ of $\mu$ and $\nu$ by $x\mapsto q(x)+V$. All three marginals have finite second moments almost everywhere, since $\mu,\nu$ have finite second moments and $\sigma$ has finite cost.

Let $g$ be a bounded Borel function of the first variable and $h$ a bounded Borel function of the second. For each Borel set $E\subset\ker C$, test the two identities in (\ref{eqn:mart}) with $g\mathbf1_{q^{-1}(E)}$ and $h\mathbf1_{q^{-1}(E)}$, respectively. Since $q(x)=q(y)$ for $\sigma$-almost every $(x,y,z)$ and $E$ is arbitrary, disintegration gives both identities for $\sigma_{\mathcal S,V}$ almost everywhere, simultaneously for all bounded Borel $g,h$ by first taking indicators of a countable generating algebra. This gives $\sigma_{\mathcal S,V}\in\Sigma(\mu_{\mathcal S,V},\nu_{\mathcal S,V})$ and, on taking $g=h=1$, the common barycentre. Finally, (\ref{eqn:optimality}) holds $\sigma_{\mathcal S,V}$-almost everywhere for $\theta_V$-almost every $\mathcal S$, which gives \ref{i:optimal}.
\end{proof}

\begin{example}\label{exa:furtherreduction}
    There is no guarantee that after performing the decomposition described in Proposition \ref{pro:reduction} the obtained conditional measures $\mu_{\mathcal{S},V},\nu_{\mathcal{S},V}$ will be such that the difference of their covariances is non-degenerate. Let us consider the following example of two measures $\mu,\nu\in\mathcal{P}_2(\mathbb{R}^2)$.
Let
\begin{equation*}
    \mu=\frac14(\delta_{(-1,-1)}+\delta_{(1,-1)}+2\delta_{(0,1)})
\end{equation*}
and let
\begin{equation*}
    \nu=\frac14(\delta_{(-1,-1)}+\delta_{(1,-1)}+\delta_{(-1,1)}+\delta_{(1,1)}).
\end{equation*}
Then $\mu,\nu$ are in convex order, so in particular they are in convex-concave order. Thus by Proposition \ref{pro:cc} the optimiser for $\mathcal{I}(\mu,\nu)$ can be taken to have isometric derivative.  The difference of the covariance matrices
\begin{equation*}
    C=\int_{\mathbb{R}^2}xx^*\, d(\nu-\mu)(x)
\end{equation*}
is degenerate -- for the second standard vector $e_2$ we have
\begin{equation*}
   \langle e_2,Ce_2\rangle=0.
\end{equation*}
Note that $\mathrm{ker}C=\mathbb{R}e_2$, as $C$ is non-zero. Note that the function $ u_{\mathbb{R}e_1,\{0\}}\in\mathcal{C}^{1,1}(\mathbb{R}^2)$ defined by the formula
\begin{equation*}
    u_{\mathbb{R}e_1,\{0\}}(x)=\frac12\langle x,e_1\rangle^2\text{ for }x\in\mathbb{R}^2
\end{equation*}
is therefore an optimiser for $\mathcal{I}(\mu,\nu)$.
However, the conditional measures with respect to the partition induced by $ u_{\mathbb{R}e_1,\{0\}}$, as in Proposition \ref{pro:reduction}, may have degenerate difference of their covariances. Indeed, let
\begin{equation*}
    \mu_1=\delta_{(0,1)}\text{ and }\mu_{-1}=\frac12(\delta_{(-1,-1)}+\delta_{(1,-1)})
\end{equation*}
and
\begin{equation*}
    \nu_1=\frac12(\delta_{(-1,1)}+\delta_{(1,1)})\text{ and }\nu_{-1}=\frac12(\delta_{(-1,-1)}+\delta_{(1,-1)}).
\end{equation*}
Then
\begin{equation*}
    \mu=\frac12\mu_{-1}+\frac12\mu_1\text{ and }\nu=\frac12\nu_{-1}+\frac12\nu_1,
\end{equation*}
is the decomposition of $\mu,\nu$ described by Proposition \ref{pro:reduction}. However, as $\mu_{-1}=\nu_{-1}$, the difference of the covariances for $\nu_{-1}$ and $\mu_{-1}$ is zero. Thus this pair can be further partitioned, unlike the pair $\mu_1,\nu_1$.
\end{example}

\section{Description of optimal plans}\label{s:descr}

The results of Section \ref{s:balance} and Section \ref{s:bimart} allow us to provide a complete description, in terms of bimartingale transports, of all optimal plans for arbitrary pairs of measures $\mu,\nu\in\mathcal{P}_2(\mathbb{R}^n)$ that are absolutely continuous with respect to the Lebesgue measure and with common barycentre for which there is an optimal $\sigma\in\Sigma(\mu,\nu)$ with absolutely continuous third marginal.

We shall employ Theorem \ref{thm:parti} to decompose the problem into problems on which the optimum in the dual problem is attained by a function with isometric derivative. Then Proposition \ref{pro:cc} will show that in this case the conditional measures are in convex-concave order with respect to some mutually complementing, orthogonal subspaces $V_1,V_2$ of the tangent spaces of the corresponding leaves. Proposition  \ref{pro:optimal}  will allow to describe these subspaces through the covariance matrices of the conditional measures. Theorem \ref{thm:bimart} will show that any bimartingale coupling of the conditional measures on the leaves will give rise, thanks to Proposition \ref{pro:bimart}, to an optimal plan.

\begin{remark}
    The requirement in Theorem \ref{thm:descr} that the difference of the covariance matrices $C_{\mathcal{S}}$ is non-degenerate on the tangent space $V(\mathcal{S})$ is a necessary condition to have the subspaces $V(\mathcal{S})_1,V(\mathcal{S})_2$ uniquely determined. Note that Proposition \ref{pro:optimal} shows that there is an ambiguity in the choice of optimal $W_1,W_2$. The ambiguity occurs exactly when $C_{\mathcal{S}}$ is degenerate.
\end{remark}

\subsection{Refinement of partition}\label{s:refine}

We start with defining the partition with respect to which we shall be disintegrating our measures.

Let $u\in\mathcal{C}^{1,1}(\mathbb{R}^n)$ have $1$-Lipschitz derivative.  In Section \ref{s:leaves} we have defined a partitioning of $\mathbb{R}^n$, up to a set of the Lebesgue measure zero, into maximal sets $\mathcal{S}$ on which $Du$ is an isometry.

Suppose that an optimiser of $\mathcal{I}(\mu,\nu)$ for a pair $\mu,\nu\in\mathcal{P}_2(\mathbb{R}^n)$ with common barycentre is an isometry. In Section \ref{s:covariance} we have shown that if the difference of the covariance matrices of $\nu,\mu$ is degenerate, this partition can be refined, and any optimal plan has to occur within this refinement.
However, as shown in Example \ref{exa:furtherreduction} it is possible that this partition is not the finest one. The idea is thus to take consecutive partitions, each as in Proposition \ref{pro:reduction}. Note the leaves of these partitions are convex sets. Moreover, when performing the partition, if the difference of the covariance matrices is degenerate, then the dimensions of these leaves drops by at least one.  Therefore, the number of partitions that we will need to do will be limited by the dimension of the ambient space.

Let us recall that the set $B(Du)$ denotes the set of points in $\mathbb{R}^n$ that belong to at least two distinct leaves of $Du$. By \cite[Corollary 2.15]{Ciosmak2021} it is contained in a Borel set $N(Du)$ on which $Du$ is not differentiable. By the Rademacher theorem $N(Du)$ is of the Lebesgue measure zero.

We also recall that $CL(\mathbb{R}^n)$ and its subspace $CC(\mathbb{R}^n)$ are defined as the sets of all non-empty closed and all non-empty closed and convex sets, respectively.

Let us fix two measures $\mu,\nu\in\mathcal{P}_2(\mathbb{R}^n)$ with common barycentre and an optimiser $u\in\mathcal{C}^{1,1}(\mathbb{R}^n)$ of $\mathcal{I}(\mu,\nu)$. Throughout this subsection we work under the assumptions of Theorem \ref{thm:descr}. In particular, we fix an optimal plan $\sigma_0\in\Sigma(\mu,\nu)$ whose third marginal is absolutely continuous. The leaves of $Du$ induce the map
\begin{equation*}
    \mathcal{S}_0\colon\mathbb{R}^n\to CC(\mathbb{R}^n)
\end{equation*}
by the formula
\begin{equation*}
    \mathcal{S}_0(x)=\{x\}\text{ if }x\in N(Du),
\end{equation*}
and
\begin{equation*}
    \mathcal{S}_0(x)\text{ is the unique leaf of }Du\text{ that contains }x.
\end{equation*}
Let $\theta_0\in\mathcal{P}(CC(\mathbb{R}^n))$ and maps
\begin{equation*}
    CC(\mathbb{R}^n)\ni\mathcal{S}_0\mapsto (\mu_{\mathcal{S}_0},\nu_{\mathcal{S}_0})\in\mathcal{P}_2(\mathbb{R}^n)\times\mathcal{P}_2(\mathbb{R}^n)
\end{equation*}
be as in Theorem \ref{thm:partim} and Theorem \ref{thm:parti}. For a leaf $\mathcal{S}_0\in CC(\mathbb{R}^n)$ the corresponding matrix
\begin{equation*}
    C_{\mathcal{S}_0}=\int_{\mathbb{R}^n}xx^*\, d(\nu_{\mathcal{S}_0}-\mu_{\mathcal{S}_0})(x)
\end{equation*}
is symmetric, but possibly degenerate.

We define a new map
\begin{equation*}
    \mathcal{S}_1\colon\mathbb{R}^n\to CC(\mathbb{R}^n)
\end{equation*}
by the formula
\begin{equation*}
    \mathcal{S}_1(x)=\mathcal{S}_0(x)\cap \big(x+(\ker C_{\mathcal{S}_0(x)})^{\perp}\big).
\end{equation*}
Clearly, if $\mathcal{S}_0(x)=\mathcal{S}_0(y)$ and $\mathcal{S}_1(x)\cap\mathcal{S}_1(y)\neq\emptyset$, then $\mathcal{S}_1(x)=\mathcal{S}_1(y)$, for each $x,y\in\mathbb{R}^n$.

Let $X$ be a measurable space. In \cite{Beer1994} it is proven that a map $f\colon X\to CL(\mathbb{R}^n)$ is measurable if and only if it is measurable as a multifunction. The latter is defined by the condition that for any open set $U\subset\mathbb{R}^n$ the set
\begin{equation*}
\{x\in X\mid f(x)\cap U\neq \emptyset\}
\end{equation*}
is measurable in $X$.

\begin{lemma}\label{lem:leafmeasurability}
    The map
    \begin{equation*}
    \mathcal{S}_1\colon\mathbb{R}^n\to CC(\mathbb{R}^n)
\end{equation*}
is Borel measurable.
\end{lemma}
\begin{proof}
Recall that $CC(\mathbb{R}^n)$ is equipped with the Wijsman topology. Measurability with respect to the Wijsman topology is equivalent to measurability as a multifunction. That is, we need to show that for any open set $U\subset\mathbb{R}^n$ the set
\begin{equation*}
    \big\{x\in\mathbb{R}^n\mid \mathcal{S}_1(x)\cap U\neq \emptyset\big\}
\end{equation*}
is Borel measurable. In \cite[Theorem 5.1]{Ciosmak2021} it is shown that $\mathcal{S}_0$ is measurable. By \cite[Lemma 18.4, 3., p. 594]{Aliprantis2006}\footnote{The compactness assumption can be relaxed, since $\mathbb{R}^n$ is $\sigma$-compact, cf. proof of \cite[Theorem 5.1]{Ciosmak2021}.}, it suffices to show that the map
\begin{equation}\label{eqn:val}
    x\mapsto x+(\mathrm{ker}C_{\mathcal{S}_0(x)})^{\perp}
\end{equation}
is Borel measurable.

Note that by property \ref{i:mees} of Theorem \ref{thm:partim} and by standard measure-theoretic arguments, for any vector $v\in\mathbb{R}^n$, the map
\begin{equation*}
  CC(\mathbb{R}^n) \ni \mathcal{S}_0\mapsto\int_{\mathbb{R}^n}y\langle y,v\rangle\, d(\nu_{\mathcal{S}_0}-\mu_{\mathcal{S}_0})(y)\in\mathbb{R}^n
\end{equation*}
is Borel measurable.
Thus, so is the map
\begin{equation*}
  \mathbb{R}^n \ni x\mapsto\int_{\mathbb{R}^n}y\langle y,v\rangle\, d(\nu_{\mathcal{S}_0(x)}-\mu_{\mathcal{S}_0(x)})(y)\in\mathbb{R}^n,
\end{equation*}
as a composition of two Borel measurable maps.  As $C_{\mathcal{S}_0}$ is symmetric, the subspace $(\mathrm{ker}C_{\mathcal{S}_0})^{\perp}$ is the range of the map
\begin{equation*}
    \mathbb{R}^n\ni v\mapsto\int_{\mathbb{R}^n}y\langle y,v\rangle\, d(\nu_{\mathcal{S}_0}-\mu_{\mathcal{S}_0})(y)\in\mathbb{R}^n.
\end{equation*}
Let us define the function $\phi\colon   \mathbb{R}^n\times\mathbb{R}^n\to\mathbb{R}^n$ by the formula
\begin{equation*}
  \mathbb{R}^n\times\mathbb{R}^n \ni (x,v)\mapsto x+\int_{\mathbb{R}^n}y\langle y,v\rangle\, d(\nu_{\mathcal{S}_0(x)}-\mu_{\mathcal{S}_0(x)})(y)\in\mathbb{R}^n.
\end{equation*}
Note that it is a Carath\'eodory function, i.e., it is Borel measurable in the first variable and continuous in the second variable.

Let $U\subset\mathbb{R}^n$ be open. The condition that
\begin{equation*}
    \big(x+(\mathrm{ker}C_{\mathcal{S}_0(x)})^{\perp}\big)\cap U\neq \emptyset
\end{equation*}
is equivalent to the existence of $v\in\mathbb{R}^n$ such that
\begin{equation*}
\phi(x,v)\in U.
\end{equation*}
If $(v_i)_{i=1}^{\infty}$ is a countable dense subset of $\mathbb{R}^n$, then
\begin{equation*}
   \bigcup\Big\{ \{x\in\mathbb{R}^n\mid \phi(x,v)\in U\}\mid v\in\mathbb{R}^n\Big\}=\bigcup_{i=1}^{\infty}   \{x\in\mathbb{R}^n\mid \phi(x,v_i)\in U\}.
\end{equation*}
Each of the sets $ \{x\in\mathbb{R}^n\mid \phi(x,v_i)\in U\}$ is Borel measurable. This shows that (\ref{eqn:val}) is measurable as a multifunction.
 Thus, so is $\mathcal{S}_1$.
\end{proof}

\subsection{Proof of Theorem \ref{thm:descr}}\label{s:proof}

\begin{proof}[Proof of Theorem \ref{thm:descr}]
Fix an optimiser $u\in\mathcal{C}^{1,1}(\mathbb{R}^n)$ of
$\mathcal{I}(\mu,\nu)$ and an optimal plan
$\sigma_0\in\Sigma(\mu,\nu)$ with absolutely continuous third marginal.

We now construct the conditional problems simultaneously with the
partitions. Here and below, by a Borel disintegration we mean that the evaluation on
every Borel set of conditional measures is a Borel function of the element of the partition.

At level $k=0$, Theorems \ref{thm:partim} and \ref{thm:parti} give
\begin{equation*}
 \theta_0=(\mathcal S_0)_{\#}\mu=(\mathcal S_0)_{\#}\nu,
 \quad
 \mu=\int\mu_{\mathcal S_0}\,d\theta_0(\mathcal S_0),
 \quad
 \nu=\int\nu_{\mathcal S_0}\,d\theta_0(\mathcal S_0).
\end{equation*}
For $\theta_0$-almost every $\mathcal S_0$, the conditional measures are
concentrated on $\mathcal S_0$, have a common barycentre, and $u$ is
optimal for
$\mathcal I(\mu_{\mathcal S_0},\nu_{\mathcal S_0})$, with $Du$
isometric on $\mathcal S_0$. Moreover, every optimal
$\sigma\in\Sigma(\mu,\nu)$ has a Borel disintegration
\begin{equation*}
 \sigma=\int\sigma_{\mathcal S_0}\,d\theta_0(\mathcal S_0),
\end{equation*}
where
$\sigma_{\mathcal S_0}\in
\Sigma(\mu_{\mathcal S_0},\nu_{\mathcal S_0})$ is concentrated on
$\mathcal S_0^3$ and is conditionally optimal.

Suppose that these objects have been constructed at level $k<n$. Set
\begin{equation}\label{eqn:leafcovariance}
 C_{\mathcal S_k(x)}
 =\int_{\mathbb R^n}yy^*\,
 d(\nu_{\mathcal S_k(x)}-\mu_{\mathcal S_k(x)})(y)
\end{equation}
and define
\begin{equation*}
 \mathcal S_{k+1}(x)
 =\mathcal S_k(x)\cap
 \Bigl(x+(\ker C_{\mathcal S_k(x)})^\perp\Bigr).
\end{equation*}
Then $x\mapsto C_{\mathcal S_k(x)}$ is a Borel function. Hence, as in the proof of Lemma
\ref{lem:leafmeasurability}, the map $\mathcal S_{k+1}$ is Borel
measurable.

For $\theta_k$-almost every $\mathcal S_k$, Proposition
\ref{pro:reduction}, applied on the affine hull of $\mathcal S_k$, shows
that the pushforwards of $\mu_{\mathcal S_k}$ and
$\nu_{\mathcal S_k}$ under $\mathcal S_{k+1}$ coincide. It also shows
that every $\sigma_{\mathcal S_k}$ is concentrated on triples belonging
to the same element of the refinement. Consequently,
\begin{equation*}
 \theta_{k+1}
 =(\mathcal S_{k+1})_{\#}\mu
 =(\mathcal S_{k+1})_{\#}\nu.
\end{equation*}
\phantomsection\label{par:refinedconditionals}Disintegrating $\mu,\nu$, and $\sigma$ with respect to
$\mathcal S_{k+1}$ gives Borel conditional measures
$\mu_{\mathcal S_{k+1}},\nu_{\mathcal S_{k+1}}$, and
$\sigma_{\mathcal S_{k+1}}$. Since $\mathcal S_{k+1}$ refines
$\mathcal S_k$, uniqueness of disintegration shows that these are
obtained by disintegrating the corresponding conditional measures at
level $k$.

Proposition \ref{pro:reduction} gives the support, common-barycentre and
martingale properties at level $k+1$. The contact equality for $u$ and
$\sigma_{\mathcal S_k}$ passes to the conditional plans; hence $u$ and
$\sigma_{\mathcal S_{k+1}}$ are conditionally optimal. Finally, $Du$
remains isometric because
$\mathcal S_{k+1}\subseteq\mathcal S_k$.

Finally, support on $\mathcal{S}_k(x)$ and the common barycentre condition
imply that $(\mathrm{ker}C_{\mathcal{S}_k(x)})^{\perp}$ is contained in the
tangent space to $\mathcal{S}_k(x)$. Indeed, if $v\in V(\mathcal{S}_k(x))^{\perp}$ then
\begin{equation*}
C_{\mathcal{S}_k(x)}v=\int_{\mathbb{R}^n}(y-x)\langle y-x,v\rangle\,d(\nu_{\mathcal{S}_k(x)}-\mu_{\mathcal{S}_k(x)})=0
\end{equation*}
since $\mu_{\mathcal{S}_k(x)},\nu_{\mathcal{S}_k(x)}$ are concentrated on $\mathcal{S}_k(x)$.
If $C_{\mathcal{S}_k}$ is non-degenerate on that tangent space, then
$(\ker C_{\mathcal{S}_k})^{\perp}$ is the whole tangent space and no
refinement occurs. Otherwise its dimension is strictly smaller, and every set
\begin{equation*}
 \mathcal{S}_{k+1}(x)=\mathcal{S}_k(x)\cap \bigl(x+(\ker C_{\mathcal{S}_k(x)})^{\perp}\bigr)
\end{equation*}
has strictly smaller dimension. Along a chain of proper refinements the
dimension decreases at every step, so after at most $n$ steps the covariance
is non-degenerate on the tangent space; in dimension zero this is understood
vacuously.

For the remainder of the proof, we will drop the subscript $n$ for the partition at level $n$. We get therefore
\begin{equation*}
 \theta=\mathcal{S}_{\#}\mu=\mathcal{S}_{\#}\nu,
 \qquad
 \mu=\int\mu_{\mathcal{S}}\,d\theta(\mathcal{S}),
 \qquad
 \nu=\int\nu_{\mathcal{S}}\,d\theta(\mathcal{S}).
\end{equation*}
It gives directly \ref{i:concen}, \ref{i:bary}, \ref{i:covar}, and
\ref{i:set} of Theorem \ref{thm:descr}.

For $\theta$-almost every $\mathcal{S}$, a quadratic affine map is optimal for $\mathcal{I}(\mu_{\mathcal{S}},\nu_{\mathcal{S}})$ and has isometric
derivative. Proposition \ref{pro:cc} therefore gives $\mu_{\mathcal{S}}\prec_{c-c}\nu_{\mathcal{S}}$.
Since $C_{\mathcal{S}}$ is non-degenerate on the tangent space
$V(\mathcal{S})$, Proposition \ref{pro:optimal} identifies
$V(\mathcal{S})_1$ and $V(\mathcal{S})_2$ with its positive and negative
eigenspaces.

Let now $\sigma\in\Sigma(\mu,\nu)$ be optimal. Then by the above induction
\begin{equation*}
 \sigma=\int\sigma_{\mathcal{S}}\,d\theta(\mathcal{S}),
\end{equation*}
where $\sigma_{\mathcal{S}}$ is concentrated on $\mathcal{S}^3$ and is
optimal in
$\Sigma(\mu_{\mathcal{S}},\nu_{\mathcal{S}})$. Setting
$\pi_{\mathcal{S}}=\mathrm{P}_{12}\sigma_{\mathcal{S}}$, Theorem
\ref{thm:bimart} gives
\begin{equation*}
 \pi_{\mathcal{S}}\in
 \Gamma_{bm}(\mu_{\mathcal{S}},\nu_{\mathcal{S}},
 V(\mathcal{S})_1,V(\mathcal{S})_2),
 \quad
 \sigma_{\mathcal{S}}=(R_{\mathcal{S}})_{\#}\pi_{\mathcal{S}}.
\end{equation*}
This proves the necessity in (\ref{eqn:sigmadis}) and
\ref{i:displan}--\ref{i:rpush}.

Conversely, if we suppose that $\sigma$ is defined as in (\ref{eqn:sigmadis}), with
\begin{equation*}
    CC(\mathbb{R}^n)\ni\mathcal{S}\mapsto \pi_{\mathcal{S}}\in\mathcal{P}(\mathbb{R}^n\times\mathbb{R}^n)
\end{equation*}
satisfying \ref{i:displan}-\ref{i:rpush}. Then $\sigma\in\mathcal{P}((\mathbb{R}^n)^3)$, as, cf. Lemma \ref{lem:vmeasur}, the map $R_{\mathcal{S}}$ is Borel measurable. Here we pick a Borel measurable map $z\colon \mathcal{S}(\mathbb{R}^n)\to\mathbb{R}^n$ such that $z(\mathcal{S})\in\mathcal{S}$. Such a map exists -- it suffices to  compose the map of Lemma \ref{lem:selection} with the map $\mathcal{S}\mapsto(\mathcal{S},\mathcal{S})$. Moreover, the fact that $\sigma\in\Sigma(\mu,\nu)$ follows by \ref{i:set}, \ref{i:rpush} and Proposition \ref{pro:bimart}.
Note that for $\theta$-almost every $\mathcal{S}$ the measure $\sigma_{\mathcal{S}}$ is optimal for $\mathcal{J}(\mu_{\mathcal{S}},\nu_{\mathcal{S}})$, as $\pi_{\mathcal{S}}$ is a bimartingale coupling, see Proposition \ref{pro:bimart}. If  $u\in\mathcal{C}^{1,1}(\mathbb{R}^n)$ is the chosen optimiser, then for $\theta$-almost every leaf $\mathcal{S}$ we have
\begin{equation*}
    \int_{\mathbb{R}^n}u\, d(\nu_{\mathcal{S}}-\mu_{\mathcal{S}})=\int_{(\mathbb{R}^n)^3}c\, d\sigma_{\mathcal{S}}.
\end{equation*}
Integrating out these equalities with respect to $\theta$, taking into account \ref{i:set} and the definition of $\sigma$, yields that $\sigma$ is indeed optimal for $\mathcal{J}(\mu,\nu)$.

To prove formula (\ref{eqn:cost}) it suffices to show that
\begin{equation*}
 \mathcal{J}(\mu_{\mathcal{S}},\nu_{\mathcal{S}})=\frac12\norm{C_{\mathcal{S}}}_1\text{ for }\theta\text{-almost every }\mathcal{S}.
\end{equation*}
This is, however, the last assertion of Proposition \ref{pro:optimal}.

The proof is complete.

\end{proof}

\section{General absolutely continuous measures}\label{s:general}

Let us consider now the case of two  absolutely continuous probability measures $\mu,\nu\in\mathcal{P}_2(\mathbb{R}^n)$ with common barycentre, without the additional assumption of the existence of an optimal $\sigma\in\Sigma(\mu,\nu)$ with absolutely continuous third marginal.
The main point of the current section is that the partition into leaves of $Du$ is still relevant also to the solution of the problem. Moreover, optimal coupling $\pi$ of $\mu$ and of $\nu$ can only assign mass to neighbouring leaves of $Du$ and then $\sigma$ is a composition of plans charging only the triples in the neighbouring leaves and in their interfaces. Here we say that two leaves $\mathcal{S}_1,\mathcal{S}_2$ are neighbouring if $\mathcal{S}_1\cap\mathcal{S}_2\neq \emptyset$.

\begin{theorem}\label{thm:noassu}
    Let $\mu,\nu\in\mathcal{P}_2(\mathbb{R}^n)$ be absolutely continuous with respect to the Lebesgue measure and have common barycentre.
    Then there exist
    \begin{enumerate}
        \item a partition of $\mathbb{R}^n$, up to a set of the Lebesgue measure zero, into closed and convex sets $\mathcal{S}\in CC(\mathbb{R}^n)$, such that for each element  $\mathcal{S}$ of the partition its tangent space $V(\mathcal{S})=V(\mathcal{S})_1\oplus V(\mathcal{S})_2$ can be decomposed into sum of two mutually orthogonal subspaces $V(\mathcal{S})_1,V(\mathcal{S})_2$,
        \item  a Borel measurable map $z\colon (\mathcal{S}(\mathbb{R}^n))^2\to\mathbb{R}^n$, such that
        \begin{equation*}
            z(\mathcal{S}_1,\mathcal{S}_2)\in\mathcal{S}_1\cap\mathcal{S}_2
        \end{equation*}
        for intersecting $\mathcal{S}_1,\mathcal{S}_2\in \mathcal{S}(\mathbb{R}^n)$, and maps
        \begin{equation*}
            R_{\mathcal{S}_1,\mathcal{S}_2}\colon(\mathbb{R}^n)^2\to(\mathbb{R}^n)^3,
        \end{equation*}
        defined for all intersecting $\mathcal{S}_1,\mathcal{S}_2\in\mathcal{S}(\mathbb{R}^n)$, and
        given by the formulae
        \begin{equation*}
            R_{\mathcal{S}_1,\mathcal{S}_2}(x,y)=\Big(x,y,P_{V(\mathcal{S}_1)_2}\big(x-z(\mathcal{S}_1,\mathcal{S}_2)\big)+P_{V(\mathcal{S}_2)_1}\big(y-z(\mathcal{S}_1,\mathcal{S}_2)\big)+z(\mathcal{S}_1,\mathcal{S}_2)\Big)
        \end{equation*}
     for $x,y\in\mathbb{R}^n$,
        \item Borel probability measures $\theta_1,\theta_2\in\mathcal{P}(CC(\mathbb{R}^n))$ and Borel probability measures
    \begin{equation*}
        \mu_{\mathcal{S}},\nu_{\mathcal{S}}\in\mathcal{P}_2(\mathbb{R}^n) \text{ defined for all }\mathcal{S}\in CC(\mathbb{R}^n)
    \end{equation*}
    \end{enumerate}
    such that:
    \begin{enumerate}[label=(\roman*)]
        \item  $\mu_{\mathcal{S}}$ is concentrated on $\mathcal{S}$  for $\theta_1$-almost every $\mathcal{S}$ and $\nu_{\mathcal{S}}$  is concentrated on $\mathcal{S}$  for $\theta_2$-almost every $\mathcal{S}$.
        \item\label{i:set2} for any Borel set $A\subset\mathbb{R}^n$ the map
    \begin{equation*}
        CC(\mathbb{R}^n)\ni\mathcal{S}\mapsto (\mu_{\mathcal{S}}(A), \nu_{\mathcal{S}}(A))\in\mathbb{R}\times\mathbb{R}
    \end{equation*}
    is Borel measurable
     \begin{equation*}
        \mu(A)=\int_{CC(\mathbb{R}^n)}\mu_{\mathcal{S}}(A)\, d\theta_1(\mathcal{S})\text{ and }        \nu(A)=\int_{CC(\mathbb{R}^n)}\nu_{\mathcal{S}}(A)\, d\theta_2(\mathcal{S}),
    \end{equation*}
    \item  a plan $\sigma\in \mathcal{P}((\mathbb{R}^n)^3)$
    belongs to $ \Sigma(\mu,\nu)$  and is optimal for $\mathcal{J}(\mu,\nu)$ if and only if it is of the form
    \begin{equation}\label{eqn:sigmadis2}
        \sigma=\int_{(CC(\mathbb{R}^n))^2}(R_{\mathcal{S}_1,\mathcal{S}_2})_{\#}\pi_{\mathcal{S}_1,\mathcal{S}_2}\, d\theta(\mathcal{S}_1,\mathcal{S}_2),
    \end{equation}
    for some map
    \begin{equation*}
        (CC(\mathbb{R}^n))^2\ni (\mathcal{S}_1,\mathcal{S}_2)\mapsto\pi_{\mathcal{S}_1,\mathcal{S}_2}\in \mathcal{P}(\mathbb{R}^n\times\mathbb{R}^n)
    \end{equation*}
    and some $\theta\in\mathcal{P}((CC(\mathbb{R}^n))^2)$ such that:
    \begin{enumerate}[label=(\alph*)]
    \item\label{i:displan2} for any Borel set $B\subset\mathbb{R}^n\times\mathbb{R}^n$ the map
    \begin{equation*}
       ( CC(\mathbb{R}^n))^2\ni (\mathcal{S}_1,\mathcal{S}_2)\mapsto \pi_{\mathcal{S}_1,\mathcal{S}_2}(B)\in\mathbb{R}
    \end{equation*}
    is Borel measurable,
\item\label{i:coupling} $\theta$ is a coupling of $\theta_1$ and of $\theta_2$, that is concentrated on the set of pairs $(\mathcal{S}_1,\mathcal{S}_2)$ such that $\mathcal{S}_1\cap\mathcal{S}_2\neq\emptyset$,
\item\label{i:munu} for $\theta_1$-almost every $\mathcal{S}_1$
\begin{equation}\label{eqn:mudis}
    \mu_{\mathcal{S}_1}=\int_{CC(\mathbb{R}^n)}\mathrm{P}_1\pi_{\mathcal{S}_1,\mathcal{S}_2}\, d\theta_{\mathcal{S}_1}(\mathcal{S}_2),
\end{equation}
and
for $\theta_2$-almost every $\mathcal{S}_2$
\begin{equation}\label{eqn:nudis}
    \nu_{\mathcal{S}_2}=\int_{CC(\mathbb{R}^n)}\mathrm{P}_2\pi_{\mathcal{S}_1,\mathcal{S}_2}\, d\theta^{\mathcal{S}_2}(\mathcal{S}_1),
\end{equation}
where $\theta_{\mathcal{S}},\theta^{\mathcal{S}}$ are probability measures in $\mathcal{P}(CC(\mathbb{R}^n))$, defined for each $S\in CC(\mathbb{R}^n)$ such that for any Borel set $C\subset CC(\mathbb{R}^n)$ the map
\begin{equation*}
    CC(\mathbb{R}^n)\ni \mathcal{S}\mapsto (\theta_{\mathcal{S}}(C),\theta^{\mathcal{S}}(C))\in\mathbb{R}^2
\end{equation*}
is Borel measurable, and for any Borel sets $C,D\subset  CC(\mathbb{R}^n)$
\begin{equation}\label{eqn:thetas}
    \theta(C\times D)=\int_{ C}\theta_{\mathcal{S}}(D)\, d\theta_1(\mathcal{S}),    \theta(C\times D)=\int_{ D}\theta^{\mathcal{S}}(C)\, d\theta_2(\mathcal{S}),
\end{equation}
 \item\label{i:bimartplan2} the martingale condition is satisfied: for all bounded, Borel functions $g,h\colon\mathbb{R}^n\to\mathbb{R}$
            \begin{equation*}
P_{V(\mathcal{S}_1)_1}\Bigg( \int_{CC(\mathbb{R}^n)}\bigg(\int_{(\mathbb{R}^n)^2}\Big(P_{V(\mathcal{S}_2)_1}(y-z(\mathcal{S}_1,\mathcal{S}_2))
-(x-z(\mathcal{S}_1,\mathcal{S}_2))\Big)g(x)\, d\pi_{\mathcal{S}_1,\mathcal{S}_2}(x,y)\bigg)\,d\theta_{\mathcal{S}_1}(\mathcal{S}_2)\Bigg)=0,
\end{equation*}
for $\theta_1$-almost every $\mathcal{S}_1$,
and
\begin{equation*}
P_{V(\mathcal{S}_2)_2} \Bigg(\int_{CC(\mathbb{R}^n)} \bigg( \int_{(\mathbb{R}^n)^2}\Big(P_{V(\mathcal{S}_1)_2}(x-z(\mathcal{S}_1,\mathcal{S}_2))
-(y-z(\mathcal{S}_1,\mathcal{S}_2))\Big)h(y)\, d\pi_{\mathcal{S}_1,\mathcal{S}_2}(x,y)\bigg)\,d\theta^{\mathcal{S}_2}(\mathcal{S}_1)\Bigg)=0.
\end{equation*}
for $\theta_2$-almost every $\mathcal{S}_2$,
 \item\label{i:vv} the triples $(x,y,z)$ are optimal: for $\theta$-almost every $(\mathcal{S}_1,\mathcal{S}_2)$ and $\pi_{\mathcal{S}_1,\mathcal{S}_2}$-almost every $(x,y)$
    \begin{equation*}
   P_{V(\mathcal{S}_2)_1} P_{V(\mathcal{S}_1)_2}(x-z(\mathcal{S}_1,\mathcal{S}_2))=0,
   \end{equation*}
   and
   \begin{equation*}
  P_{V(\mathcal{S}_1)_2} P_{V(\mathcal{S}_2)_1}(y-z(\mathcal{S}_1,\mathcal{S}_2))=0.
\end{equation*}
    \end{enumerate}
\end{enumerate}
\end{theorem}

\begin{lemma}\label{lem:selection}
    Let
    \begin{equation*}
        U=\{(\mathcal{S}_1,\mathcal{S}_2)\in (CC(\mathbb{R}^n))^2\mid \mathcal{S}_1\cap\mathcal{S}_2\neq\emptyset\}.
    \end{equation*}
    There is a Borel measurable map $z\colon U\to\mathbb{R}^n$ such that
    \begin{equation*}
        z(\mathcal{S}_1,\mathcal{S}_2)\in \mathcal{S}_1\cap\mathcal{S}_2
    \end{equation*}
    for all $(\mathcal{S}_1,\mathcal{S}_2)\in  U$.
\end{lemma}
\begin{proof}
    We shall use the Kuratowski--Ryll-Nardzewski theorem, see, e.g., \cite[Theorem 18.13, p. 600]{Aliprantis2006}. To this aim we recall that $CC(\mathbb{R}^n)$ is a Polish space, see \cite{Beer1991}. We need to verify the assignment
    \begin{equation*}
        U \ni (\mathcal{S}_1,\mathcal{S}_2)\mapsto\mathcal{S}_1\cap\mathcal{S}_2\in CC(\mathbb{R}^n)
    \end{equation*}
    is weakly measurable. It is immediate that the assignments
       \begin{equation*}
U \ni (\mathcal{S}_1,\mathcal{S}_2)\mapsto\mathcal{S}_i\in CC(\mathbb{R}^n)
    \end{equation*}
    for $i=1,2$ are weakly measurable. For $r>0$ we set
    \begin{equation*}
              U_r=\{(\mathcal{S}_1,\mathcal{S}_2)\in (CC(\mathbb{R}^n))^2\mid \mathcal{S}_1\cap\mathcal{S}_2\cap B(0,r)\neq\emptyset\},
    \end{equation*}
    where $B(0,r)$ is the closed ball of radius $r$ centred at $0$. Now,  \cite[Lemma 18.4(3), p. 594]{Aliprantis2006} together with \cite[Theorem 18.13, p. 600]{Aliprantis2006} show that $U_r$ is Borel measurable and that there is a Borel measurable selection
    \begin{equation*}
        z_r\colon U_r \to\mathbb{R}^n
    \end{equation*}
    such that $z_r(\mathcal{S}_1,\mathcal{S}_2)\in \mathcal{S}_1\cap\mathcal{S}_2\cap B(0,r)$ for each $(\mathcal{S}_1,\mathcal{S}_2)\in U_r$. Since
    \begin{equation*}
        U=\bigcup\{U_r\mid r\in\mathbb{N}\},
    \end{equation*}
    we can define the selection
        \begin{equation*}
        z\colon U \to\mathbb{R}^n
    \end{equation*}
    by the formula
    \begin{equation*}
        z=\sum_{r=1}^{\infty}z_r\mathbf{1}_{U_{r}\setminus U_{r-1}},
    \end{equation*}
where we put $U_0=\emptyset$.
\end{proof}

\begin{lemma}\label{lem:vmeasur}
 Let $u\in\mathcal{C}^{1,1}(\mathbb{R}^n)$ have $1$-Lipschitz derivative.   The maps
 \begin{equation*}
    \{ (x,y)\in (\mathbb{R}^n)^2\mid \mathcal{S}(x)\cap\mathcal{S}(y)\neq\emptyset\}\ni(x,y)\mapsto P_{V(\mathcal{S}(x))_1}\big(x-z(\mathcal{S}(x),\mathcal{S}(y))\big)\in\mathbb{R}^n,
 \end{equation*}
 and
  \begin{equation*}
    \{ (x,y)\in (\mathbb{R}^n)^2\mid \mathcal{S}(x)\cap\mathcal{S}(y)\neq\emptyset\}\ni(x,y)\mapsto P_{V(\mathcal{S}(y))_2}\big(y-z(\mathcal{S}(x),\mathcal{S}(y))\big)\in\mathbb{R}^n,
 \end{equation*}
are Borel measurable.
\end{lemma}
\begin{proof}
    Let us recall that, for $x\notin N(Du)$, the set $\mathcal{S}(x)$ is the unique leaf containing $x$, and $\mathcal{S}(x)=\{x\}$ for $x\in N(Du)$. Then $V(\mathcal{S})$ is the tangent space to $\mathcal{S}$ and
    \begin{equation*}
        V(\mathcal{S})=V(\mathcal{S})_1\oplus V(\mathcal{S})_2
    \end{equation*}
    is its orthogonal decomposition into subspaces $V(\mathcal{S})_1, V(\mathcal{S})_2$ such that
    \begin{equation*}
        Du(x)-Du(y)=P_{V(\mathcal{S})_1}(x-y)-P_{V(\mathcal{S})_2}(x-y)\text{ for }x,y\in\mathcal{S},
    \end{equation*}
    cf. Lemma \ref{lem:isoderma}.
    Moreover, for $(x,y)\in \{ (x,y)\in (\mathbb{R}^n)^2\mid \mathcal{S}(x)\cap\mathcal{S}(y)\neq\emptyset\}$
    \begin{equation*}
        z(\mathcal{S}(x),\mathcal{S}(y))\in\mathcal{S}(x)\cap\mathcal{S}(y),
    \end{equation*}
    where $z$ is as in Lemma \ref{lem:selection} applied for $v=Du$.
    Note now that
    \begin{equation*}
2P_{V(\mathcal{S}(x))_1}\big(x-z(\mathcal{S}(x),\mathcal{S}(y))\big)=Du(x)-Du(z(\mathcal{S}(x),\mathcal{S}(y)))+\big(x-z(\mathcal{S}(x),\mathcal{S}(y))\big)
    \end{equation*}
    and
        \begin{equation*}
2P_{V(\mathcal{S}(y))_2}\big(y-z(\mathcal{S}(x),\mathcal{S}(y))\big)=-\big(Du(y)-Du(z(\mathcal{S}(x),\mathcal{S}(y)))\big)+\big(y-z(\mathcal{S}(x),\mathcal{S}(y))\big).
    \end{equation*}
This immediately implies the assertion of the lemma.
\end{proof}

\begin{proof}[Proof of Theorem \ref{thm:noassu}]
  Let $u\in\mathcal{C}^{1,1}(\mathbb{R}^n)$ be an optimiser of $\mathcal{I}(\mu,\nu)$.  Since both $\mu$ and $\nu$ are absolutely continuous, the set of points that belong to at least two distinct leaves of $Du$ is of the Lebesgue measure zero, see \cite[Corollary 2.15]{Ciosmak2021}, and thanks to the measurability of the map $x\mapsto \mathcal{S}(x)$, see \cite[Section 5]{Ciosmak2021}, we may disintegrate $\mu$ and $\nu$ with respect to the partition of $\mathbb{R}^n$ into the leaves of $Du$, see Theorem \ref{thm:disintegration}. We obtain Borel probability measures $\theta_1,\theta_2\in\mathcal{P}(CC(\mathbb{R}^n))$ and Borel probability measures
    \begin{equation*}
        \mu_{\mathcal{S}},\nu_{\mathcal{S}}\in\mathcal{P}_2(\mathbb{R}^n) \text{ defined for all }\mathcal{S}\in CC(\mathbb{R}^n)
    \end{equation*}    such that $\mu_{\mathcal{S}}$ is concentrated on $\mathcal{S}$  for $\theta_1$-almost every $\mathcal{S}$ and $\nu_{\mathcal{S}}$  is concentrated on $\mathcal{S}$  for $\theta_2$-almost every $\mathcal{S}$. Moreover
  \begin{equation*}
    \mu=\int_{CC(\mathbb{R}^n)}\mu_{\mathcal{S}}\, d\theta_1(\mathcal{S})\text{ and }     \nu=\int_{CC(\mathbb{R}^n)}\nu_{\mathcal{S}}\, d\theta_2(\mathcal{S}).
  \end{equation*}
On each element $\mathcal{S}$ of the partition, the derivative $Du$ is isometric and, as proven in Lemma \ref{lem:isoderma}, it can be written as a difference of two orthogonal projections $P_{V(\mathcal{S})_1}-P_{V(\mathcal{S})_2}$, where $V(\mathcal{S})_1\oplus V(\mathcal{S})_2=V(\mathcal{S})$ gives an orthogonal decomposition of the tangent space $V(\mathcal{S})$ to $\mathcal{S}$. More precisely,
\begin{equation*}
    Du(x)-Du(x_0)=P_{V(\mathcal{S})_1}(x-x_0)-P_{V(\mathcal{S})_2}(x-x_0)\text{ for }x,x_0\in\mathcal{S}.
\end{equation*}

Let now $\sigma\in\Sigma(\mu,\nu)$ be optimal for $\mathcal{J}(\mu,\nu)$.
As the first two marginals of $\sigma$ are absolutely continuous with respect to the Lebesgue measure, and the map $(x,y,z)\mapsto (\mathcal{S}(x),\mathcal{S}(y))$ is Borel measurable, by Theorem~\ref{thm:disintegration}\ref{i:disintegrationfibre} we can disintegrate $\sigma$ with respect to this map and obtain probabilities $\sigma_{\mathcal{S}_1,\mathcal{S}_2}\in\mathcal{P}((\mathbb{R}^n)^3)$, and $\theta\in\mathcal{P}((CC(\mathbb{R}^n))^2)$ such that  $\sigma_{\mathcal{S}_1,\mathcal{S}_2}$ is concentrated on $\mathcal{S}_1\times\mathcal{S}_2\times\mathbb{R}^n$ for $\theta$-almost every $(\mathcal{S}_1,\mathcal{S}_2)$ and that for any Borel $E\subset(\mathbb{R}^n)^3$ the map
\begin{equation*}
   (CC(\mathbb{R}^n))^2\ni (\mathcal{S}_1,\mathcal{S}_2)\mapsto \sigma_{\mathcal{S}_1,\mathcal{S}_2}(E)\in\mathbb{R}
\end{equation*}
is Borel measurable. Moreover
\begin{equation}\label{eqn:sigma}
    \sigma=\int_{(CC(\mathbb{R}^n))^2}\sigma_{\mathcal{S}_1,\mathcal{S}_2}\, d\theta(\mathcal{S}_1,\mathcal{S}_2).
\end{equation}
For $\theta$-almost every $(\mathcal{S}_1,\mathcal{S}_2)$ and $\sigma_{\mathcal{S}_1,\mathcal{S}_2}$-almost every $(x,y,z)$, Lemma~\ref{lem:three} applied to~\eqref{eqn:optimality} and uniqueness of the leaves at $x$ and $y$ give $z\in\mathcal{S}_1\cap\mathcal{S}_2$.
 Let us take a map $z$, as in Lemma \ref{lem:selection}.

 By $\pi$, let us denote the projection of $\sigma$ onto the first two copies of $\mathbb{R}^n$ and for all $(\mathcal{S}_1,\mathcal{S}_2)$ let us set  $\pi_{\mathcal{S}_1,\mathcal{S}_2}$ to be the projection of $\sigma_{\mathcal{S}_1,\mathcal{S}_2}$ onto the first two copies of $\mathbb{R}^n$. By (\ref{eqn:sigma}) it follows that
\begin{equation*}
    \pi(D)=\int_{(CC(\mathbb{R}^n)^2}\pi_{\mathcal{S}_1,\mathcal{S}_2}(D)\, d\theta(\mathcal{S}_1,\mathcal{S}_2)
\end{equation*}
for any Borel $D\subset\mathbb{R}^n\times\mathbb{R}^n$, and by Corollary \ref{col:twozet} we see that the optimal $\sigma$ is given by the formula
\begin{equation}\label{eqn:sigmaform}
    \sigma(C)=\int_{(CC(\mathbb{R}^n)^2}(R_{\mathcal{S}_1,\mathcal{S}_2})_{\#}\pi_{\mathcal{S}_1,\mathcal{S}_2}(C)\, d\theta(\mathcal{S}_1,\mathcal{S}_2),
\end{equation}
where
\begin{equation*}
    R_{\mathcal{S}_1,\mathcal{S}_2}\colon (\mathbb{R}^n)^2\to(\mathbb{R}^n)^3
\end{equation*}
is defined by the formula
\begin{equation*}
     R_{\mathcal{S}_1,\mathcal{S}_2}(x,y)=\Big(x,y,P_{V(\mathcal{S}_1)_2}(x-z(\mathcal{S}_1,\mathcal{S}_2))+P_{V(\mathcal{S}_2)_1}(y-z(\mathcal{S}_1,\mathcal{S}_2))+z(\mathcal{S}_1,\mathcal{S}_2)\Big),
\end{equation*}
with $\pi_{\mathcal{S}_1,\mathcal{S}_2}$ concentrated on pairs $(x,y)$ that satisfy \ref{i:vv}.

Let $\theta_{\mathcal{S}},\theta^{\mathcal{S}}$ be the probabilities on $CC(\mathbb{R}^n)$, defined for all $\mathcal{S}\in CC(\mathbb{R}^n)$, that are disintegration kernels of $\theta$ with respect to the first and with respect to the second coordinate respectively, obtained through applications of Theorem \ref{thm:disintegration}. Let us denote the pushforward measures by $\theta^1,\theta^2$ respectively, so that (\ref{eqn:thetas}) is satisfied with $\theta^i$ in place of $\theta_i$ for $i=1,2$.

Since $\sigma\in\Sigma(\mu,\nu)$, the marginals of $\pi$ are $\mu$ and $\nu$.
Therefore, if $C\subset CC(\mathbb{R}^n)$, $A\subset\mathbb{R}^n$ are Borel sets, then
\begin{align*}
& \int_{CC(\mathbb{R}^n)}\mathbf{1}_C(\mathcal{S}_1)\Bigg(\int_{CC(\mathbb{R}^n)}\mathrm{P}_1\pi_{\mathcal{S}_1,\mathcal{S}_2}(A)\, d\theta_{\mathcal{S}_1}(\mathcal{S}_2)\Bigg)\, d\theta^1(\mathcal{S}_1)= \\
&\int_{(CC(\mathbb{R}^n))^2}\mathrm{P}_1\pi_{\mathcal{S}_1,\mathcal{S}_2}(A\cap\mathcal{S}^{-1}(C))\,d\theta(\mathcal{S}_1,\mathcal{S}_2)=\pi((A\cap\mathcal{S}^{-1}(C))\times \mathbb{R}^n)=\\
&\mu(A\cap \mathcal{S}^{-1}(C))=\int_{CC(\mathbb{R}^n)}\mu_{\mathcal{S}}(A\cap\mathcal{S}^{-1}(C))\, d\theta_1(\mathcal{S})=\int_{CC(\mathbb{R}^n)}\mathbf{1}_C(\mathcal{S})\mu_{\mathcal{S}}(A)\, d\theta_1(\mathcal{S}).
\end{align*}
Taking $A=\mathbb{R}^n$ shows that $\theta^1=\theta_1$. Moreover, applying the identity with Borel $A\subset\mathbb{R}^n$, and arguing as in the proof of Theorem \ref{thm:parti}, we see that (\ref{eqn:mudis}) holds true.
The fact that $\theta^2=\theta_2$ and the identity (\ref{eqn:nudis}) are proved analogously, and this completes the proof of \ref{i:munu} and shows that $\theta$ is a coupling of $\theta_1$ and of $\theta_2$.

The condition \ref{i:coupling} is an immediate corollary of the fact that
\begin{equation*}
    \sigma\bigg(\Big(\bigcup\{(\mathcal{S}_1,\mathcal{S}_2)\in CC(\mathbb{R}^n))^2\mid \mathcal{S}_1\cap\mathcal{S}_2=\emptyset\}\times\mathbb{R}^n\Big)\bigg)=0.
\end{equation*}

We know also that $\sigma$ has to satisfy the martingale constraints.
These are equivalent to requiring that for $\theta_1$-almost every $\mathcal{S}_1$ and for $\theta_2$-almost every $\mathcal{S}_2$ and for all Borel, bounded $g,h\colon\mathbb{R}^n\to\mathbb{R}$ we have
\begin{equation*}
\int_{CC(\mathbb{R}^n)}\Bigg(\int_{(\mathbb{R}^n)^2}\Big(P_{V(\mathcal{S}_2)_1}(y-z(\mathcal{S}_1,\mathcal{S}_2))
-P_{V(\mathcal{S}_1)_1}(x-z(\mathcal{S}_1,\mathcal{S}_2))\Big)g(x)\, d\pi_{\mathcal{S}_1,\mathcal{S}_2}(x,y)\Bigg)\,d\theta_{\mathcal{S}_1}(\mathcal{S}_2)=0,
\end{equation*}
and
\begin{equation*}
\int_{CC(\mathbb{R}^n)} \Bigg( \int_{(\mathbb{R}^n)^2}\Big(P_{V(\mathcal{S}_1)_2}(x-z(\mathcal{S}_1,\mathcal{S}_2))
-P_{V(\mathcal{S}_2)_2}(y-z(\mathcal{S}_1,\mathcal{S}_2))\Big)h(y)\, d\pi_{\mathcal{S}_1,\mathcal{S}_2}(x,y)\Bigg)\,d\theta^{\mathcal{S}_2}(\mathcal{S}_1)=0.
\end{equation*}
By \ref{i:vv} we see that this is equivalent to
\begin{equation*}
P_{V(\mathcal{S}_1)_1}\Bigg( \int_{CC(\mathbb{R}^n)}\bigg(\int_{(\mathbb{R}^n)^2}\Big(P_{V(\mathcal{S}_2)_1}(y-z(\mathcal{S}_1,\mathcal{S}_2))
-(x-z(\mathcal{S}_1,\mathcal{S}_2))\Big)g(x)\, d\pi_{\mathcal{S}_1,\mathcal{S}_2}(x,y)\bigg)\,d\theta_{\mathcal{S}_1}(\mathcal{S}_2)\Bigg)=0,
\end{equation*}
and
\begin{equation*}
P_{V(\mathcal{S}_2)_2} \Bigg(\int_{CC(\mathbb{R}^n)} \bigg( \int_{(\mathbb{R}^n)^2}\Big(P_{V(\mathcal{S}_1)_2}(x-z(\mathcal{S}_1,\mathcal{S}_2))
-(y-z(\mathcal{S}_1,\mathcal{S}_2))\Big)h(y)\, d\pi_{\mathcal{S}_1,\mathcal{S}_2}(x,y)\bigg)\,d\theta^{\mathcal{S}_2}(\mathcal{S}_1)\Bigg)=0.
\end{equation*}
This is to say, \ref{i:bimartplan2} holds true.

Conversely, suppose we have Borel probability measures $\pi_{\mathcal{S}_1,\mathcal{S}_2}\in\mathcal{P}((\mathbb{R}^n)^2)$, defined for all $(\mathcal{S}_1,\mathcal{S}_2)\in (CC(\mathbb{R}^n))^2$,  a coupling $\theta$  of $\theta_1$ and $\theta_2$ that satisfy \ref{i:displan2}, \ref{i:coupling}, \ref{i:munu}, \ref{i:bimartplan2}, \ref{i:vv}.
Let $\sigma$ be defined as in (\ref{eqn:sigmadis2}). Then $\sigma\in\mathcal{P}((\mathbb{R}^n)^3)$ is a well-defined probability measure, thanks to Lemma \ref{lem:selection} and Lemma \ref{lem:vmeasur}. Moreover, thanks to \ref{i:munu}, its projection onto the first two coordinates is a coupling of $\mu$ and of $\nu$.

We now justify that the third coordinate prescribed by
$R_{\mathcal{S}_1,\mathcal{S}_2}$ belongs to both leaves. By
\ref{i:coupling}, \ref{i:munu}, and the
absolute continuity of $\mu,\nu$, for $\theta$-almost every
$(\mathcal{S}_1,\mathcal{S}_2)$ and $\pi_{\mathcal{S}_1,\mathcal{S}_2}$-almost every $(x,y)$, the points $x$
and $y$ belong to the unique leaves $\mathcal{S}_1$ and $\mathcal{S}_2$, respectively, with non-empty intersection. Fix such $(\mathcal{S}_1,\mathcal{S}_2,x,y)$ and write
\begin{align*}
 z_0&=z(\mathcal{S}_1,\mathcal{S}_2),\\
x'&=z_0+P_{V(\mathcal{S}_1)_2}(x-z_0),\\
y'&=z_0+P_{V(\mathcal{S}_2)_1}(y-z_0).
\end{align*}
Then $z'=x'+y'-z_0$ is the third coordinate of $R_{\mathcal{S}_1,\mathcal{S}_2}(x,y)$. The argument in the proof of Corollary \ref{col:subspaces} shows that   $Du$ is isometric on $\{y, x'\}$ and on $\{x, y'\}$.
By the uniqueness of the leaves containing $x,y$  we get
\begin{equation*}
  x'\in\mathcal{S}_2,\quad y'\in\mathcal{S}_1.
\end{equation*}
Since $x'-z_0$ belongs to $V(\mathcal{S}_1)_2$ -- the negative directions of the derivative on $\mathcal{S}_2$ -- also we must have $x'-z_0\in V(\mathcal{S}_2)_2$. Similarly, $y'-z_0\in V(\mathcal{S}_1)_1$.
Hence, for $i=1,2$, Corollary \ref{col:square} gives
\begin{equation*}
z'= x'+P_{V(\mathcal{S}_i)_1} (y'-x')\in\mathcal{S}_i.
\end{equation*}
Thus $ z'\in\mathcal{S}_1\cap\mathcal{S}_2$, which is the common
face by Lemma \ref{lem:intersection}. Condition \ref{i:vv} and the converse
implication of Corollary \ref{col:twozet}, applied with the reference triple
$(z_0,z_0,z_0)$, now imply that
(\ref{eqn:optimality}) holds $\sigma$-almost everywhere.

Finally, the preceding argument also gives, pointwise,
\begin{align*}
  z'-x
 &=P_{V(\mathcal{S}_1)_1}\Big(
   P_{V(\mathcal{S}_2)_1}(y-z_0)-(x-z_0)\Big),\\
  z'-y
 &=P_{V(\mathcal{S}_2)_2}\Big(
   P_{V(\mathcal{S}_1)_2}(x-z_0)-(y-z_0)\Big).
\end{align*}
Therefore \ref{i:bimartplan2}, followed by integration with respect to
$\theta_1$ and $\theta_2$, gives, for all bounded Borel
$g,h\colon\mathbb{R}^n\to\mathbb{R}$,
\begin{equation*}
 \int_{(\mathbb{R}^n)^3}(z-x)g(x)\,d\sigma(x,y,z)=0,
 \qquad
 \int_{(\mathbb{R}^n)^3}(z-y)h(y)\,d\sigma(x,y,z)=0.
\end{equation*}
Together with \ref{i:munu}, these are precisely the martingale constraints.
Thus $\sigma\in\Sigma(\mu,\nu)$ is optimal for
$\mathcal{J}(\mu,\nu)$.
\end{proof}

\begin{example}\label{exa:atomicthird}
Let us consider the one-dimensional case -- let
$\mu,\nu\in\mathcal{P}_2(\mathbb{R})$ be absolutely continuous. In that case, the leaves are closed intervals, possibly unbounded or reduced to a point. If a leaf containing $x$ is a singleton, then $z=x=y$ for
$\sigma$-almost every such $(x,y,z)$. We therefore remove the trivial diagonal component supported on singleton leaves and retain the notation $\pi=\mathrm{P}_{12}\sigma$ and $\theta$ for the
restrictions to non-trivial leaves.

For $i=1,2$, let
\begin{equation*}
    \mathcal{S}(\mathbb{R})_i=\bigl\{ \mathcal{S}\in\mathcal{S}(\mathbb{R}) \mid V(\mathcal{S})_i=\mathbb{R}\bigr\}.
\end{equation*}
Thus, if $\mathcal{S}\in\mathcal{S}(\mathbb{R})_1$, then
\begin{equation*}
   V(\mathcal{S})_1=\mathbb{R},\quad V(\mathcal{S})_2=\{0\},
\end{equation*}
whereas, if $\mathcal{S}\in\mathcal{S}(\mathbb{R})_2$, then
\begin{equation*}
    V(\mathcal{S})_1=\{0\},\quad V(\mathcal{S})_2=\mathbb{R}.
\end{equation*}
The set of non-trivial leaves is at most countable. Moreover, two distinct intersecting non-trivial leaves belong to different families: if $D^2u$ had the same sign on their union, then they would be contained in a common leaf, by maximality.

There is also a restriction on the direction in which distinct leaves can be coupled. Suppose that
\begin{equation*}
    \mathcal{S}_1\in\mathcal{S}(\mathbb{R})_2, \quad \mathcal{S}_2\in\mathcal{S}(\mathbb{R})_1, \quad \mathcal{S}_1\cap\mathcal{S}_2\neq\emptyset.
\end{equation*}
Then Theorem~\ref{thm:noassu}, \ref{i:vv}, gives $x=z(\mathcal{S}_1,\mathcal{S}_2)=y$ for $\pi_{\mathcal{S}_1,\mathcal{S}_2}$-almost every $(x,y)$.
Since $\mu,\nu$ are absolutely continuous, it follows that $\theta(\mathcal{S}_1,\mathcal{S}_2)=0$ for every such off-diagonal pair. Consequently, apart from the diagonal components, mass can only be transported from a leaf in $\mathcal{S}(\mathbb{R})_1$ to a neighbouring leaf in $\mathcal{S}(\mathbb{R})_2$. For the components that can have positive mass, the third coordinate in Theorem~\ref{thm:noassu} is therefore given by
\begin{equation*}
    z=\begin{cases}y,
&\mathcal{S}_1=\mathcal{S}_2
 \in\mathcal{S}(\mathbb{R})_1,\\
 x,
&\mathcal{S}_1=\mathcal{S}_2
 \in\mathcal{S}(\mathbb{R})_2,\\
z(\mathcal{S}_1,\mathcal{S}_2),
&\mathcal{S}_1\in\mathcal{S}(\mathbb{R})_1,\quad
 \mathcal{S}_2\in\mathcal{S}(\mathbb{R})_2,\quad
 \mathcal{S}_1\neq\mathcal{S}_2.
\end{cases}
\end{equation*}

We use the shorthand $\theta(\mathcal{S}_1,\mathcal{S}_2)=\theta\bigl(\{(\mathcal{S}_1,\mathcal{S}_2)\}\bigr)$.
Since the set of non-trivial leaves is at most countable,
\begin{equation*}
    \theta =\sum_{\substack{ \mathcal{S}_1,\mathcal{S}_2\in\mathcal{S}(\mathbb{R})\\
\mathcal{S}_1\cap\mathcal{S}_2\neq\emptyset}}\theta(\mathcal{S}_1,\mathcal{S}_2)\delta_{(\mathcal{S}_1,\mathcal{S}_2)}.
\end{equation*}
The fact that the marginals of \(\theta\) are
\(\theta_1,\theta_2\) shows that
\begin{equation*}
\mu(\mathcal{S})=\theta_1(\mathcal{S})=\sum_{\substack{\mathcal{S}_2\in\mathcal{S}(\mathbb{R})\\\mathcal{S}_2\cap\mathcal{S}\neq\emptyset}}\theta(\mathcal{S},\mathcal{S}_2),
\end{equation*}
and
\begin{equation*}
\nu(\mathcal{S})=\theta_2(\mathcal{S})=\sum_{\substack{ \mathcal{S}_1\in\mathcal{S}(\mathbb{R})\\ \mathcal{S}_1\cap\mathcal{S}\neq\emptyset}}\theta(\mathcal{S}_1,\mathcal{S}).
\end{equation*}
In view of the allowed directions described above,
\begin{align*}
&\mu(\mathcal{S})=\theta(\mathcal{S},\mathcal{S})+\sum_{\substack{ \mathcal{S}_2\in\mathcal{S}(\mathbb{R})_2\\ \mathcal{S}\cap\mathcal{S}_2\neq\emptyset}}\theta(\mathcal{S},\mathcal{S}_2),\quad\nu(\mathcal{S})=\theta(\mathcal{S},\mathcal{S})\text{, if }\mathcal{S}\in\mathcal{S}(\mathbb{R})_1,\\
&\mu(\mathcal{S})=\theta(\mathcal{S},\mathcal{S}),\quad \nu(\mathcal{S})=\theta(\mathcal{S},\mathcal{S})+\sum_{\substack{ \mathcal{S}_1\in\mathcal{S}(\mathbb{R})_1\\\mathcal{S}_1\cap\mathcal{S}\neq\emptyset}}\theta(\mathcal{S}_1,\mathcal{S}),
\text{ if }\mathcal{S}\in\mathcal{S}(\mathbb{R})_2.
\end{align*}

Let us now spell out the equalities \ref{i:bimartplan2}. We write $\theta_{\mathcal{S}_1}(\mathcal{S}_2)$ for $\theta_{\mathcal{S}_1}(\{\mathcal{S}_2\})$, and similarly for $\theta^{\mathcal{S}_2}(\mathcal{S}_1)$.

Fix $\mathcal{S}_1\in\mathcal{S}(\mathbb{R})_1$. For every bounded Borel function $g\colon\mathbb{R}\to\mathbb{R}$, the first equality in \ref{i:bimartplan2} becomes
\begin{equation*}
\theta_{\mathcal{S}_1}(\mathcal{S}_1)
\int_{\mathbb{R}^2}
(y-x)g(x)\,
d\pi_{\mathcal{S}_1,\mathcal{S}_1}(x,y)
=-\sum_{\substack{ \mathcal{S}_2\in\mathcal{S}(\mathbb{R})_2\\\mathcal{S}_1\cap\mathcal{S}_2\neq\emptyset}}\theta_{\mathcal{S}_1}(\mathcal{S}_2)\int_{\mathbb{R}^2}\bigl(z(\mathcal{S}_1,\mathcal{S}_2)-x\bigr)g(x)\, d\pi_{\mathcal{S}_1,\mathcal{S}_2}(x,y).
\end{equation*}
Indeed, on the diagonal component the third coordinate is \(y\),
whereas on an off-diagonal component it is
\(z(\mathcal{S}_1,\mathcal{S}_2)\).

Fix now $\mathcal{S}_2\in\mathcal{S}(\mathbb{R})_2$. For every bounded Borel function $h\colon\mathbb{R}\to\mathbb{R}$, the second equality in \ref{i:bimartplan2} becomes
\begin{equation*}
\theta^{\mathcal{S}_2}(\mathcal{S}_2)
\int_{\mathbb{R}^2}
(x-y)h(y)\,
d\pi_{\mathcal{S}_2,\mathcal{S}_2}(x,y)
=-\sum_{\substack{ \mathcal{S}_1\in\mathcal{S}(\mathbb{R})_1\\ \mathcal{S}_1\cap\mathcal{S}_2\neq\emptyset}}\theta^{\mathcal{S}_2}(\mathcal{S}_1)\int_{\mathbb{R}^2}\bigl(z(\mathcal{S}_1,\mathcal{S}_2)-y\bigr)h(y)\,d\pi_{\mathcal{S}_1,\mathcal{S}_2}(x,y).
\end{equation*}
For a source leaf in \(\mathcal{S}(\mathbb{R})_2\), the first equality
in \ref{i:bimartplan2} is void. Similarly, for a target leaf in
\(\mathcal{S}(\mathbb{R})_1\), the second equality is void.

Thus, \ref{i:bimartplan2} reduces to the usual bimartingale identity on a leaf if the corresponding row of $\theta$ or
column of $\theta$ has no off-diagonal mass.

For a concrete example, let $T$ be uniform on $[1,2]$ and independently choose one of the triples
\begin{equation*}
 (X,Y,Z)=(-2T,-T,-2T),\qquad (T,-T,0),\qquad (T,2T,2T)
\end{equation*}
with probability $1/3$ each. The endpoint laws are absolutely continuous. Averaging the last two alternatives at $X=T$ and the first two at $Y=-T$ gives $\mathbb{E}[Z\mid X]=X$ and $\mathbb{E}[Z\mid Y]=Y$, and elsewhere the corresponding equality is pointwise. Each triple satisfies (\ref{eqn:optimality}) for $u(t)=t|t|/2$, so the plan and $u$ are optimal. For any optimal plan with these marginals, Lemma~\ref{lem:three} with the same $u$ gives $|Z|-Z=|X|-X$ and $|Z|+Z=|Y|+Y$. Thus $\mathbb{P}(Z<0)=\mathbb{P}(X<0)=1/3$ and $\mathbb{P}(Z>0)=\mathbb{P}(Y>0)=1/3$. Every optimal third marginal has an atom of mass $1/3$ at the common endpoint of the two leaves.
\end{example}

\subsection*{Funding}
The author gratefully acknowledges the support of the University of Toronto and of the Beijing Institute of Mathematical Sciences and Applications. This project was partially funded by the Polish National Agency for Academic Exchange (NAWA) under the Polish Returns 2025 programme, agreement no.\ \mbox{BPN/PPO/2025/1/00055}, and by the National Natural Science Foundation of China, grant no.\ \mbox{W2633033}.
\subsection*{Competing interests}
None.
\subsection*{Data availability statement.}
No data was used for the research described in the article.

\providecommand{\bysame}{\leavevmode\hbox to3em{\hrulefill}\thinspace}
\providecommand{\MR}{\relax\ifhmode\unskip\space\fi MR }
\providecommand{\MRhref}[2]{%
  \href{http://www.ams.org/mathscinet-getitem?mr=#1}{#2}
}
\providecommand{\href}[2]{#2}
\bibliographystyle{amsplain}

\end{document}